\documentclass[12pt,reqno]{amsart}
\usepackage[margin=1in]{geometry}
\usepackage{amsmath,amssymb,mathrsfs}
\usepackage{amsthm}
\usepackage{color}
\usepackage{array}
\usepackage{tikz}
\usetikzlibrary{positioning,arrows.meta,decorations.pathreplacing}
\usepackage[colorlinks=true,citecolor=red,urlcolor=blue,linkcolor=blue]{hyperref}
\usepackage{float}
\newcommand{\beq}[1]{\begin{equation} \label{#1}}
\newcommand{\eeq}{\end{equation}}
\newcommand{\bea}{\bed\begin{array}{rl}}
\newcommand{\eea}{\end{array}\eed}
\newcommand{\bed}{\begin{displaymath}}
\newcommand{\eed}{\end{displaymath}}
\newcommand{\barray}{\begin{array}{ll}}
\newcommand{\earray}{\end{array}}
\newcommand{\disp}{\displaystyle}
\newcommand{\ad}{&\!\disp}
\newcommand{\aad}{&\!\disp}
\newcommand{\beqa}[1]{\begin{equation}\label{#1}\barray}
\newcommand{\eeqa}{\earray\end{equation}}
\newcommand{\al}{\alpha}
\newcommand{\e}{\varepsilon}

\newcommand{\la}{\lambda}
\newcommand{\La}{\Lambda}
\newcommand{\sg}{\sigma}

\newcommand{\ga}{\gamma}

\newcommand{\Ga}{\Gamma}
\newcommand{\dl}{\delta}
\newcommand{\Dl}{\Delta}
\newcommand{\cd}{(\cdot)}

\newcommand{\sqe}{\sqrt{\e}}

\def\indi{{\bf 1}}
\def\half{\frac{1}{2}}

\newcommand{\CA}{{\mathcal A}}
\newcommand{\CF}{{\mathcal F}}
\newcommand{\CM}{\mathcal{M}}
\newcommand{\CU}{\mathcal{U}}

\newcommand{\CO}{\mathcal{O}}
\newcommand{\CI}{\mathcal{I}}
\newcommand{\CB}{\mathcal{B}}
\newcommand{\CN}{\mathcal{N}}
\newcommand{\CL}{\mathcal{L}}
\newcommand{\CP}{\mathcal{P}}
\newcommand{\CE}{\mathcal{E}}
\newcommand{\CK}{\mathcal{K}}
\newcommand{\CD}{\mathcal{D}}
\newcommand{\CR}{\mathcal{R}}
\newcommand{\CS}{\mathcal{S}}

\newcommand{\CT}{\mathcal{T}}

\newcommand{\CW}{\mathcal{W}}
\newcommand{\CJ}{\mathcal{J}}
\newcommand{\CG}{\mathcal{G}}
\newcommand{\EE}{{\mathbb E}}
\newcommand{\PP}{{\mathbb P}}

\newcommand{\NN}{{\mathbb N}}
\newcommand{\rr}{{\mathbb R}}
\newcommand{\MM}{{\mathbb M}}

\newcommand{\UU}{{\mathbb U}}
\newcommand{\II}{{\mathbb I}}

\newcommand{\SK}{\mathscr{K}}

\newcommand{\SB}{\mathscr{B}}
\newcommand{\SE}{\mathscr{E}}
\newcommand{\ST}{\mathscr{T}}
\newcommand{\SD}{\mathscr{D}}

\newcommand{\SF}{\mathscr{F}}

\newcommand{\SO}{\mathscr{O}}

\newcommand{\fM}{\mathfrak{M}}
\newcommand{\fT}{\mathfrak{T}}

\newcommand{\fU}{\mathfrak{U}}

\newcommand{\fg}{\mathfrak{g}}
\newcommand{\Scv}{\mathsf{S}}
\newcommand{\notrow}[3]{\noindent
 \parbox[t]{0.19\textwidth}{\raggedright\strut #1\strut}\hfill
 \parbox[t]{0.52\textwidth}{\raggedright\strut #2\strut}\hfill
 \parbox[t]{0.23\textwidth}{\raggedright\strut #3\strut}\par}
\newcommand{\lbar}{\overline}
\newcommand{\wdt}{\widetilde}
\newcommand{\wdh}{\widehat}
\newcommand{\lf}{\lfloor}
\newcommand{\rf}{\rfloor}

\newcommand{\qv}[1]{\langle #1 \rangle}
\newcommand{\bqv}[1]{\big\langle #1 \big\rangle}
\newcommand{\Bqv}[1]{\Big\langle #1 \Big\rangle}

\newcommand{\feh}{\frac{1}{\sqrt{\e}h(\e)}}

\newcommand{\Ran}{\operatorname{Ran}}
\numberwithin{equation}{section}

\theoremstyle{plain}
\newtheorem{thm}{Theorem}[section]
\newtheorem{lem}[thm]{Lemma}
\newtheorem{cor}[thm]{Corollary}
\newtheorem{prop}[thm]{Proposition}
\theoremstyle{definition}

\newtheorem{ass}[thm]{Assumption}
\newtheorem{exm}[thm]{Example}
\theoremstyle{definition}
\newtheorem{rem}[thm]{Remark}

\newcommand{\thmref}[1]{Theorem~{\rm \ref{#1}}}
\newcommand{\lemref}[1]{Lemma~{\rm \ref{#1}}}

\newcommand{\corref}[1]{Corollary~{\rm \ref{#1}}}
\newcommand{\propref}[1]{Proposition~{\rm \ref{#1}}}

\newcommand{\remref}[1]{Remark~{\rm \ref{#1}}}

\newcommand{\assmref}[1]{Assumption~{\rm \ref{#1}}}
\newcommand{\secref}[1]{Section~{\rm \ref{#1}}}
\newcommand{\appref}[1]{Appendix~{\rm \ref{#1}}}

\newcommand{\mixrate}{\rho}
\newcommand{\mixer}{\phi_{\zeta}}

\newcommand{\Bpi}{B_{\pi}}
\newcommand{\Pcond}{\mathfrak{P}}

\begin{document}

\title[MDP for SEEs in fast Markovian environments]{Moderate deviations for semilinear stochastic evolution equations in fast Markovian environments}

\author[H. Qian]{Hongjiang Qian}
\address{Department of Mathematics, Wayne State University}
\email{hjqian@wayne.edu}

\author[G. Yin]{George Yin}
\address{Department of Mathematics, University of Connecticut}
\email{gyin@uconn.edu}


\begin{abstract}
This work establishes moderate deviation principles for semilinear stochastic evolution equations with coefficients depending on fast Markov environments. The environment is either an exponentially ergodic Markov process on a Polish space, independent of the noise, or a finite-state Markov chain whose transition rates depend on time and on a finite-dimensional observation of the state. The rate function is quadratic in the residual of the linearized dynamics. Its covariance is the sum of the averaged covariance of the Wiener noise and the asymptotic covariance generated by the environment, given by a Green--Kubo formula
in the independent case and by a Poisson equation in the switching case. The rate function is written with the Moore--Penrose inverse of the square root of this covariance and requires the residual to lie in its range, not merely in the closure. This allows degenerate covariances and nonclosed ranges. In the switching case, the dependence of the invariant distribution on the slow state adds an extra term to the linearized drift. Since the stochastic convolution is not a semimartingale, the proofs combine the regeneration split chain technique with Galerkin approximation. The assumptions are verified for stochastic reaction--diffusion equations, including space-time white noise in one spatial dimension.
\end{abstract}

\maketitle

\setcounter{tocdepth}{2}
\tableofcontents

\section{Introduction}\label{sec:intro}

In this paper, we consider the following semilinear
stochastic evolution equation in a fast-varying Markovian environment:
\beq{ab-srde}\left\{\barray
dX^\e(t)=[AX^\e(t)+F(X^\e(t),\zeta_{t/\e})]dt
       +\sqrt\e\,G(X^\e(t),\zeta_{t/\e})dW^{\mathsf{Q}}(t),\; \forall\, t\in [0,1], \\
X^\e(0)=x_0\in H,
\earray\right.
\eeq
where $A$ is the generator of a $C_0$-semigroup $\{e^{tA}\}_{t\geq 0}$ and $W^{\mathsf{Q}}$ is a cylindrical $\mathsf{Q}$-Wiener process. The drift and diffusion coefficients
are modulated either by an independent exponentially ergodic Markov
process $\zeta$ on a general Polish space $\SE$ or by a state-dependent regime-switching process $\al^\e$ with generator $Q(t,\Upsilon X^\e(t))/\e$ and state space $\CM=\{1,\dots,m_0\}$. Here, $\Upsilon:H\to\rr^{d_\Upsilon}$ is a bounded linear operator. Example~\ref{ex1} and Example~\ref{ex2}
near the end of this section illustrate the origins of these two environments. The process $X^\e$ evolves on the slow time scale and $\zeta_{\cdot/\e}$ or $\al^\e$ evolves on the fast time scale.

Stochastic reaction--diffusion equations (SRDEs)  provide a standard class of examples of \eqref{ab-srde}. They describe diffusive
systems in combustion, population dynamics, ecology, and chemical
transport. We refer to Da Prato and Zabczyk \cite{DZ92} and Cerrai
\cite{Cer01bk} for a comprehensive study. On a bounded domain $\SO\subset\rr^d$, $d\geq1$, consider
\beq{r-srde}\left\{\barray
\disp \frac{\partial X^\e}{\partial t}(t,\xi)=\CA X^\e(t,\xi)+f(\xi,X^\e(t,\xi),\zeta_{t/\e})+\sqrt{\e} g(\xi,X^\e(t,\xi),\zeta_{t/\e}) \frac{\partial^2 W}{\partial t\partial \xi}(t,\xi), \\
X^\e(0,\xi)=x_0(\xi),\;  \xi\in \SO,\quad
\CN X^\e(t,\xi)=0,\;  t\geq 0,\; \xi\in \partial \SO,
\earray\right.\eeq
where $\CA$ is a uniformly elliptic
second-order differential operator with regular real coefficients,
and $\CN$ specifies Dirichlet or Neumann boundary conditions. The
reaction $f:\mathscr{O}\times\rr\times\SE\to\rr$ is measurable and satisfies
the conditions specified in \secref{sec:rd-applications}, and $x_0\in H:=L^2(\SO;\rr)$.  The
driving Wiener noise $W$ is white in time, with spatial covariance
specified in Section~\ref{sec:rd-applications}. Space-time white noise is
allowed for $d=1$, while higher dimensions require spatial regularity. Here $\zeta$ is the independent Markov environment in \eqref{ab-srde}.

 As $\e\to 0$, $X^\e$ converges to $\bar X$, the solution of the following averaged equation:
\beq{avg-srde}\frac{\partial \bar X}{\partial t}(t,\xi)=\CA\bar X(t,\xi)+\bar f(\xi,\bar X(t,\xi)), \quad \bar X(0,\xi)=x_0(\xi),\xi\in \SO,
\eeq
where $\CN\bar X(t,\xi)=0$ on $\partial\SO$, $\pi$ is the invariant measure of $\zeta$, and $\bar f(\xi,r):=\int_\SE f(\xi,r,z)\pi(dz)$. This is the averaging principle in the small noise setting. Related averaging results for stochastic reaction--diffusion equations are given in \cite{Gao19,Cer09,Cer11,CF09,WR13,CL17}; see also \cite{Kha68,Kha80,Fre78,FW98,KY04,RX21} for the corresponding finite-dimensional results.

For illustration, let us concentrate on the independent Markov process $\zeta$. Define
\beqa{erg}\ad
\bar F(x)=\int_\SE F(x,z)\pi(dz).
\eeqa
The abstract averaged equation is
\beq{ab-avg}
d\bar X(t)=[A\bar X(t)+\bar F(\bar X(t))]dt,
\qquad \bar X(0)=x_0\in H.
\eeq

This paper aims to study the asymptotic behavior of the $H$-valued fluctuation
\begin{equation}\label{def:eta}
\eta^\e(t)=\frac{X^\e(t)-\bar X(t)}{\sqrt\e h(\e)},\quad t\in[0,1],
\end{equation}
as $\e \to 0$, where $h(\e)$ is the deviation scale satisfying
\beq{dev-s}
h(\e)\to +\infty,\quad \sqe h(\e)\to 0, \quad \text{ as } \e\to 0.
\eeq
Results concerning the decay rate of probabilities for $\eta^\e$ are called moderate deviation principles (MDPs). They bridge the gap between central limit theorems (CLTs) and large deviation principles (LDPs). If $h(\e)=1$, the limiting behavior of $\eta^\e$ falls in the domain of CLTs. For $h(\e)=\e^{-1/2}$, it leads to large deviation principles. Large and moderate deviation principles for stochastic evolution equations have been studied in \cite{CR04,Fre85,SBD19,Sow92,WZ15}. Results for slow--fast stochastic reaction--diffusion systems include \cite{WRD12,HSS19,GSS23}. In Euclidean space, Guillin \cite{Gui03} established an MDP for stochastic differential equations in such a random environment. We study the infinite-dimensional analogue on a
separable Hilbert space $H$ and also consider state-dependent regime-switching as our random environment.

The main results have two parts.
\begin{itemize}
\item[(i)] \thmref{thm:main} establishes an MDP in $C^0([0,1];H)$ for semilinear stochastic evolution equations whose drift and diffusion coefficients depend on an independent exponentially ergodic Markov process on a general Polish space.
\item[(ii)] \thmref{thm:feedback-main} establishes an MDP for finite-state regime switching with transition rates depending on a finite-dimensional observation of the slow component. The dependence of the invariant distribution on the slow state contributes an additional term to the linearized drift, given in \eqref{add-term}. The time-inhomogeneous Markov chain case follows as a corollary.
\end{itemize}
In both cases, the rate function is expressed through the combined covariance of the environmental and Wiener fluctuations, which can be degenerate. The characterization uses the Moore--Penrose inverse of its square root and an exact range condition, without requiring invertibility, closed range, or constant rank. The proofs use regeneration split chains and Poisson equations to deal with the fast environment. We also use dimension-uniform Galerkin estimates to transfer the finite-dimensional fluctuation principles to $H$. The approach is probabilistic and does not require a control representation of the rate function as in the weak convergence method \cite{DE11}.

By expanding $\eta^\e$ as in \eqref{exp-eta}, one finds that in addition to the linearization around $\bar X$ for the drift $F$, the following two integrals essentially contribute to the MDPs:
\beq{H-int}
\int_0^t e^{(t-s)A}
\bigl[F(\bar X(s),\zeta_{s/\e})-\bar F(\bar X(s))\bigr]\,ds
\quad\text{and}\quad
\int_0^t e^{(t-s)A}
G(\bar X(s),\zeta_{s/\e})\,dW^{\mathsf{Q}}(s).
\eeq
The first term represents fluctuations of the fast random environment. The second is a stochastic convolution driven by cylindrical Wiener noise with a diffusion coefficient depending on the environment. Both are $H$-valued convolution-type processes, and in both cases the semigroup $e^{(t-s)A}$ is part of the integrand. This is the main structural difference from the finite-dimensional setting considered in \cite{Gui03}. For the first term, define
\beq{F-tz}
F_c(s,z):=
F(\bar X(s),z)-\bar F(\bar X(s)),
\quad s\in[0,1],\ z\in\SE,
\eeq
where the subscript ``c'' in $F_c$ means ``centered''.
In \cite{Gui03}, the environmental contribution is described through the MDP of additive functionals of the fast-varying Markov process. Here the presence of the semigroup leads instead to an $H$-valued convolution of such a random environment, which yields a Volterra-type mapping.
Thus one must first understand how semigroup convolution changes the moderate deviation principles.
The second term in \eqref{H-int} creates a different difficulty. The term $
\int_0^t e^{(t-s)A}
G(\bar X(s),\zeta_{s/\e})\,dW^{\mathsf{Q}}(s)$, in general,
is not an $H$-valued semimartingale. Consequently, the finite-dimensional semimartingale arguments used in \cite{Gui03} cannot be applied directly. To overcome the difficulty, we first project the stochastic convolution onto the span of the first $n$ eigenfunctions of $A$. If $\Pi_n$ denotes this projection, the projected stochastic convolution is finite dimensional and its limiting behavior can be identified from its quadratic variation. The infinite-dimensional limit is then recovered by proving that the terms containing $I-\Pi_n$ are exponentially negligible, with estimates uniform in $n$.

The two random terms in \eqref{H-int} together determine the quadratic part of the rate function. For the environmental fluctuation, let $\{P_t\}_{t\geq0}$ be the transition semigroup of $\zeta$ and define
\beqa{def-U}\ad
U(t,z):=\int_0^\infty P_sF_c(t,\cdot)(z)\,ds.
\eeqa
For $x,y,u\in H$, let $
(x\otimes y)u:=\langle y,u\rangle_Hx$.
The asymptotic covariance generated by the fast-varying random environment is given by
\beqa{CF}\ad
C_F(t)
:=
\int_\SE
\bigl[
F_c(t,z)\otimes U(t,z)
+
U(t,z)\otimes F_c(t,z)
\bigr]
\,\pi(dz).
\eeqa
On the other hand, the stochastic convolution contributes an averaged covariance
\beqa{CG}\ad
\langle C_G(t)u,v\rangle_H
:=
\int_\SE\langle
\mathsf{Q}^* G(\bar X(t),z)^* u,\,
\mathsf{Q}^* G(\bar X(t),z)^* v\rangle_H
\,\pi(dz),
\qquad u,v\in H.
\eeqa
Hence the two fluctuation mechanisms combine through $
C(t):=C_F(t)+C_G(t)$.

Define $
\SD_A:=\{
\gamma\in W^{1,2}([0,1];H)\cap L^2([0,1];\CD(A)):
\gamma(0)=0\}$.
For $v\in\Ran(C(t)^{1/2})$,
$[C(t)^{1/2}]^\dagger v$ denotes the unique vector $
u\in[\ker C(t)^{1/2}]^\perp$
such that $C(t)^{1/2}u=v$. Equivalently,
$[C(t)^{1/2}]^\dagger v$ is the minimum-norm preimage of $v$ under $C(t)^{1/2}$. For $x,v\in H$, set
$
\overline{\nabla F}(x)v
:=
\int_\SE \nabla F(x,z)v\,\pi(dz).
$
Theorem~\ref{thm:main} gives the rate function in the explicit form
\beqa{intro-MP-rate}\ad
S(\gamma)
=
\frac12
\int_0^1
\big|
\bigl[C(t)^{1/2}\bigr]^\dagger
\bigl[
\dot\gamma(t)-A\gamma(t)
-\overline{\nabla F}(\bar X(t))\gamma(t)
\bigr]
\big|_H^2
\,dt,
\eeqa
provided that $\gamma\in\SD_A$, $
\dot\gamma(t)-A\gamma(t)
-\overline{\nabla F}(\bar X(t))\gamma(t)
\in\Ran(C(t)^{1/2})$
for almost every $t$, and the integral in \eqref{intro-MP-rate} is finite. Otherwise,
$S(\gamma)=+\infty$.
Note that the residual must belong to the actual range of $C(t)^{1/2}$, rather than only to its closure.

The state-dependent regime-switching case has the same two basic fluctuation mechanisms, but an additional term appears in the linearized drift. Indeed, if
$
\pi(t,\Upsilon X)=\{\pi_i(t,\Upsilon X)\}_{i\leq m_0} $
is the invariant distribution of $Q(t,\Upsilon X)$ and the  generator of the switching is $
\e^{-1}Q(t,\Upsilon X^\e(t))$,
then the dependence of $\pi_i(t,\Upsilon X)$ on the slow variable contributes
\beqa{add-term}\ad
\sum_{i=1}^{m_0}
\nabla_y\pi_i(t,\Upsilon\bar X(t))
[\Upsilon v]\,
F(\bar X(t),i),\quad v\in H.
\eeqa
Section~\ref{sec:feedback} gives the corresponding rate function. If $Q$ does not depend on $\Upsilon X^\e$, this additional term vanishes, and the result reduces to the case where the random environment is an independent time-inhomogeneous finite-state Markov chain.

The following examples describe how the two random environments arise in applications.
\begin{exm}[A fast exogenous environment]
\label{ex1}
Consider a population or chemical concentration field $X^\e(t,\xi)$ whose local reaction is influenced by an external reservoir. Suppose that the reservoir evolves independently of the field according to
 $
d\zeta_s=-\kappa(\zeta_s-\bar z)\,ds+\sigma\,dB_s^{\rm env}, \kappa,\sigma>0$,
with invariant distribution $N(\bar z,\sigma^2/(2\kappa))$. On the slow time scale, the field is driven by the accelerated process $\zeta^\e(t):=\zeta_{t/\e}$.
The environmental state may enter the reaction coefficient directly.  For example, $
 f(\xi,r,z)
 =
 f_0(\xi,r)+\psi(\xi)\tanh z$,
where $f_0$ describes the reaction in the absence of environmental fluctuations and the second term represents the effect of the reservoir. The Brownian motion $B^{\rm env}$ and the cylindrical Wiener process driving the slow field are assumed to be independent. This gives the fast exogenous Markov environment considered in Section~\ref{sec:erg-mc}.  Example~\ref{ex:OU-environment} provides an Ornstein--Uhlenbeck realization with explicitly computable environmental covariance.
\end{exm}

\begin{exm}[Fast ion-channel gating]
\label{ex2}
Let $X^\e(t,\xi)$ denote the membrane voltage along a spatially extended excitable cell, and let $\alpha^\e(t)\in\CM$ describe the configuration of a rapidly switching population of ion channels.
 A bounded-current version of the stochastic cable equation is
$
 \partial_t X^\e(t,\xi)
=
 \mathcal A X^\e(t,\xi)+I_{\rm app}(\xi)
 -g_{\rm ch}(\xi,\alpha^\e(t))
 \Theta\bigl(X^\e(t,\xi)-E_{\rm rev}(\alpha^\e(t))\bigr)
 +\sqrt\e\,
 G(X^\e(t),\alpha^\e(t))\mathsf{Q}\dot W_0(t,\xi)$.
Here $I_{\rm app}$ is the applied current, while  $g_{\rm ch}(\xi,i)$ and $E_{\rm rev}(i)$ denote the conductance and reversal potential associated with channel configuration $i$.  The function $\Theta$ agrees with the identity on the voltage range of interest and is smoothly truncated outside that range.
The channel configuration evolves as a finite-state Markov chain whose
transition rates depend on a finite-dimensional observation of the
voltage.
Here
$
 \Upsilon X
 = (
 \int_{\SO}\upsilon_k(\xi)X(\xi)\,d\xi)_{k=1}^{d_\Upsilon}
$
records a finite collection of spatially averaged voltages.  Thus the slow voltage affects the fast switching rates through
$\Upsilon X^\e$, while the channel state feeds back into the evolution
of $X^\e$. This gives the state-dependent regime-switching model studied
in Section~\ref{sec:feedback}. If the transition rates are independent of
$\Upsilon X^\e$, the switching process reduces to an exogenous fast
Markov environment. Related stochastic ion-channel models and their
infinite-dimensional averaging limits are studied in
\cite{Aus08,BR11,GT12}.
\end{exm}
The rest of the paper is organized as follows. Section~\ref{sec:assump} introduces the abstract
assumptions and the common analytic preliminaries, including the Galerkin approximation of the quadratic rate function
and its Moore--Penrose representation. Section~\ref{sec:erg-mc} establishes
the MDP for an independent exponentially ergodic Markov environment.
Section~\ref{sec:feedback} establishes the result for the state-dependent regime-switching process with a finite state space. The state-independent time-inhomogeneous Markov chain follows as \corref{cor:feedback-independent}.
Section~\ref{sec:rd-applications} collects the reaction--diffusion realizations of the abstract assumptions. Several technical proofs are postponed to the appendices.

\section{Functional setting and analytical  tools}\label{sec:assump}

\subsection{Notation and assumptions on the slow dynamics} Throughout the paper, we work on the time interval $[0,1]$. Thus, $T$ will not represent a terminal time. Let $H$ be a separable real Hilbert space with inner product
$\langle\cdot,\cdot\rangle_H$ and norm $|\cdot|_H$. We denote by
$\CL(H)$ and $\CL_2(H)$ the bounded and Hilbert--Schmidt operators
on $H$, respectively, with the Hilbert--Schmidt norm
$\|T\|_{\CL_2(H)}^2:=\operatorname{Tr}(T^*T)$ for $T\in \CL_2(H)$.
Let $W_0$ be a cylindrical Wiener process on $H$ with identity covariance.
For $\mathsf{Q}\in\CL(H)$, we use $W^{\mathsf{Q}}:=\mathsf{Q}W_0$ for the cylindrical noise with covariance $\mathsf{Q}\mathsf{Q}^*$. For a linear operator \(T\), we denote its kernel and range by \(\ker T\) and \(\Ran T\), respectively. $T^{\dagger}$
denotes its Moore--Penrose inverse, with domain
$\CD(T^{\dagger})=\Ran T\oplus(\Ran T)^\perp$.
If $v=T w+z$, where $w\in(\ker T)^\perp$ and $z\in(\Ran T)^\perp$, then
$T^{\dagger} v=w$. Let $C^0([0,1];H)$ be the
space of continuous paths starting at zero,
with sup-norm $\|\cdot\|_\infty:=\sup_{t\leq 1}|\cdot|_H$, and $D^0([0,1];H)$ be the set of all $H$-valued c\`adl\`ag paths with initial value $0$, equipped with the uniform norm. We use $W^{1,2}([0,1];H)$ for absolutely continuous $H$-valued paths whose weak derivative belongs to $L^2([0,1];H)$. The letter $K$, with and without subscripts, represents a generic constant whose value may change from line to line. We write $\NN_+=\{1,2,\dots\}$.

\begin{ass}\label{ass:A}\rm
The operator $A:\CD(A)\subset H\to H$ is negative, self-adjoint and generates a $C_0$-semigroup
$e^{tA}$ for $t\geq0$. Moreover, there exist an orthonormal basis $\{e_k\}_{k\geq1}$ of $H$
and numbers
$0<\alpha_1\leq\alpha_2\leq\cdots\uparrow\infty$ such that
$Ae_k=-\alpha_ke_k$ for $k\geq1$.
\end{ass}

For $\delta\geq0$, we define the fractional Sobolev space $\mathcal H_\delta:=\CD((-A)^\delta)$ and its norm
$ \|x\|_\delta^2:=\sum_{k\geq1}(1+\alpha_k^{2\delta})
|\langle x,e_k\rangle_H|^2$.
It then gives (see \cite{Paz83}) for any $0\leq\delta\leq1$,
\begin{align}
 \|(-A)^\delta e^{tA}\|_{\CL(H)}&\leq K_\delta t^{-\delta}, \label{frac-A}\\
 |(e^{tA}-I)x|_H&\leq K_\delta(t\wedge1)^\delta\|x\|_\delta. \label{e-I}
\end{align}
For $n\in\NN_+$, let $\Pi_n$ be the orthogonal projection onto
$H_n:=\operatorname{span}\{e_1,\ldots,e_n\}$, and set $A_n:=A|_{H_n}$.
We assume $F:H\times\SE\to H$ and $G:H\times\SE\to\CL(H)$ are strongly measurable.

\begin{ass}\label{ass:fg}\rm We also assume that there exist constants $K_0,\bar K_G,K_G,K_F',K_G'<\infty$ and
$0\leq\vartheta<1$ such that the following statements hold uniformly
in $z\in\SE$.
\begin{itemize}
\item[(H1)] The map $F(\cdot,z)$ is
$K_F'$-Lipschitz, uniformly in $z$, that is, $|F(x,z)-F(y,z)|_H\leq K_F'|x-y|_H$.
Moreover, $\sup_z|F(0,z)|_H\leq K_0$.
It is G\^ateaux differentiable with
$\sup_{x,z}\|\nabla F(x,z)\|_{\CL(H)}\leq K_F'$.
Moreover, the mapping $(x,v)\mapsto \nabla F(x,z)v$ is uniformly continuous on
every compact subset of $H\times H$, uniformly in $z\in\SE$.

\item[(H2)] For any $x,y\in H, z\in \SE$, and $0<r\leq 1$,
\begin{align}
 \|G(x,z)\mathsf{Q}\|_{\CL(H)} & \leq \bar K_G(1+|x|_H), \label{G-op}\\
 \|e^{rA} G(0,z)\mathsf{Q}\|_{\CL_2(H)}^2 & \leq K_G^2 r^{-\vartheta}, \label{G0}\\
 \|e^{rA}[G(x,z)-G(y,z)]\mathsf{Q}\|_{\CL_2(H)}^2
 &\leq (K_G')^2\, r^{-\vartheta}|x-y|_H^2.                         \label{G-diff}
\end{align}
\end{itemize}
\end{ass}
\begin{rem}\label{rem:G-bounds}\rm
By the triangle inequality and the bounds \eqref{G0} and \eqref{G-diff} in (H2), we have $\|e^{rA} G(x,z)\mathsf{Q}\|_{\CL_2(H)}^2 \leq (K_G+K_G'|x|_H)^2 r^{-\vartheta}$. Using $e^{rA}=e^{rA/2}e^{rA/2}$ and $\|(I-\Pi_n)e^{rA/2}\|_{\CL(H)} =e^{-\al_{n+1}r/2}$, we obtain
\beq{abstract-Gtail}
\|(I-\Pi_n)e^{rA}G(x,z)\mathsf{Q}\|_{\CL_2(H)}^2 \leq 2^{\vartheta}(K_G+K_G'|x|_H)^2 e^{-\al_{n+1}r}r^{-\vartheta}.
\eeq
For each $i$, the last bound \eqref{G-diff} in (H2) at $r=1$ and $e^{A}e_i=e^{-\al_i}e_i$ give $|\mathsf{Q}^*[G(x,z)-G(y,z)]^* e_i|_H \leq e^{\al_i} K_G'|x-y|_H$. On the ball $|x|_H\leq R$, \eqref{G0} and \eqref{G-diff} at $r=1$ give $|\mathsf{Q}^* G(x,z)^* e_i|_H \leq e^{\al_i}(K_G+K_G'R)$. Subtracting two inner products and applying Cauchy--Schwarz shows that each entry of $G(x,z)\mathsf{Q}\mathsf{Q}^*G(x,z)^*$ in the basis $\{e_k\}$, $x\mapsto \qv{\mathsf{Q}^* G(x,z)^* e_i, \mathsf{Q}^* G(x,z)^*e_j}_H$, is Lipschitz on this ball, uniformly in $z$, with constants $2e^{\al_i+\al_j}K_G'(K_G+K_G'R)$.
\end{rem}

\subsection{Preliminaries} We present some preliminary results that will be used in the sequel.
The following lemma concerns the exponential estimates for stochastic convolutions.

\begin{lem}\label{lem:conv-exp}
Let $\Phi$ be predictable. Suppose that there exists a deterministic constant
$K_\Phi<\infty$ such that for all $s$ and
$0<t\leq1$, $
\|e^{tA}\Phi(s)\|_{\CL_2(H)}^2
\leq K_\Phi^2\, t^{-\vartheta}$. Then the stochastic convolution $\int_0^{\cdot} e^{(\cdot-s)A} \Phi(s)dW_0(s)$  has a
continuous version in $H$, and
\begin{equation}\label{exp-conv}
 \PP\Big(\sup_{t\leq1}\Big|\int_0^te^{(t-s)A}\Phi(s)\,dW_0(s)\Big|_H> r \Big)
 \leq K_1\exp\{-r^2/(K_2 K_\Phi^2)\}.
\end{equation}
The constants $K_1$ and $K_2>0$ depend only on $\vartheta$
and the semigroup $e^{tA}$. In particular, they are independent of $\Phi$, $K_\Phi$, and $r$.
\end{lem}
The above lemma follows from Peszat \cite[Theorem~2.2]{Pes94} with $E=H$, exponent $a<(1-\vartheta)/2$ and bound $K_\Phi^2/(1-2a-\vartheta)$; see also \cite[Theorem~3.2]{CR04}. Let $C$ be a bounded nonnegative self-adjoint operator on $H$. Define the associated local rate function by
\bea\ad
L_C(v):=L(C,v)=\sup_{u\in H}\{\qv{v,u}-\tfrac12\langle Cu,u\rangle\}.
\eea
For a strongly measurable, uniformly bounded, and nonnegative self-adjoint field
$C:[0,1]\to \CL(H)$, set $L_C(t,v):=L_{C(t)}(v)$ and define
\beqa{quad-L}
 \CI_C(x)\ad :=\begin{cases}
\int_0^1L_C(t,\dot x(t)-Ax(t))\,dt,&\text{ if } x\in\SD_A,\\
 +\infty,&\text{otherwise}.
 \end{cases}
\eeqa
No inverse of $C(t)$ is assumed.   The next lemma identifies finite cost of $L_C(t,v)$ with membership in $\Ran(C(t)^{1/2})$.

\begin{lem}\label{lem:MP-energy} Let $\CU,H$ be Hilbert spaces, $T\in \CL(\CU,H)$ and $C=TT^*$. For $v\in H$,
\beq{LC}
L_C(v)=\sup_{y\in H}\{\qv{v,y}_H-\half \qv{Cy, y}_H\}=\begin{cases}
\half \min\{\|u\|_{\CU}^2: T u =v\}, \quad v\in \Ran T,\\
+\infty, \quad \text{otherwise},
\end{cases}
\eeq
and the minimum is attained at a unique $u\in\ker(T)^\perp$. In particular, for the field $C$ in \eqref{quad-L} and $T=C(t)^{1/2}$, we have
\beq{LC-t}
L_C(t,v)=\begin{cases}
      \half |[C(t)^{1/2}]^\dagger v|_H^2, \quad v\in \Ran(C(t)^{1/2}),\\
      +\infty, \quad \text{otherwise}.
\end{cases}
\eeq
Moreover, if $v(t)$ is strongly measurable such that $\int_0^1 L_C(t,v(t))dt<\infty$, then the minimum-norm preimage $u(t)=[C(t)^{1/2}]^\dagger v(t)$ is strongly measurable and belongs to $L^2([0,1];H)$.
\end{lem}

\begin{cor}\label{cor:add}
Let $C_1,C_2$ be bounded nonnegative self-adjoint operators on $H$ and $\CK\in \CL(H,H')$ for a Hilbert space $H'$. Then for $w\in H'$ and $v\in H$, we have
\begin{align}
      \inf\{L_{C_1}(v_1): \CK v_1=w\}&=L_{\CK C_1 \CK^*}(w) \label{LC1} \\
\inf_{v_1+v_2=v}\{L_{C_1}(v_1)+L_{C_2}(v_2)\}&=L_{C_1+C_2}(v), \label{LC1+C2}
\end{align}
and the infima are attained whenever they are finite.
\end{cor}
Proofs of \lemref{lem:MP-energy} and \corref{cor:add} are postponed to \appref{app:lem-MP-energy} and \ref{app:cor-add}.

\begin{lem}\label{lem:quad-lift}
Let $C:[0,1]\to \CL(H)$ be a strongly measurable, nonnegative self-adjoint field with uniform bounds. For $n\in \NN_+$, define $C_n(t):=\Pi_n C(t)\Pi_n|_{H_n}$ and $\CI_n$ on $C^0([0,1];H_n)$
by \eqref{quad-L} with $A_n,C_n$ in place of $A,C$. Then, for every $x\in C^0([0,1];H)$
\begin{equation}\label{projective-rate}
 \sup_{n\geq 1} \CI_n(\Pi_n x)=\CI_C(x),
\end{equation}
and $\CI_C$ is a good rate function on $C^0([0,1];H)$ with uniform topology.
Moreover, for every $R<\infty$, the set $\cup_{n\geq 1}\{z\in C^0([0,1];H_n): \CI_{n}(z)\leq R\}$ is relatively compact in $C^0([0,1];H)$.
\end{lem}

The proof is given in \appref{app:quad-lift}. The following corollary transfers LDPs from the projections $\Pi_nY^\e$ to $Y^\e$.

\begin{cor}[Galerkin approximation]\label{thm:abstract-lift}
Let $C,\CI_C$, and $\CI_n$ be as in \lemref{lem:quad-lift} and let $Y^\e$ take values in
$C^0([0,1];H)$. Suppose that
\begin{itemize}
\item[(i)] for
each $n\geq 1$, $\Pi_nY^\e$ satisfies an LDP in $C^0([0,1];H_n)$ with the uniform topology, speed $h^{-2}(\e)$, and good rate function $\CI_n$.
\item[(ii)] for
every $r>0$,
\begin{equation}\label{abstract-lift-tail}
 \lim_{n\to\infty}\limsup_{\e \to 0}h^{-2}(\e)
 \log\PP(\|(I-\Pi_n)Y^\e\|_\infty>r)=-\infty.
\end{equation}
\end{itemize}
Then $Y^\e$ satisfies an LDP in $C^0([0,1];H)$, with speed $h^{-2}(\e)$ and rate $\CI_C$ in \eqref{quad-L}.
\end{cor}

\begin{proof}By \eqref{abstract-lift-tail}, the laws of $\Pi_nY^\e$ are
exponentially good approximations of the laws of $Y^\e$ in the
uniform metric. Since $C^0([0,1];H)$ is Polish and the projected rate
functions are good, \cite[Theorem~4.2.16 (a) and
Exercise~4.2.29 (a)]{DZ09} yield an LDP for $Y^\e$ in $C^0([0,1];H)$
with a good rate function $\CJ$. By the contraction principle and uniqueness of the rate function
\cite[Theorem~4.2.1 and Lemma~4.1.4]{DZ09}, for every $x\in C^0([0,1];H)$
and $n\ge1$, $
\CI_n(\Pi_n x)
=\inf\bigl\{\CJ(y):y\in C^0([0,1];H),\ \Pi_n y=\Pi_n x\bigr\}.$
For fixed $x$, the constraint sets on the right are closed and
decrease to $\{x\}$. Since $\CJ$ is good,
\cite[Lemma~4.1.6 (a)]{DZ09} gives
$\CJ(x)=\sup_{n\ge1}\CI_n(\Pi_n x)$, which equals $\CI_C(x)$
by \eqref{projective-rate}.
\end{proof}

For $n\in\NN_+$ and $z\in C^0([0,1];H_n)$, define
\beq{def-Phi-n}
(\Phi_nz)(t):=z(t)+\int_0^tA_ne^{(t-s)A_n}z(s)\,ds,\qquad t\in[0,1].
\eeq
 Its inverse is given by $(\Phi_n^{-1}x)(t)=x(t)-\int_0^tA_nx(s)\,ds$. Since $A_n$ is bounded, $\Phi_n$ and $\Phi_n^{-1}$ are continuous in the uniform norm. Moreover, $\Phi_n$ maps $\{x\in W^{1,2}([0,1];H_n):x(0)=0\}$ bijectively onto itself, and $\frac{d}{dt}(\Phi_n^{-1}x)(t)=\dot x(t)-A_nx(t)$ on this set.

\section{Exponentially ergodic Markov process on a general Polish space}\label{sec:erg-mc}
\begin{figure}[t]
\centering
\begin{tikzpicture}[scale=0.85, transform shape,
 box/.style={draw,rounded corners=2pt,align=center,font=\footnotesize,inner sep=3pt},
 mod/.style={box,dashed},
 arr/.style={-{Stealth[length=4pt]},semithick},
 lab/.style={font=\scriptsize,fill=white,inner sep=1pt,align=center},
 node distance=6mm and 3mm]
\node[box] (eta) {$\eta^\e$\\ \thmref{thm:main}};
\node[box,below=13mm of eta] (heta) {$\wdh\eta^\e=\CT(R^\e)$\\ \propref{prop:cp}};
\node[mod,below=of heta] (R) {$R^\e=\la^\e+M^\e$\\ \lemref{lem:local-R}};
\node[mod,below left=7mm and 18mm of R] (la) {$\la^\e$\\ \thmref{thm:la}};
\node[box,below right=7mm and 30mm of R] (M) {$M^\e$\\ \thmref{thm:mdp-comb}};
\node[box,below=9mm of la] (D) {$D_\e^N$\\ \lemref{lem:neg-D}};
\node[box,left=of D] (L) {$L_\e^N$\\ \lemref{lem:LDP-L}};
\node[box,right=of D] (A) {$A_\e^N$\\ \lemref{lem:neg-A}};
\node[box,below left=9mm and 1.5mm of M.south] (PM) {$\Pi_nM^\e$\\ \lemref{new:finite}};
\node[box,below right=9mm and 1.5mm of M.south] (TM) {$(I-\Pi_n)M^\e$\\ \lemref{new:tail}};
\draw[arr] (heta) -- node[lab,draw,dashed,right=2mm]{exponential equivalence\\ \propref{prop:eta-equivalence}} (eta);
\draw[arr] (R) -- node[lab,right=2mm]{contraction} (heta);
\draw[arr] (la) -- (R);
\draw[arr] (M) -- (R);
\draw[arr] (L) -- node[lab,left]{LDP} (la);
\draw[arr] (D) -- node[lab]{negligible} (la);
\draw[arr] (A) -- (la);
\draw[arr] (PM) -- node[lab,left]{\corref{thm:abstract-lift}} (M);
\draw[arr] (TM) -- node[lab,right]{tail} (M);
\end{tikzpicture}
\caption{Scheme of the proof of \thmref{thm:main}. The arrows point from the ingredients to the conclusion. The dashed boxes are the steps modified in Section~\ref{sec:feedback}: $\la^\e$ is rewritten through a Poisson equation (Section~\ref{sec:sd-decomp}), the LDP of $R^\e$ is replaced by the joint LDP in \propref{prop:sd-driving}, and the exponential equivalence by \propref{prop:equi-eta-al}.}
\label{fig:roadmap}
\end{figure}
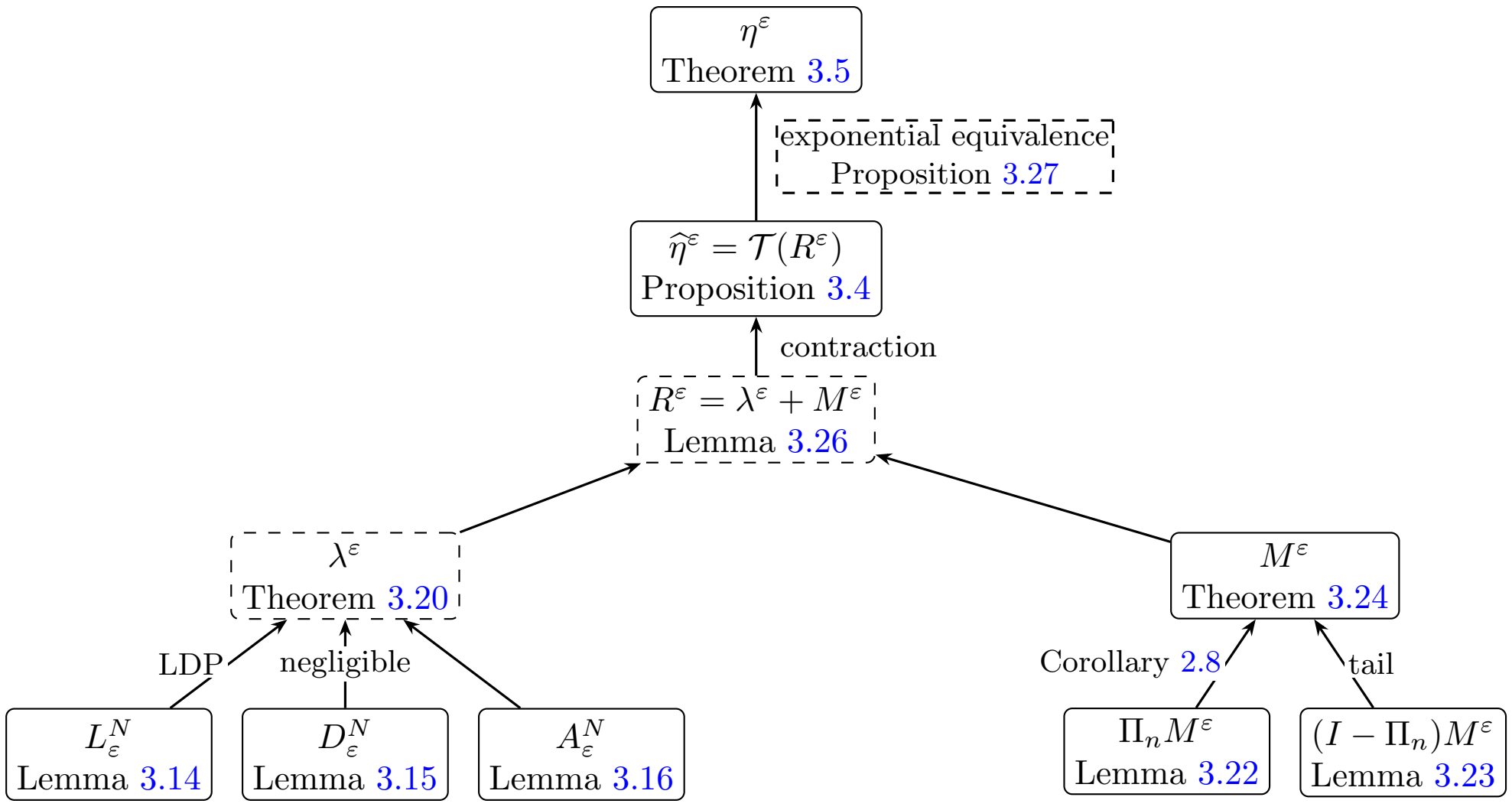

In this section, we consider the random environment to be an independent exponentially ergodic Markov process on a general Polish space. Figure~\ref{fig:roadmap} summarizes the proof architecture of \thmref{thm:main}. The argument separates the environmental and stochastic convolution fluctuations, combines their LDPs, and then obtains the main result  by exponential equivalence.

\subsection{Assumptions on the random environment}\label{ass:zeta}
Let $\zeta=\{\zeta_t\}_{t\geq0}$ be a c\`adl\`ag Markov process
on a general Polish space $\SE$, defined on
$(\Omega^\zeta,\CF^\zeta,(\CF_t^\zeta),\PP^\zeta)$, with transition
semigroup $\{P_t\}_{t\geq0}$ and invariant probability measure $\pi$. We assume that $\zeta$ is $\pi$-irreducible, aperiodic, and exponentially ergodic: there exist a finite $\phi_\zeta \in L^1(\pi)$  and a constant $0<\mixrate<1$ such that
\beq{exp-erg}
\|P_t(z,\cdot)-\pi\|_{\text{var}} \leq \mixer(z)\mixrate^t,\quad \forall\ t \geq 0.
\eeq
Denote by $\mu$ the initial law of $\zeta_0$. We also  assume that
\beqa{cond-mu}
K_\phi:=\int_\SE \mixer(z)\,\mu(dz)<\infty.
\eeqa
For details on exponentially ergodic Markov processes, we refer to  \cite{Num84,DMT95} and references therein. We can then take the stochastic basis $(\Omega,\CF,\{\CF_t\}_{t\geq 0},\PP)$ as the product probability space of the random environment $\zeta$ and the cylindrical Wiener process $W_0$.

Let $\SF$ be a separable Hilbert space and let $\ST:\SE\to\SF$ be a bounded map with $\pi(\ST)=0$. Define its \emph{Green--Kubo covariance} $C_{\ST}\in\CL(\SF)$ by
\beqa{CF-GK}
 \langle C_\ST u, v\rangle_\SF \ad =\int_0^\infty\int_\SE
 \big[\langle \ST(z),u\rangle_\SF P_s\langle \ST(\cdot),v\rangle_\SF(z) \\
 \aad\qquad\qquad + \langle \ST(z),v\rangle_\SF P_s\langle \ST(\cdot),u\rangle_\SF(z)\big]\pi(dz)ds, \quad u,v\in \SF.
 \eeqa
 The integrals are absolutely convergent by the exponential mixing of $\zeta$. To verify the nonnegativity, fix $u\in\mathscr F$ and set $f_u(z):=\langle\mathscr T(z),u\rangle_{\mathscr F}$. Since $\pi(f_u)=0$ and $\zeta$ is stationary and Markov under $\pi$, $\mathbb E_\pi[f_u(\zeta_s)f_u(\zeta_t)]=\pi(f_uP_{t-s}f_u)$ for $s\le t$, hence
\[
R^{-1}\mathbb E_\pi\Big|\int_0^R f_u(\zeta_s)\,ds\Big|^2
=\frac2R\int_0^R\!\!\int_0^t \pi(f_uP_{t-s}f_u)\,ds\,dt
=2\int_0^R\Big(1-\frac sR\Big)\pi(f_uP_sf_u)\,ds .
\]
The exponential mixing estimate gives an integrable bound on $\pi(f_uP_sf_u)$, so dominated convergence yields
$
\langle C_{\mathscr T}u,u\rangle_{\mathscr F}
=\lim_{R\to\infty}R^{-1}\,\mathbb E_\pi|\int_0^R\langle\mathscr T(\zeta_s),u\rangle_{\mathscr F}\,ds|^2\ge0 .
$
Hence $C_{\mathscr T}$ is a bounded nonnegative self-adjoint operator. The representation \eqref{CF-GK} is the Green--Kubo formula for the covariance of the additive functional generated by $\ST$; see, for example, \cite{Pav10} for its relation to the associated Poisson equation.

Define the averaged drift $\bar F$ by \eqref{erg}. By Assumption~\ref{ass:fg}, it is well-defined and is globally
Lipschitz with constant $K_F'$.

\begin{thm}\label{thm:XbarX}
Suppose that Assumptions~\ref{ass:A} and \ref{ass:fg} hold. We have
\begin{itemize}
\item[(1)]
for each $\e>0$ and $x_0\in H$, equation
\eqref{ab-srde} has a unique mild solution $X^\e(\cdot)\in L^2(\Omega;C([0,1];H))$ satisfying the following
form:
\bea\ad
X^\e(t)=e^{tA}x_0+\int_0^t e^{(t-s)A}F(X^\e(s),\zeta_{s/\e})ds+\sqrt{\e} \int_0^t e^{(t-s)A}G(X^\e(s),\zeta_{s/\e})dW^{\mathsf{Q}}(s).
\eea
\item[(2)]
for any $x_0 \in H$, equation \eqref{ab-avg} has a unique mild solution $\bar X\in C([0,1];H)$ such that
\bea\ad
\bar X(t)= e^{tA} x_0 +\int_0^t e^{(t-s)A} \bar F(\bar X(s))ds.
\eea
Moreover, $\|\bar X\|_\infty \leq K(1+|x_0|_H)<\infty$.
If $x_0\in \CD((-A)^{\theta'}), \theta' \in [0,1)$, we obtain
\bea
\disp \sup_{t\in [0,1]}\|\bar X(t)\|_{\theta'}
\ad \leq K_{\theta'}(1+\|x_0\|_{\theta'})<\infty,\\
|\bar X(t)-\bar X(s)|_H \ad \leq K_{\theta'}\, (t-s)^{\theta'} (1+\|x_0\|_{\theta'}), \quad\forall\, 0\leq s < t \leq 1.
\eea
\item[(3)]
there exist $K<\infty$ and $\kappa>0$ such that
\begin{equation}\label{solution-exp}
 \limsup_{\e\to0}\e\log\PP(\|X^\e\|_\infty> K)\leq-\kappa.
\end{equation}
The same statement applies when the environment has finite state space $\CM$ in Section~\ref{sec:feedback}.
\end{itemize}
\end{thm}

\begin{proof}
Statements (1) and (2) follow from standard fixed-point arguments for mild solutions; see \cite{DZ92} and Cerrai and R\"ockner \cite[Theorem~2.2]{CR04}. Thus, the details are omitted. For (3), let $R>0$, define the stopping time $\tau^\e:=\inf\{t\leq1:|X^\e(t)|_H>R\}\wedge1$ and $\CW_R^\e(t):=\sqe\int_0^te^{(t-s)A}\indi_{\{s\leq\tau^\e\}}G(X^\e(s),\zeta_{s/\e})dW^{\mathsf{Q}}(s)$. By \remref{rem:G-bounds} and \lemref{lem:conv-exp}, $\PP(\|\CW_R^\e\|_\infty>r)\leq K_1\exp\{-r^2/(\e K_2(K_G+K_G' R)^2)\}$ for $r>0$. On $\{\|\CW_R^\e\|_\infty\leq r\}$, (H1) and Gronwall's inequality give $\sup_{t\leq\tau^\e}|X^\e(t)|_H\leq(|x_0|_H+K_0+r)e^{K_F'}$. Take $r=e^{-K_F'}R/2$ and $R>2(|x_0|_H+K_0)e^{K_F'}$. Then $\sup_{t\leq\tau^\e}|X^\e(t)|_H<R$, hence $\tau^\e=1$, and \eqref{solution-exp} holds with $K=R$ and $\kappa=e^{-2K_F'}R^2/(4K_2(K_G+K_G' R)^2)$.
\end{proof}

\begin{lem}\label{lem:Fc-U}
Suppose that Assumptions~\ref{ass:A} and~\ref{ass:fg} hold and
$x_0\in\CD((-A)^\theta)$ for some $\theta\in(1/2,1)$.
Then $F_c$ defined in \eqref{F-tz} is centered under $\pi$ and satisfies
\beqa{Fc-H-bounds}\ad
 \sup_{t\in[0,1],\,z\in\SE}|F_c(t,z)|_H\le K,\qquad
 \sup_{z\in\SE}|F_c(t,z)-F_c(s,z)|_H\le K|t-s|^\theta.
\eeqa
The Bochner integral $U(t,z)$ in \eqref{def-U}
converges in $H$, and $\pi(U(t,\cdot))=0$,
\beqa{U-H-bounds}
 |U(t,z)|_H\le K\mixer(z),\qquad
 |U(t,z)-U(s,z)|_H\le K\mixer(z)|t-s|^\theta.
\eeqa
Consequently, $C_F(t)$ is the Green--Kubo covariance \eqref{CF-GK} of $\ST=F_c(t,\cdot)$ with $\SF=H$. It is a bounded, nonnegative,
self-adjoint operator and
\beq{CF-time-holder}
 \|C_F(t)-C_F(s)\|_{\CL(H)}\le K|t-s|^\theta.
\eeq
\end{lem}

\begin{proof}
By the definition of $F_c$, we have $\pi(F_c(t,\cdot))=\pi(F(\bar X(t),\cdot))-\bar F(\bar X(t))=0$. Moreover, by Theorem~\ref{thm:XbarX} and Assumption~\ref{ass:fg} (H1), we obtain
\bea
 \ad |F_c(t,z)|_H \le|F(\bar X(t),z)|_H
 +\int_\SE|F(\bar X(t),y)|_H\,\pi(dy) \le 2K(1+|\bar X(t)|_H)\le K,\\
\ad |F_c(t,z)-F_c(s,z)|_H
\le 2 K_F'|\bar X(t)-\bar X(s)|_H
 \le K|t-s|^\theta.
\eea
For $r\ge0$, \eqref{exp-erg} and the centering of $F_c(t,\cdot)$ give
\beqa{U-P-bound}
 |P_rF_c(t,\cdot)(z)|_H
 \le\int_\SE|F_c(t,y)|_H
 |P_r(z,dy)-\pi(dy)|
 \le K\mixer(z)\mixrate^r.
\eeqa
By replacing $F_c(t,\cdot)$ by $F_c(t,\cdot)-F_c(s,\cdot)$, we obtain
$
 |P_r(F_c(t,\cdot)-F_c(s,\cdot))(z)|_H
 \le K\mixer(z)\mixrate^r|t-s|^\theta.
$
Since $0<\mixrate<1$, $|U(t,z)|_H
 \le\int_0^\infty K\mixer(z)\mixrate^rdr
 =\frac{K\mixer(z)}{-\log\mixrate}$ and $|U(t,z)-U(s,z)|_H \le\frac{K\mixer(z)}{-\log\mixrate}|t-s|^\theta$. Moreover, we have $\pi(U(t,\cdot))
=\int_0^\infty\pi(P_rF_c(t,\cdot))dr=0$.
By \eqref{def-U} and Fubini's theorem, which applies by \eqref{U-P-bound} and $\pi(\mixer)<\infty$, for $u,v\in H$,
\bea
\qv{C_F(t)u,v}_H\ad=\int_0^\infty\!\!\int_\SE\big[\qv{F_c(t,z),u}_HP_s\qv{F_c(t,\cdot),v}_H(z)\\
\aad\qquad\qquad+\qv{F_c(t,z),v}_HP_s\qv{F_c(t,\cdot),u}_H(z)\big]\pi(dz)ds.
\eea
This is \eqref{CF-GK} with $\ST=F_c(t,\cdot)$ and $\SF=H$, so $C_F(t)$ is bounded, nonnegative and self-adjoint.
Finally,
$\|u\otimes v\|_{\CL(H)}\le|u|_H|v|_H$ yields
\begin{align*}
 &\|C_F(t)-C_F(s)\|_{\CL(H)}
 \le\int_\SE\Bigl(
 |F_c(t,z)-F_c(s,z)|_H
 [|U(t,z)|_H+|U(s,z)|_H]\\
 &\; +[|F_c(t,z)|_H+|F_c(s,z)|_H]
 |U(t,z)-U(s,z)|_H\Bigr)\pi(dz) \le K\pi(\mixer)|t-s|^\theta.
\end{align*}
\end{proof}

\subsection{\texorpdfstring{Decomposition of $\eta^\e$, main result, and auxiliary lemmas}{Decomposition of eta epsilon, main result, and auxiliary lemmas}}\label{sec:decomp-eta}
\subsubsection{\texorpdfstring{Decomposition of $\eta^\e$}{Decomposition of eta epsilon}} For any $t\in[0,1]$, we have
\beqa{exp-eta}\disp
\eta^\e(t) =\frac{X^\e(t)-\bar X(t)}{\sqrt{\e}h(\e)}
\ad = \frac{1}{\sqrt{\e}h(\e)}\int_0^t e^{(t-s)A}\big[ F(X^\e(s),\zeta_{s/\e})-\bar F(\bar X(s)) \big]ds \\
\aad\quad + h^{-1}(\e) \int_0^t e^{(t-s)A} G(X^\e(s),\zeta_{s/\e})dW^{\mathsf{Q}}(s).
\eeqa
The goal is to prove the LDP for $\eta^\e$ in $C^0([0,1];H)$.
We introduce an auxiliary process $\wdh \eta^\e$ as
\beqa{hat-eta}\disp
\wdh\eta^\e(t)
\ad:=\la^\e(t)+M^\e(t)+\int_0^t e^{(t-s)A} \overline{\nabla F}(\bar X(s))\wdh \eta^\e(s)ds,
\eeqa
where $\la^\e$ and $M^\e$ are defined as, for any $t\in[0,1]$,
\begin{align}
\la^\e(t)&:=\frac{1}{\sqrt{\e}h(\e)}\int_0^t e^{(t-s)A}\big[F_c(s,\zeta_{s/\e})\big]ds, \label{la} \\
M^\e(t)&:=h^{-1}(\e)\int_0^t e^{(t-s)A} G(\bar X(s),\zeta_{s/\e})dW^{\mathsf{Q}}(s) \label{M}
\end{align}
Here $F_c$ is defined in \eqref{F-tz}.  For our purpose, we prove the LDP of $\wdh\eta^\e$ and then establish the
exponential equivalence between $\eta^\e$ and $\wdh \eta^\e$, that is,
\beq{eta-equiv}
\lim_{\e\to 0} h^{-2}(\e) \log\PP(\|\eta^\e- \wdh \eta^\e\|_\infty >r) =-\infty, \quad \forall\,  r>0.
\eeq
 This is \propref{prop:eta-equivalence} in
Section~\ref{sec:equiv-etah}. Let
\beq{def-R}
R^\e(t):=\la^\e(t)+M^\e(t).
\eeq
Then $\wdh \eta^\e=\CT(R^\e)$, where $\CT:C^0([0,1];H)\to C^0([0,1];H)$ is the Volterra solution map:
\beq{Tu}
 \CT(x)(t)=x(t)+\int_0^te^{(t-s)A}\overline{\nabla F}(\bar X(s))\CT(x)(s)ds.
\eeq

\begin{lem}\label{lem:Tc}
The map $\CT$ is a continuous linear homeomorphism, and $
 \|\CT(x)-\CT(y)\|_\infty\leq e^{K_F'}\|x-y\|_\infty$.
Its inverse is $\CT^{-1}(y)=y-\int_0^\cdot e^{(\cdot-s)A}\overline{\nabla F}(\bar X(s))y(s)ds$.
\end{lem}

\begin{proof}
Let $(J y)(t):=\int_0^t e^{(t-s)A}\overline{\nabla F}(\bar X(s))y(s)ds$. Then $\|J^n\|\leq (K_F')^n/n!$ on $C^0([0,1];H)$, so the Neumann series gives $\CT=(I-J)^{-1}=\sum_{n\geq0} J^n$ and $\|\CT\|\leq e^{K_F'}$.
\end{proof}

Lemma~\ref{lem:Tc} and the contraction principle in \cite{DZ09} now give the following proposition.

\begin{prop}\label{prop:cp}
If $R^\e$ satisfies an LDP in $C^0([0,1];H)$ with speed $h^{-2}(\e)$ and a good rate function $\CI$, then
$\widehat\eta^\e=\CT(R^\e)$ satisfies an LDP with speed $h^{-2}(\e)$ and a good rate function $S$ given by
$
 S(x)=\CI(x-\int_0^\cdot e^{(\cdot-s)A}\overline{\nabla F}(\bar X(s))x(s)ds).
$
\end{prop}

The rate function of the first main result involves the averaged linearization $\overline{\nabla F}$, the environmental covariance $C_F$ in \eqref{CF}, and the Wiener covariance $C_G$ in \eqref{CG}. By \lemref{lem:Fc-U}, $C_F(t)$ is a bounded, nonnegative, self-adjoint operator and $C_F$ is H\"older continuous in time with exponent larger than $1/2$. By \eqref{G-op} and \thmref{thm:XbarX} (2), $C_G(t)$ is nonnegative and self-adjoint with $\sup_{t\in[0,1]}\|C_G(t)\|_{\CL(H)}\leq \bar K_G^2(1+\|\bar X\|_\infty)^2$. The first main result is as follows.

\begin{thm}\label{thm:main}
Suppose that Assumptions~\ref{ass:A} and~\ref{ass:fg}, and the
environmental conditions of Section~\ref{ass:zeta}, including \eqref{cond-mu}, hold.
Let $x_0\in\CD((-A)^\theta)$ be deterministic for some
$\theta\in(1/2,1)$.
Let $h(\e)$ satisfy \eqref{dev-s}. Then, $\eta^\e=(X^\e-\bar X)/(\sqrt{\e}h(\e))$
satisfies an LDP in $C^0([0,1];H)$ with speed $h^{-2}(\e)$
and good rate function given by
\begin{equation}\label{main-rate}
 S(x)=\begin{cases}
\disp \int_0^1L_{C_F+C_G}
 (t,\dot x(t)-Ax(t)-\overline{\nabla F}(\bar X(t))x(t))dt,&\text{ if }x\in \SD_A,\\
 +\infty,&\text{otherwise}.
 \end{cases}
\end{equation}
Equivalently, $X^\e$ satisfies the corresponding moderate deviation
principle. Let $C:=C_F+C_G$ and
$v_x(t):=\dot x(t)-Ax(t)-\overline{\nabla F}(\bar X(t))x(t)$. Then $S(x)$ has the following explicit form
\beqa{main-MP-rate}\ad
S(x)=\frac12\int_0^1
\big|\big[C(t)^{1/2}\big]^{\dagger}v_x(t)\big|_H^2dt,
\eeqa
when $x\in\SD_A$, $v_x(t)\in\Ran(C(t)^{1/2})$ a.e., and the integral is finite. Otherwise, $S(x)=+\infty$.
\end{thm}
The proof of
Theorem~\ref{thm:main} is given at the end of this section, which is established by investigating the LDP of $\la^\e$ and $M^\e$ in \secref{sec:mdp-la} and \secref{sec:mdp-M}, separately.

\begin{rem}
  Compared with the convolution type integral of the random environment for $\la^\e$ in \eqref{la}, in \cite{Gui01}, Guillin studied the deviation behavior of the
inhomogeneous additive functional of the random environment:
\bea\ad
\frac{1}{\sqrt{\e}h(\e)} \int_0^{\cdot}\hbar(s,\zeta_{s/\e})ds,
\eea
where for every $z\in\SE$, $\hbar(\cdot,z)\in
C_b([0,1];\rr^d)$ and $\pi(\hbar(t,\cdot))=0$ for every $t\geq0$.
The time dependence of $\hbar$ is controlled by the following condition on the modulus of continuity of $\hbar$:
$
\omega_{\hbar}(\dl):=\sup_{|t-s|\leq\dl,\,z\in\SE}
|\hbar(t,z)-\hbar(s,z)|$
satisfies $\lim_{\dl\to0}\omega_{\hbar}(\dl)/\sqrt{\dl}=0$. In our paper, the additive functional is replaced by a Hilbert
space-valued convolution, which is a special Volterra mapping. This result partially answers the question in
\cite[Remark~3.5]{Gui01}.
\end{rem}

\subsubsection{Auxiliary lemmas} In what follows,
we will need the following moderate deviation principle for
$m$-dependent random variables in a separable Banach space $E$; see Chen \cite{Che97}. Let
$\{a_k\}_{k\geq 1}$ be a centered, strictly stationary, $m$-dependent random sequence in
$E$ with norm $\|\cdot\|_{E}$. Denote by $E^*$ the topological dual of $E$, and for $N\in \NN_+$, denote by $E^N$ the product space with product topology.

Let $b_n>0$ be non-decreasing and satisfy $b_n/\sqrt n\to\infty$ and $b_n/n\to0$. We  assume that
$a_1\in\mathrm{WM}_0^2$, meaning that $\EE V(a_1)=0$ and
$\EE V(a_1)^2<\infty$ for every $V\in E^*$, and that there exists a $\kappa_a>0$ such that
$\EE\exp(\kappa_a \|a_1\|_E)<\infty$. Define
\bea\ad
\UU_n(t):=\frac{1}{b_n} \sum_{k=1}^{\lfloor nt \rfloor
}a_k,\quad t\in [0,1].
\eea

\begin{lem}\label{lem:A1}
Let $\CP:0=t_0<t_1<\cdots<t_N=1$ be a partition of $[0,1]$. Suppose that $\UU_n(1) \to 0$ in probability. Then
$\PP((\UU_n(t_1),\ldots,\UU_n(t_N))\in \cdot)$ satisfies an LDP in $E^N$ with speed $n/b_n^2$ and a good rate function $\CI^\CP:E^N \to [0,+\infty]$ given by
\bea\ad
\CI^\CP(z)=\sum_{i=1}^N(t_i-t_{i-1})
\sup_{V\in E^*}\Big\{V\Big(\frac{z_i-z_{i-1}}{t_i-t_{i-1}}\Big)
-\half \sigma_V^2\Big\},
\eea
for $z=(z_1,\dots,z_N), z_0 =0$, and
$ \sg_V^2=\EE V(a_1)^2+ 2\sum_{j=1}^m\EE[V(a_1)V(a_{1+j})]$.
\end{lem}

\begin{rem}\label{rem:A1}
\lemref{lem:A1} was proved in Hu and Lee \cite[Theorem~1]{HL03}
for i.i.d.\ separable Banach-valued variables $\{a_k\}_{k\geq 1}$. We use the stationary $m$-dependent partial-sum result of Chen \cite[Theorem~1.1]{Che97} to extend to this trajectory result. Its proof is postponed to Appendix~\ref{app:A1} for completeness.

\end{rem}
\begin{rem}
The condition $\EE\exp(\kappa_a \|a_1\|_E)<+\infty$ for some
$\kappa_a>0$ can be relaxed to the following tail condition: there is
$K_a>0$ such that, for every $r>0$,
\bea\ad
\limsup_{n\to\infty}n b_n^{-2} \log\PP\big(\|a_1\|_E > r b_n\big)\leq -r^2/K_a.
\eea
We refer to Guillin \cite[Lemma~5]{Gui00} and Hu and Lee
\cite[Section~4]{HL03} for detailed discussions. In this paper, we assume this exponential integrability to keep the
argument focused.
\end{rem}

To proceed, we extend Guillin \cite[Theorem 1]{Gui00} from finite-dimensional functionals to maps with values in a separable Banach space. For a fixed time, the MDP in a separable Banach space is due to Chen \cite[Chapter~IV, Theorems~2.1 and~4.1]{Che99}, and the lemma below gives its path version. Let $\{X_k\}_{k\geq 0}$ be a geometrically ergodic Markov chain on $\SE$ with transition probability $P$.
For a bounded measurable map $\fT:\SE\to E$ with $\pi(\fT)=0$, define
\bea\ad
\MM_n(t)(\fT)=\frac{1}{b_n} \sum_{k=0}^{\lfloor nt\rfloor -1}\fT(X_k), \quad t\in [0,1].
\eea

\begin{lem}\label{lem:B-mdp}
Let $\MM_n(1)(\fT)\to0$ in $\PP_\pi$-probability and let the initial law $\mu$ of $X_0$ satisfy
\begin{equation}\label{chain-cond-mu}
 \exists\; \varsigma>1:\; \sum_{n\geq1}\varsigma^n\,\|\mu P^n-\pi\|_{\rm var}<\infty.
\end{equation}
Then $\PP_\mu(\MM_n(\cdot)(\fT)\in \cdot)$ satisfies an LDP in $D^0([0,1];E)$ equipped with the uniform metric, with speed $n/b_n^2$ and good rate function
\begin{equation}\label{rate}
 \wdt J_{\fT}(\ga)=\begin{cases}
\int_0^1 J_\fT(\dot \ga(t))dt,& \text{ where }
 \ga \in \mathcal{AC},\\
 +\infty,&\text{otherwise}.
 \end{cases}
\end{equation}
Here $\mathcal{AC}:=\{\ga\in D([0,1];E):$ there exists $\dot\ga\in L^1([0,1];E)$ such that $\ga(t)=\int_0^t \dot{\ga}(s)ds$, $\forall\, t\in [0,1]\}$, where $L^1([0,1];E)$ is the space of $E$-valued Bochner integrable functions on the interval $[0,1]$. Moreover, $J_\fT(z)$ is defined as
\begin{equation}\label{JT}
 J_\fT(z)=\sup_{V\in E^*}\{V(z)-\tfrac12\sigma^2(V(\fT))\},\quad z\in E,
\end{equation}
where $\sg^2(V(\fT))$ is characterized by
\beqa{iden}\ad
 \sigma^2(V(\fT))=\pi[(V(\fT))^2] +2\sum_{j=1}^\infty\pi[V(\fT) P^jV(\fT)].
\eeqa
\end{lem}

The proof uses the regeneration split chain technique and
Lemma~\ref{lem:A1}; see Guillin \cite[Theorem 1]{Gui00} and de Acosta \cite[Proposition 2.2]{Aco97}. Since this technique plays a crucial role in the rest of the paper, the proof is given in Appendix~\ref{app:Bmdp} for
completeness. For details of this technique, we refer to
Nummelin \cite{Num84} and references therein.

\begin{cor}\label{cor:zeta-path}
Let $\SF$ be a separable Hilbert space and let $\ST:\SE\to\SF$ be a bounded map with $\pi(\ST)=0$. Then $
(\sqe h(\e))^{-1} \int_0^{\cdot}\ST (\zeta_{s/\e})ds $ satisfies an LDP in $D^0([0,1];\SF)$ with the uniform metric, with speed $h^{-2}(\e)$ and good rate $\int_0^1 L_{C_\ST}(\dot{\ga}(t))dt$ if $\ga\in \mathcal{AC}$; $+\infty$ otherwise. Here $C_{\ST}$ is the Green--Kubo covariance \eqref{CF-GK}.
\end{cor}
\begin{proof}
Let $\zeta^{(k)}:=\{\zeta_{k+s}\}_{s\in [0,1]}$, $k\geq0$, be the unit-time path chain, with values in the space $\wdh \SE=D([0,1];\SE)$ of c\`{a}dl\`{a}g $\SE$-valued segments, equipped with its coordinate $\sg$-field, and with transition kernel $\wdh P$. For a generic segment $\omega\in \wdh \SE$, define $\wdh \ST(\omega)=\int_0^1 \ST(\omega(s))ds$. Let $\wdh \pi$ be the stationary segment law and $\wdh\mu$ the initial segment law. Since $\zeta$ is stationary under $\PP_\pi$, $\wdh \pi(\wdh \ST)=\int_0^{1} \pi(\ST)ds =0$. Given $\zeta^{(0)}=\omega$, the segment $\zeta^{(j)}$ starts from the endpoint $\omega(1)$ after an additional time $j-1$. Hence, for $j \geq 1$, the Markov property and \eqref{exp-erg} give
\bea
\|\operatorname{Law}(\zeta^{(j)}\mid \zeta^{(0)}=\omega)-\wdh \pi\|_{\rm var} \ad \leq \phi_\zeta(\omega(1))\rho^{j-1} \\
\|\wdh \mu \wdh P^j- \wdh \pi\|_{\rm var} \ad \leq \|\mu P_j -\pi\|_{\rm var}\leq \mu(\phi_\zeta) \rho^j.
\eea
Thus, the path chain is geometrically ergodic and satisfies \eqref{chain-cond-mu}. Since $\SF$ is a Hilbert space, stationarity and \eqref{exp-erg} give $\EE_\pi|\int_0^n\ST(\zeta_s)ds|_\SF^2\leq Kn$, so $b_n^{-1}\sum_{k=0}^{n-1}\wdh\ST(\zeta^{(k)})\to0$ in $\PP_\pi$-probability. We then apply \lemref{lem:B-mdp} to $\wdh \ST$. Take $n=\lfloor \e^{-1}\rfloor$ and $b_n=h(\e)/\sqe$, then $b_n/\sqrt{n}\to \infty, b_n/n\to 0$ and $n/b_n^2 = h^{-2}(\e)(1+o(1))$. Passing to subsequences of $\e\to0$ and interpolating log-linearly, we may assume that $\{b_n\}$ is non-decreasing. Moreover,
\bea\ad
\sup_{0\leq t\leq 1}\Big\|\frac{1}{\sqe h(\e)} \int_0^t \ST(\zeta_{s/\e})ds - \frac{1}{b_n} \sum_{k=0}^{\lfloor nt \rfloor -1} \wdh \ST(\zeta^{(k)})\Big\|_{\SF} \\
\aad\quad =\sup_{0\leq t\leq 1}\frac{1}{b_n}\Big\|\int_{\lfloor nt \rfloor}^{t/\e} \ST(\zeta_s)ds \Big\|_{\SF} \leq \frac{2\sup_z\|\ST(z)\|_{\SF}}{b_n} \to 0.
\eea
Hence the two processes are exponentially equivalent in the uniform metric. Finally, for every $u\in \SF$, \eqref{iden} and the absolute summability implied by \eqref{exp-erg} give
\bea\ad
\sg^2(\qv{\wdh \ST,u}_\SF)=\lim_{R\to \infty} \frac{1}{R}\EE_\pi \Big|\int_0^R \qv{\ST(\zeta_s),u}_\SF ds\Big|^2=\qv{C_\ST u, u}_\SF.
\eea
Substituting it into \eqref{JT} gives $L_{C_\ST}$ and then proves the stated rate function.
\end{proof}
We also need the following consequence of \corref{cor:zeta-path}.

\begin{lem}\label{lem:averaging}
Let $b:[0,1]\times\SE\to\rr$ be bounded and centered under $\pi$, with
$\sup_z|b(t,z)-b(s,z)|\to0$ as $|t-s|\to0$.
Then, for every $r>0$,
\beq{superavg-cont}
 \limsup_{\e\to0}h^{-2}(\e)\log\PP_\mu
 \Big(\sup_{t\leq1}\Big|\int_0^tb(s,\zeta_{s/\e})\,ds\Big|>r\Big)
 =-\infty.
\eeq
\end{lem}

\begin{proof}
First suppose that $b(t,z)=b(z)$ does not depend on $t$. By \corref{cor:zeta-path} with $\SF=\rr$ and $\ST=b$, $\Ga_\e^{\rm hom}(b,\cdot):=(\sqe h(\e))^{-1}\int_0^{\cdot}b(\zeta_{s/\e})ds$ satisfies an LDP in $D^0([0,1];\rr)$ with speed $h^{-2}(\e)$ and rate function $\int_0^1L_{C_b}(\dot\ga(t))dt$. By \eqref{LC}, $L_{C_b}(v)\geq\sup_{y\in\rr}\{vy-\half(C_b+1)y^2\}=v^2/(2(C_b+1))$, hence $\|\ga\|_\infty^2\leq\int_0^1|\dot\ga(t)|^2dt\leq 2(C_b+1)\int_0^1L_{C_b}(\dot\ga(t))dt$ for $\ga\in\mathcal{AC}$.
Fix $R>0$. Since $\sqe h(\e)\to0$, the event in \eqref{superavg-cont} is contained in $\{\|\Ga_\e^{\rm hom}(b,\cdot)\|_\infty\geq R\}$ for small $\e$, and the LDP upper bound on this closed set gives
\bea\ad
\limsup_{\e\to0}h^{-2}(\e)\log\PP_\mu\big(\|\Ga_\e^{\rm hom}(b,\cdot)\|_\infty\geq R\big)\leq-\inf_{\|\ga\|_\infty\geq R}\int_0^1L_{C_b}(\dot\ga(t))dt\leq-\frac{R^2}{2(C_b+1)}.
\eea
Letting $R\to\infty$ gives \eqref{superavg-cont}.

For general $b$, let $\omega_b(\dl):=\sup_{|t-s|\leq\dl,\,z\in\SE}|b(t,z)-b(s,z)|$ and choose $p\geq1$ with $\omega_b(1/p)\leq r/2$. Replacing $b(t,z)$
on $[j/p,(j+1)/p)$ by $b(j/p,z)$ and splitting the integral at the points $j/p$ gives
$\sup_{t\leq1}|\int_0^tb(s,\zeta_{s/\e})ds|\leq\omega_b(1/p)+2\sum_{j=0}^{p-1}\sup_{t\leq1}|\int_0^tb(j/p,\zeta_{s/\e})ds|$.
Hence the event in \eqref{superavg-cont} is contained in $\bigcup_{j=0}^{p-1}\{\sup_{t\leq1}|\int_0^tb(j/p,\zeta_{s/\e})ds|>r/(4p)\}$, and the first part, applied to each $b(j/p,\cdot)$, completes the proof.
\end{proof}

We use the following Bernstein inequality to control the maximal partial sums of the Hilbert-valued martingale arising in the regeneration argument.

\begin{lem}[Bernstein inequality in a Hilbert space]\label{lem:Hil-Ber} Let $\{\beta_k\}_{k=1}^{N_1}$ be a
martingale difference sequence
in a separable Hilbert space $H$ with respect to the filtration $\{\CF_k^\beta\}_{k=1}^{N_1}$. Suppose that there are constants $\kappa_0,\kappa_v >0$ and nonnegative $\CF_{k-1}^\beta$-measurable random variables $v_k^\beta$ such that
for every integer $p\geq 2$ and $1\leq k \leq N_1$,
  \beqa{Hil-Ber-mom}\ad
 \EE(|\beta_k|^p\mid\CF_{k-1}^\beta)
 \leq p!\, v_k^\beta \, \kappa_0^{p-2}/2,
 \quad \sum_{k=1}^{N_1} v_k^\beta \leq \kappa_v \ \hbox{ a.s.}
\eeqa
Then for every $r>0$,
\begin{equation}\label{Hil-Ber-tail}
 \PP\Big(\max_{1\leq j \leq N_1}\Big|\sum_{k=1}^j \beta_k\Big|>r \Big)
 \leq2\exp\Big(-\frac{r^2}{2 \kappa_v+ 2\kappa_0 r}\Big).
\end{equation}
In particular, if $|\beta_k|_H \leq \kappa_\beta$ almost surely for every $k$, where $\kappa_\beta>0$, then
\bea\ad
 \PP\Big(\max_{1\leq j\leq N_1}\Big|\sum_{i=1}^{j} \beta_i \Big|>r\Big)
 \leq2\exp\Big(-\frac{r^2}{2(N_1 \kappa_\beta^2+ \kappa_\beta r)}\Big).
\eea
\end{lem}

\begin{proof}
For $j\ge0$, set $f_j=\sum_{k=1}^{j\wedge N_1}\beta_k$, with $f_0=0$. Then $\{f_j\}_{j\ge0}$ is a martingale with respect to $\{\CF_{j\wedge N_1}^\beta\}_{j\ge0}$. Its increments are $\beta_j$ for $j\le N_1$ and zero thereafter. Hence $
\sup_{j\ge0}|f_j|_H
=\max_{1\le j\le N_1}|\sum_{k=1}^j\beta_k|_H,
$ and, for every integer $p\ge2$,
\bea\ad
\Big\|\sum_{j\ge1}\mathbb E\bigl(
|f_j-f_{j-1}|_H^p\mid\CF_{(j-1)\wedge N_1}^\beta \bigr)\Big\|_{L^\infty}
=\Big\|\sum_{j=1}^{N_1}\EE\bigl(
|\beta_j|_H^p\mid\CF_{j-1}^\beta\bigr)\Big\|_{L^\infty}
\le \frac{p!}{2}\kappa_v\kappa_0^{p-2}.
\eea
The Hilbert norm is $(2,1)$-smooth by the parallelogram identity. Thus the moment condition in Pinelis \cite[Theorem 3.3]{Pin94} holds with $D=1$, $B^2=\kappa_v$, and $\Gamma=\kappa_0$. Using the corrected bound in \cite{Pin99}, we obtain
\bea\ad
\PP\Big(\sup_{j\ge0}|f_j|_H>r\Big)
\le 2\exp\Big(
-\frac{r^2}{\kappa_v+\kappa_0r+
\sqrt{\kappa_v(\kappa_v+2\kappa_0r)}}\Big).
\eea
Since $\sqrt{\kappa_v(\kappa_v+2\kappa_0r)}\le\kappa_v+\kappa_0r$, this gives the claimed Bernstein bound \eqref{Hil-Ber-tail}. If $|\beta_k|_H\le\kappa_\beta$, then $
\EE(|\beta_k|_H^p\mid\CF_{k-1}^\beta)
\le \kappa_\beta^{p-2} \EE(|\beta_k|_H^2\mid\CF_{k-1}^\beta)$. Taking $v_k^\beta=\mathbb E(|\beta_k|_H^2\mid\CF_{k-1}^\beta)$, $\kappa_0=\kappa_\beta$, and $\kappa_v=N_1 \kappa_\beta^2$ verifies the moment assumptions and proves the bounded-increment case.
\end{proof}

\subsection{\texorpdfstring{The LDPs of $\la^\e$, $M^\e$, and their sum}{The LDPs of lambda epsilon, M epsilon, and their sum}}
\subsubsection{\texorpdfstring{The LDP of $\la^\e$}{The LDP of lambda epsilon}}\label{sec:mdp-la}

Fix $N\in\NN_+$ and set $t_k=k/N$, $0\leq k\leq N$ and $\Dl:=1/N$. The $\la^\e$ can be decomposed as
\beq{decomp-la}\barray
\la^\e(t)
\ad = \feh \sum_{k=0}^{\lfloor Nt \rfloor-1}\int_{t_k}^{t_{k+1}} \big[e^{(t-s)A} F_c(s,\zeta_{s/\e})-e^{(t-t_k) A}  F_c(t_k,\zeta_{s/\e})\big]ds \\
\aad+\feh \sum_{k=0}^{\lfloor Nt \rfloor-1} \int_{t_k}^{t_{k+1}} e^{(t-t_k) A}F_c(t_k,\zeta_{s/\e})ds\!+\!\feh \int_{\lfloor Nt \rfloor/N}^t e^{(t-s)A}  F_c(s,\zeta_{s/\e})ds\\
\ad=: D_{\e}^{N}(t)+L_\e^{N} (t)+A_\e^{N}(t)
\earray
\eeq
In what follows,
we will prove the LDP for $L_\e^N$ in \lemref{lem:LDP-L} and that both $D_\e^N$ and $A_\e^N$ are negligible with respect to the large deviation principles as $N$ goes to infinity in \lemref{lem:neg-D} and \lemref{lem:neg-A}, respectively.

\begin{lem}[LDP for $L_\e^N$]\label{lem:LDP-L}
  For fixed $N\in \NN_+$, $L_\e^N$ satisfies an LDP in $D^0([0,1];H)$ with speed $h^{-2}(\e)$ and a good rate function given by
\beq{IN}
\CI_N^{\rm hom}(\ga)=\sum_{k=0}^{N-1}\Dl\cdot L\Big(e^{\Dl A} C_F(t_k)e^{\Dl A^*}, \frac{\ga(t_{k+1})-e^{\Dl A} \ga(t_k)}{\Dl}\Big)
\eeq
for $\ga \in\mathscr{S}_N$, where $\mathscr{S}_N:=\{\ga\in D^0([0,1];H):\ga(t)=e^{(t-t_k)A}\ga(t_k), \text{ for } t\in [t_k,t_{k+1}), 0\leq k\leq N-1\}$. The rate is
$+\infty$ outside $\mathscr{S}_N$.

\end{lem}

\begin{proof}
For $L_\e^N$ defined in \eqref{decomp-la}, we note that it can be rewritten as
\bea
L_\e^{N}(t)
\ad = \sum_{k=0}^{\lfloor Nt \rfloor
-1} e^{(t-t_k)A}\Big[\feh \int_0^{t_{k+1}} F_c(t_k,\zeta_{s/\e}) ds -\feh \int_0^{t_k} F_c(t_k,\zeta_{s/\e})ds\Big] \\
\ad=: \sum_{k=0}^{\lfloor Nt \rfloor
-1} e^{(t-t_k)A} \big[\Ga_\e^{\rm hom}(F_k,t_{k+1})-\Ga_\e^{\rm hom}(F_k, t_k)\big],
\eea
where $F_k(z):=F_c(t_k,z)$ for $z\in\SE$ and
\beq{Ga-h}
\Ga_\e^{\rm hom}(F_k, \cdot): =\feh \int_0^{\cdot} F_k(\zeta_{s/\e})ds.
\eeq
Let $F^N:=(F_0,\ldots,F_{N-1}):\SE\to H^N$.  By
Guillin \cite[Theorem~1, Section~2.1]{Gui01} and
\corref{cor:zeta-path} with $\SF=H^N$ and $\ST=F^N$, $\Ga_\e^{\rm hom}(F^N,\cdot)$ satisfies an LDP in
$D^0([0,1];H^N)$ with speed $h^{-2}(\e)$ and rate
\beqa{J-N}
J_{N}^{\rm hom}(\ga)=\int_0^1 \sup_{\phi \in H^N}\big\{\qv{\dot \ga(t), \phi}_{H^N}-\frac{1}{2} \sg^2 \big(\qv{ F^N,
\phi}_{H^N}\big)\big\}dt,
\eeqa
for absolutely continuous $H^N$-valued paths $\ga$ with
$\ga(0)=0$, otherwise $J_N^{\rm hom}(\ga)=+\infty$. Here $\sg^2(\qv{F^N,\phi}_{H^N})=\qv{C_{F^N}\phi,\phi}_{H^N}$ by the proof of \corref{cor:zeta-path}.

Define the mapping
$\Psi^N:D^0([0,1];H^N)\to D^0([0,1];H)$ by
\beq{Phi-map}
\ga(t)\mapsto \Psi^N(\ga)(t)= \sum_{k=0}^{\lfloor Nt \rfloor
-1} e^{(t-t_k)A}\big[\ga_k(t_{k+1})-\ga_k(t_k)\big].
\eeq
The image of $\Psi^N$ lies in $\mathscr{S}_N$. Indeed, for $\wdt\ga\in D^0([0,1];H^N)$ and $t\in [t_j,t_{j+1})$, we have $\lfloor Nt \rfloor
=j$ and the semigroup property gives
$
\Psi^N(\wdt \ga)(t)=e^{(t-t_j)A} \sum_{k=0}^{j-1} e^{(t_j-t_k)A}[\wdt \ga_k(t_{k+1})-\wdt \ga_k(t_k)]=e^{(t-t_j)A} \Psi^N(\wdt \ga)(t_j)$. Hence $L_\e^N=\Psi^N(\Ga_\e^{\rm hom}(F^N,\cdot))\in \mathscr{S}_N$ almost surely and for $\ga\notin \mathscr{S}_N$, the set $\{\wdt \ga: \Psi^N(\wdt \ga)=\ga\}$ is empty, so $\CI_N^{\rm hom}(\ga)=\inf\emptyset=+\infty$.  Moreover, the map $\Psi^N$ is continuous in the uniform topology.  The contraction
principle therefore gives an LDP for $L_\e^N$ in $D^0([0,1];H)$,
with speed $h^{-2}(\e)$ and a rate
 function $\CI_N^{\rm hom}(\ga)$ given by
\bea
\CI_{N}^{\rm hom}(\ga)
\ad =\inf\big\{J_N^{\rm hom}(\wdt \ga): \Psi^N(\wdt \ga)=\ga \big\}.
\eea

Fix $\ga\in \mathscr{S}_N$ and define $w_k:=[\ga(t_{k+1})-e^{\Dl A} \ga(t_k)]/\Dl$. For an absolutely continuous $\wdt \ga$ starting at zero, the constraint $\Psi^N(\wdt \ga)=\ga$ is equivalent to
\beqa{Dl-int-ga}
e^{\Dl A}\Dl^{-1}\int_{t_k}^{t_{k+1}}\dot{\wdt \ga_k}(s)ds =w_k, \quad 0\leq k <N.
\eeqa
Indeed, \eqref{Phi-map} gives the grid recursion \eqref{Dl-int-ga}, and the recursion determines all grid values from $\ga(0)=0$. Both paths then agree between grid points by the definition of $\mathscr{S}_N$.

Let $\psi_0,\dots,\psi_{N-1}\in H$ be arbitrary. Take the step function $\phi(t)=(\phi_0(t),\dots,\phi_{N-1}(t))\in H^N$ given by $\phi_i(t)=0$ for $i\neq k$ and $\phi_k(t)=e^{\Dl A}\psi_k$ for $t\in [t_k,t_{k+1})$. By definition of $J_N^{\rm hom}$ in \eqref{J-N}, we have
\bea
J_N^{\rm hom}(\wdt\ga) \ad =\int_0^1 \sup_{\phi\in H^N}\big\{\bqv{\dot{\wdt \ga}(t),\phi}_{H^N}-\half \sg^2\big(\bqv{ F^N, \phi}_{H^N}\big)\big\}dt \\
\aad\geq \sum_{k=0}^{N-1} \int_{t_k}^{t_{k+1}} \bqv{\dot{\wdt\ga_k}(t), e^{\Dl A}\psi_k}_H -\half\sg^2\big(\bqv{ F_k,e^{\Dl A}\psi_k}_H \big) dt \\
\aad = \sum_{k=0}^{N-1}\Dl  \Big[\Bqv{\frac{e^{\Dl A}(\wdt\ga_k(t_{k+1})-\wdt \ga_k(t_k))}{\Dl }, \psi_k }_H - \half\sg^2\big(\bqv{e^{\Dl A} F_k, \psi_k}_H \big)\Big] \\
\aad =\sum_{k=0}^{N-1} \Dl \big[\qv{w_k, \psi_k}_H -\half \sg^2\big(\bqv{e^{\Dl A} F_k, \psi_k}_H \big) \big],
\eea
where the first equality uses $A=A^*$, and the last comes from \eqref{Dl-int-ga}. Taking the supremum over each $\psi_k$ and then the infimum over preimages gives
\bea\ad
\CI_N^{\rm hom}(\ga)\geq \sum_{k=0}^{N-1}\Dl \sup_{\psi\in H}\{\qv{w_k,\psi}_H-\half \sg^2(\qv{e^{\Dl A}F_k, \psi}_H)\}.
\eea
 By the proof of \corref{cor:zeta-path} and \lemref{lem:Fc-U},
\bea\ad
\sg^2(\qv{e^{\Dl A}F_k,\psi}_H)=\qv{C_F(t_k)e^{\Dl A^*}\psi,e^{\Dl A^*}\psi}_H=\qv{e^{\Dl A}C_F(t_k)e^{\Dl A^*}\psi,\psi}_H,
\eea
so the right-hand side equals the right-hand side of \eqref{IN}.

For the reverse inequality, suppose that $\ga\in \mathscr{S}_N$ and that the right-hand side of \eqref{IN} is finite.
The covariance of the $k$-th coordinate of $F^N$ is $C_{F_k}=C_F(t_k)$ by \lemref{lem:Fc-U}. Therefore, applying \eqref{LC1} directly to the map $v\mapsto e^{\Dl A} v_k$ gives
\bea
\ad \inf_{\substack{v\in H^N\\e^{\Delta A}v_k=w_k}}
 L_{C_{F^N}}(v)=L_{e^{\Delta A}C_F(t_k)e^{\Delta A^*}}(w_k).
\eea
Since the right-hand side is finite, \corref{cor:add} gives $v^{(k)}\in H^N$ such that $e^{\Dl A} v_k^{(k)}=w_k$ and $L_{C_{F^N}}(v^{(k)})=L_{e^{\Dl A} C_F(t_k)e^{\Dl A^*}}(w_k)$. Define $\wdt \ga(t)=\int_0^t \sum_{k=0}^{N-1} v^{(k)} \indi_{[t_k,t_{k+1})}(s) ds$. Then $e^{\Dl A} \Dl^{-1}\int_{t_k}^{t_{k+1}} \dot{\wdt \ga}_k (s)ds = e^{\Dl A} v_k^{(k)} =w_k$. Hence \eqref{Dl-int-ga} gives $\Psi^N(\wdt \ga)=\ga$. Consequently,
\bea\ad
\CI_N^{\rm hom}(\ga)\leq J_N^{\rm hom}(\wdt \ga)=\sum_{k=0}^{N-1}\Dl L_{C_{F^N}}(v^{(k)})=\sum_{k=0}^{N-1}\Dl L_{e^{\Dl A} C_F(t_k) e^{\Dl A^*}}(w_k),
\eea
which is the required upper bound. The proof is complete.

\end{proof}

\begin{lem}[Negligibility of $D_\e^N$]\label{lem:neg-D} For every $r>0$, we have
\beq{D}
\lim_{N\to \infty} \limsup_{\e\to 0} h^{-2}(\e) \log\PP_\mu(\|D_\e^{N}\|_\infty >r)=-\infty.
\eeq
\end{lem}

\begin{lem}[Negligibility of $A_\e^N$]\label{lem:neg-A}
      For every $r>0$, we have
\beq{A}
\lim_{N\to \infty}\limsup_{\e\to 0}h^{-2}(\e)\log\PP_\mu(\|A_\e^{N}\|_\infty>r) = -\infty.
\eeq
\end{lem}
Proofs of the above two lemmas are postponed to \appref{app:additive-holder}.

\begin{lem}\label{lem:la-app-LDP}
Let $B(\ga,\dl)=\{z:\sup_{t\in[0,1]}
|z(t)-\ga(t)|_H<\dl\}$. Then $\la^\e$ satisfies a weak LDP in
$D^0([0,1];H)$, with speed $h^{-2}(\e)$ and rate given by
\beq{I}
\CI_\zeta(\ga):=\sup_{\dl>0} \liminf_{N\to \infty} \inf_{z\in B(\ga, \dl)}\CI_N^{\rm hom}(z) =\sup_{\dl>0}\limsup_{N\to \infty}\inf_{z\in B(\ga,\dl)}\CI_N^{\rm hom}(z).
\eeq
\end{lem}
\begin{proof}
\lemref{lem:LDP-L}--\lemref{lem:neg-A} and \cite[Theorem~4.2.16 (a)]{DZ09} give a weak LDP for $\la^\e$ in $D^0([0,1];H)$, with speed $h^{-2}(\e)$ and rate defined in \eqref{I}.
\end{proof}

We next characterize the rate function $\CI_\zeta$ in \lemref{lem:la-app-LDP}.

\begin{lem}\label{lem:CI-CF}
With the uniform topology on $D^0([0,1];H)$, if $\ga_N\to \ga$, then
\beq{CF-leq}
\CI_{C_F}(\ga) \leq \liminf_{N\to \infty}\CI_N^{\rm hom}(\ga_N).
\eeq
Conversely, for every $\ga \in C^0([0,1];H)$, there are paths $\ga_N$ such that $\ga_N \to \ga$ and
\beq{CF-geq}
\limsup_{N\to \infty}\CI_N^{\rm hom}(\ga_N) \leq \CI_{C_F}(\ga).
\eeq
When $\CI_{C_F}(\ga)<\infty$, the paths can be chosen in $\mathscr{S}_N$. Moreover, if $\sup_{j} \CI_{N_j}^{\rm hom} (\ga_{N_j})<\infty$ and $N_j\to \infty$,  then $\{\ga_{N_j}\}$ is relatively compact, and all its limit points belong to $C^0([0,1];H)$.
\end{lem}

The proof is postponed to \appref{app:lem-CI-CF}.

\begin{lem}[%
{Identification of $\CI_\zeta$}]\label{lem:ind-I-zeta}
We have $\CI_\zeta(\ga)=\CI_{C_F}(\ga)$ for every $\ga\in D^0([0,1];H)$, with the explicit form given by \eqref{quad-L} and \eqref{LC-t}.
\end{lem}
\begin{proof}
We first show $\CI_\zeta \leq \CI_{C_F}$. Suppose $\CI_{C_F}(\ga)<\infty$, so $\ga\in \SD_A$, thus $\ga\in C^0([0,1];H)$. By \eqref{CF-geq}, there exists $\ga_N$ with $\|\ga_N-\ga\|_\infty \to 0$ and $\limsup_{N\to \infty}\CI_N^{\rm hom}(\ga_N)\leq \CI_{C_F}(\ga)$. For each $\dl>0$, we have $\ga_N \in B(\ga,\dl)$ for all large $N$, hence
\bea\ad
\liminf_{N\to \infty}\inf_{z\in B(\ga,\dl)}\CI_N^{\rm hom}(z)\leq \limsup_{N\to \infty}\CI_N^{\rm hom}(\ga_N) \leq \CI_{C_F}(\ga),
\eea
and the supremum over $\dl$ gives $\CI_\zeta(\ga) \leq \CI_{C_F}(\ga)$. Conversely, let $\CI_\zeta(\ga)<\infty$ and fix $\eta,\dl>0$. By the definition of $\CI_\zeta$, $\liminf_{N\to \infty} \inf_{B(\ga,\dl)}\CI_N^{\rm hom} \leq \CI_\zeta(\ga)$, so there is a subsequence $N_j\to \infty$ such that $\inf_{B(\ga,\dl)}\CI_{N_j}^{\rm hom}<\CI_\zeta(\ga)+\eta/2$ for every $j$. For each $j$, choose $\ga_{N_j}\in B(\ga,\dl)$ with cost at most $\eta/2$ above this infimum. Then
$
\sup_{j}\CI_{N_j}^{\rm hom}(\ga_{N_j}) \leq \CI_\zeta(\ga)+\eta<\infty.
$
By the last assertion in \lemref{lem:CI-CF}, there exists a further subsequence (still denoted by $\ga_{N_j}$ with certain abuse of notation), converging uniformly to some $\ga_\dl\in C^0([0,1];H)$, and $\|\ga_\dl-\ga\|_\infty \leq\dl$. The proof of \eqref{CF-leq} only uses $\ga_N\to \ga$ along the indices, so it applies along $N_j$ and gives
\bea\ad
\CI_{C_F}(\ga_\dl)\leq \liminf_{j\to\infty} \CI_{N_j}^{\rm hom}(\ga_{N_j}) \leq \CI_\zeta(\ga)+\eta.
\eea
Now let $\dl\to 0$, then $\ga_\dl \to \ga$ uniformly, so $\ga\in C^0([0,1];H)$ because $C^0([0,1];H)$ is closed in $D^0([0,1];H)$. Since $\CI_{C_F}$ is a good rate function on $C^0([0,1];H)$ by \lemref{lem:quad-lift}, it is lower-semicontinuous, and
\bea\ad
\CI_{C_F}(\ga)\leq \liminf_{\dl\to 0} \CI_{C_F}(\ga_\dl) \leq \CI_\zeta(\ga)+\eta.
\eea
Letting $\eta\downarrow 0$ gives $\CI_{C_F}\leq \CI_\zeta$. The proof is complete.
\end{proof}

Lemma~\ref{lem:la-app-LDP} gives the weak LDP for $\lambda^\e$ with rate $\CI_\zeta$, and Lemma~\ref{lem:ind-I-zeta} identifies this rate as $\CI_{C_F}$. The full LDP follows.

\begin{thm}\label{thm:la}
Under Assumptions~\ref{ass:A} and~\ref{ass:fg}, the environmental
conditions of Section~\ref{ass:zeta}, and
$x_0\in\CD((-A)^\theta)$ for some $\theta\in(1/2,1)$, the family
$\lambda^\e$ satisfies an LDP in $C^0([0,1];H)$ with speed $h^{-2}(\e)$ and good rate
$
 \CI_\zeta(x)=\CI_{C_F}(x).$
Equations \eqref{quad-L}--\eqref{LC-t} give
the explicit
quadratic expression and the
full domain.
\end{thm}
\begin{proof}
By \lemref{lem:ind-I-zeta} and \lemref{lem:quad-lift}, $\CI_\zeta=\CI_{C_F}$ is good. Let $\CO\subset D^0([0,1];H)$ be closed and let $\kappa>\limsup_{N\to\infty}\inf_{\CO}\CI_N^{\rm hom}$. Choose $N_j\to\infty$ and $\ga_{N_j}\in\CO$ with $\CI_{N_j}^{\rm hom}(\ga_{N_j})\leq\kappa$. By \lemref{lem:CI-CF}, a subsequence converges uniformly to some $\ga\in\CO$ with $\CI_{C_F}(\ga)\leq\kappa$. Hence $\inf_{\CO}\CI_\zeta\leq\limsup_{N\to\infty}\inf_{\CO}\CI_N^{\rm hom}$, and \lemref{lem:la-app-LDP} and \cite[Theorem~4.2.16 (b)]{DZ09} give the full LDP in $D^0([0,1];H)$. Since $\la^\e$ has continuous paths and $\CI_\zeta=+\infty$ outside the closed subset $C^0([0,1];H)$ of $D^0([0,1];H)$, \cite[Lemma~4.1.5 (b)]{DZ09} gives the LDP in $C^0([0,1];H)$.
\end{proof}

\subsubsection{\texorpdfstring{The LDP of $M^\e$}{The LDP of M epsilon}}\label{sec:mdp-M} We establish the LDP of $M^\e$ in $C^0([0,1];H)$. The main difficulty is that $M^\e$ is not an $H$-valued semimartingale due to the semigroup $e^{(t-s)A}$ in the integrand of $M^\e$. Thus, the finite-dimensional predictable characteristics argument of Guillin \cite[Proposition 1]{Gui03} using results of Puhalskii \cite{Puh94} cannot be applied directly. In finite-dimensional space, large deviation principles for small perturbations of random evolution equations
\bea
dx^\e(t)=b(x^\e(t), y(t))dt+ \sqe \sg(x^\e(t))dw(t)
\eea
have been proved by Bezuidenhout \cite{Bez87}, Hu \cite{Hu97}, and Mellouk \cite{Mel00}, where $y$ is a random process whose support is required to be compact.

In our case, $\sqe h(\e)M^\e$ is the mild solution of the random evolution equation $dm(t)=Am(t)dt+\sqrt\e G(\bar X(t),\zeta_{t/\e})dW^{\mathsf{Q}}(t)$, $m(0)=0$, in infinite-dimensional space.

As outlined in the introduction, we project the stochastic convolution onto the spectral subspaces $H_n$, establish the LDP of the projected martingale from its predictable quadratic variation, and recover the projected convolution by the continuous Volterra map $\Phi_n$ in \eqref{def-Phi-n}.

\begin{lem}\label{new:averaging}
For every fixed $n\in \NN_+$ and $r>0$, with $C_{G,n}(s):=\Pi_nC_G(s)\Pi_n|_{H_n}$,
\beqa{new:qv}\disp
 \limsup_{\e\to0}h^{-2}(\e)\log\PP_\mu\Big(
 \sup_{t\leq1}\big\|\int_0^t
 \big[\ad \Pi_n G(\bar X(s),\zeta_{s/\e})\mathsf{Q}\mathsf{Q}^*G(\bar X(s),\zeta_{s/\e})^*\Pi_n\\
 \aad\quad -C_{G,n}(s)\big] \,ds
 \big\|>r\Big)=-\infty.
\eeqa
\end{lem}
The proof is postponed to \appref{app:new-averaging}. The next lemma gives the LDP of the projected martingale and then reintroduces the semigroup.
\begin{lem}\label{new:finite}
For each $n\in \NN_+$, we have the following results.
\begin{enumerate}
\item[(a)] The continuous $H_n$-valued martingale
$
Z_{G,n}^\e(t):=h^{-1}(\e)\int_0^t \Pi_n G(\bar X(s),\zeta_{s/\e})dW^{\mathsf{Q}}(s)$
satisfies an LDP in $C^0([0,1];H_n)$ with speed $h^{-2}(\e)$ and good
rate function
\begin{equation*}
 J_{G,n}(z)=\begin{cases}
\int_0^1 L_{C_{G,n}}(t,\dot z(t))\,dt,
 &\text{ if } z\in W^{1,2}([0,1];H_n),\ z(0)=0,\\
 +\infty,&\text{otherwise}.
 \end{cases}
\end{equation*}
\item[(b)] $\Pi_n M^\e$ satisfies an LDP in $C^0([0,1];H_n)$ with
speed $h^{-2}(\e)$ and good rate function
\begin{equation}\label{new:IN}
\CI_{G,n}(x)=\begin{cases}
\int_0^1L_{C_{G,n}}(t,\dot x(t)-A_nx(t))dt,
 &\text{if } x\in W^{1,2}([0,1];H_n),x(0)=0\\
 +\infty,&\text{otherwise}.
 \end{cases}
\end{equation}
\end{enumerate}
\end{lem}

The proof is postponed to \appref{app:new-finite}.
\begin{lem}\label{new:tail}
Suppose that Assumptions~\ref{ass:A} and~\ref{ass:fg} hold. Then $\Pi_n M^\e$
is an exponentially good approximation of $M^\e$
in $C^0([0,1];H)$. That is, for every $r>0$,
\beq{new:exp-good}
\lim_{n\to\infty}\limsup_{\e\to0}
h^{-2}(\e)\log\PP\big(
\|(I-\Pi_n)M^\e\|_\infty>r\big)=-\infty.
\eeq
Moreover, for $\delta\in(0,1-\vartheta)$, there exist constants $K_1,\kappa_\delta>0$ such that
\beq{new:cond-tail}
\PP\big(\|(I-\Pi_n)M^\e\|_\infty>r
\,\big|\,\zeta\big)
\leq K_1\exp\big(-\kappa_\delta h^2(\e)r^2\alpha_{n+1}^{\delta}\big)\quad\text{a.s., for all } n\in \NN_+,\ r\geq0.
\eeq
\end{lem}
\begin{proof}
  Fix $\delta\in(0,1-\vartheta)$. Since $\bar X$ is bounded by \thmref{thm:XbarX} (2), \eqref{abstract-Gtail} gives
\bea
\|e^{uA}(I-\Pi_n)G(\bar X(s),z)\mathsf Q\|_{\CL_2(H)}^2
\le K e^{-\alpha_{n+1}u}u^{-\vartheta}
\le K_\delta\alpha_{n+1}^{-\delta}u^{-(\vartheta+\delta)}
\eea
uniformly in $s,z$ and $0<u\le1$, where the second inequality uses $\sup_{x>0}x^\delta e^{-x}<\infty$.

By \eqref{M}, $h(\e)(I-\Pi_n)M^\e$ is the stochastic convolution driven by $W_0$ with integrand $(I-\Pi_n)G(\bar X(s),\zeta_{s/\e})\mathsf Q$. Given $\zeta$, the integrand is deterministic and $W_0$ retains the law of a Wiener process by independence. Apply \lemref{lem:conv-exp} with exponent $\vartheta+\delta$ and threshold $h(\e)r$. With $\kappa_\delta:=(K_2K_\delta)^{-1}$, this gives
\bea\ad
\PP\bigl(\|(I-\Pi_n)M^\e\|_\infty>r\mid\zeta\bigr)
\le K_1\exp\bigl(-\kappa_\delta h^2(\e)r^2\alpha_{n+1}^{\delta}\bigr),
\qquad r\ge0.
\eea
Taking expectations gives the same bound unconditionally. For $r>0$, the upper exponential rate is therefore at most $-\kappa_\delta r^2\alpha_{n+1}^{\delta}$. Since $\alpha_{n+1}\to\infty$, this proves \eqref{new:exp-good}.
\end{proof}

Lemma~\ref{new:averaging} and Lemma~\ref{new:finite} (b) verify condition~(i) of
Corollary~\ref{thm:abstract-lift}, and
Lemma~\ref{new:tail} verifies condition~(ii).  The corollary therefore
gives an LDP for $M^\e$ with rate $\CI_{C_G}$.

\begin{thm}\label{thm:mdp-comb}
Under Assumptions~\ref{ass:A} and~\ref{ass:fg}, suppose the Markov
process satisfies the path, independence, and mixing conditions in
Section~\ref{ass:zeta}.  If \eqref{dev-s}
holds, then $M^\e$ satisfies an LDP with good rate function in $C^0([0,1];H)$, with speed
$h^{-2}(\e)$ and explicit rate
\begin{equation}\label{new:IW}
 \CI_w(x)=\CI_{C_G}(x).
\end{equation}
The explicit cost and its domain are given by \eqref{quad-L} and \eqref{LC-t}.
\end{thm}

\subsubsection{\texorpdfstring{The LDP for the sum of two fluctuations}{The LDP of R epsilon}}\label{sec:mdp-R}
We prove the LDP for $R^\e=\la^\e+M^\e$ by establishing its exponential
tightness and a local LDP, which implies the full LDP.

Since $C^0([0,1];H)$ is Polish and the rate functions in Theorems~\ref{thm:la} and~\ref{thm:mdp-comb} are good, $\la^\e$ and $M^\e$ are exponentially tight, and hence so is their sum $R^\e$. Define the rate by
\begin{equation}\label{new:sumrate}
 \CI(x)=\inf_{y\in C^0([0,1];H)}\{\CI_\zeta(x-y)+\CI_w(y)\}.
\end{equation}
Since $W_0$ is independent of $\zeta$, conditionally on $\zeta$ the process $M^\e$ is the stochastic convolution of the deterministic integrand $s\mapsto h^{-1}(\e)G(\bar X(s),\zeta_{s/\e})\mathsf{Q}$. By the definition of $Z_{G,n}^\e$, conditionally on $\zeta$, $Z^\e_{G,n}$ is a continuous Gaussian martingale with covariance $h^{-2}(\e)\int_0^\cdot\Pi_n G(\bar X(s),\zeta_{s/\e})\mathsf{Q}\mathsf{Q}^*G(\bar X(s),\zeta_{s/\e})^*\Pi_n\,ds$. By \lemref{new:averaging}, $h^2(\e)$ times this covariance is uniformly close to $\int_0^\cdot C_{G,n}(s)\,ds$, outside an event of superexponentially small probability. In the finite-dimensional setting of \cite{Gui03} this is exactly the situation of the conditional Schilder theorem \cite[Lemma~2]{Gui03}, which gives two-sided small-ball estimates for $M^\e$ given the environment. The next lemma is its analogue in $C^0([0,1];H)$. Below, $B(\phi,a)$ and $\overline B(\phi,a)$ denote the open and the closed ball in $C^0([0,1];H)$ with center $\phi$ and radius $a$.

\begin{lem}[Conditional Schilder estimate]\label{lem:conditional-M}
There exist $\CF^\zeta$-measurable events $\Omega_\e$ with
$
 \lim_{\e\to0}h^{-2}(\e)\log\PP(\Omega_\e^c)=-\infty
$
such that the following holds. Let $\phi\in C^0([0,1];H)$ satisfy $\CI_w(\phi)<\infty$, and let $\eta>0$. There is $a_0=a_0(\phi,\eta)>0$ such that for every $a\in(0,a_0]$ and all sufficiently small $\e$, depending on $\phi$, $\eta$ and $a$, almost surely on $\Omega_\e$,
\begin{align}
 \PP\big(M^\e\in B(\phi,a)\,\big|\,\zeta\big)
 &\geq\exp\big\{-h^2(\e)\big(\CI_w(\phi)+\eta\big)\big\},\label{cond-Sch-low}\\
 \PP\big(M^\e\in\overline B(\phi,a)\,\big|\,\zeta\big)
 &\leq\exp\big\{-h^2(\e)\big(\CI_w(\phi)-\eta\big)\big\}.\label{cond-Sch-up}
\end{align}
The lower bound \eqref{cond-Sch-low} holds for every $a>0$.
\end{lem}

The proof is postponed to \appref{app:lem-conditional-M}.

\begin{lem}\label{lem:local-R}
The family $R^\e$ satisfies an LDP in
$C^0([0,1];H)$ with speed $h^{-2}(\e)$ and good rate
$\CI$ defined in \eqref{new:sumrate}.  Set $C=C_F+C_G$.  Then $
 \CI(x)=\CI_C(x)$ for every $x\in C^0([0,1];H)$, with the explicit form given by \eqref{quad-L} and \eqref{LC-t}.
\end{lem}

\begin{proof}
Fix $\ga\in C^0([0,1];H)$ and $a,\eta>0$. To prove the lower bound of the local LDP of $R^\e$, we suppose $\CI(\ga)<\infty$, otherwise, it is trivial. By the definition of $\CI(\ga)$ in \eqref{new:sumrate}, for any $\eta>0$, we can choose $\phi\in C^0([0,1];H)$ such that $\CI_\zeta(\ga-\phi)+\CI_w(\phi)<\CI(\ga)+\eta$. In particular, $\CI_w(\phi)<\infty$.
Let $\widetilde a:=a/4$. By \eqref{cond-Sch-low}, which holds for every radius, for all sufficiently small $\e$, we have
\bea
\PP(M^\e\in B(\phi,\widetilde a)\,|\,\zeta)
\geq \exp\{-h^2(\e)(\CI_w(\phi)+\eta)\}\ \text{ a.s. on } \Omega_\e.
\eea
Since $B(\ga-\phi,a/2)+B(\phi,\widetilde a)\subset B(\ga,a)$, and since
$\la^\e$ and $\Omega_\e$ are $\CF^\zeta$-measurable,
\beq{low-est}\barray
\PP\big(R^\e\in B(\ga,a)\big)
\ad\geq\EE\Big[\indi_{\{\la^\e\in B(\ga-\phi,a/2)\}}\indi_{\Omega_\e}\,
 \PP\big(M^\e\in B(\phi,\widetilde a)\,\big|\,\zeta\big)\Big]\\
\aad \geq \exp\big\{-h^2(\e)(\CI_w(\phi)+\eta)\big\}
 \big\{\PP\big(\la^\e\in B(\ga-\phi,a/2)\big)-\PP(\Omega_\e^c)\big\}.
\earray\eeq
The LDP lower bound for $\la^\e$ and
$h^{-2}(\e)\log\PP(\Omega_\e^c)\to-\infty$ imply that
\bea\ad
 \liminf_{\e\to0}h^{-2}(\e)\log\PP\big(R^\e\in B(\ga,a)\big)
 \geq-\CI_\zeta(\ga-\phi)-\CI_w(\phi)-\eta
 \geq-\CI(\ga)-2\eta.
\eea
Letting $\eta\to0$ proves the local lower bound.
\smallskip

To prove the upper bound of the local LDP, for any $\kappa >0$, the level set
$\CO_\kappa=\{\phi:\CI_w(\phi)\leq \kappa \}$ is compact due to the goodness of
$\CI_w$.  Fix $0<\eta<1/2$. By \eqref{cond-Sch-up}, for every $\phi\in \CO_\kappa$ there
exists $\widetilde a(\phi)>0$ such that for all sufficiently small
$\e$,
\bea
 \PP(M^\e\in\overline B(\phi,\widetilde a(\phi))\,|\,\zeta)
 \leq\exp\{-h^2(\e)(\CI_w(\phi)-\eta)\}\ \text{ a.s. on } \Omega_\e.
\eea
Moreover, the LDP of $\la^\e$ implies that there exists
$\widehat a(\phi)>0$ such that
\bea\ad
 \PP\big(\la^\e\in\overline B(\ga-\phi,\widehat a(\phi))\big)
 \leq\exp\big\{-h^2(\e)\chi_\eta \big(\CI_\zeta(\ga-\phi)\big)\big\},
\eea
for all sufficiently small $\e$, where $\chi_\eta(x)=(x-\eta)\wedge\eta^{-1}$.  Here $\chi_\eta$ is introduced to avoid the term $\CI_{\zeta}(\ga-\phi)$ being possibly infinite. By the compactness of
$\CO_\kappa$, we can choose $\phi_1,\ldots,\phi_{N_\kappa}\in \CO_\kappa$ such that
\bea\ad
 \CO_\kappa\subset \CO_{N_\kappa}^\kappa=\bigcup_{k=1}^{N_\kappa}B(\phi_k,a_k),\qquad
 a_k=\widetilde a(\phi_k)\wedge\frac{\widehat a(\phi_k)}2.
\eea
Let $0<a\leq\min_ka_k$ and fix $k$. On the event $\{R^\e\in\overline B(\ga,a),\ M^\e\in B(\phi_k,a_k)\}$ we have $\|\la^\e-(\ga-\phi_k)\|_{\infty}\leq a+a_k\leq\widehat a(\phi_k)$ and $M^\e\in\overline B(\phi_k,\widetilde a(\phi_k))$. Conditioning on the environment and using the two preceding estimates, we obtain that  for all sufficiently small $\e$,
\bea
 \ad\PP\big(R^\e\in\overline B(\ga,a), M^\e\in B(\phi_k,a_k),\Omega_\e\big)\\
 \aad\leq\EE\Big[\indi_{\{\la^\e\in\overline B(\ga-\phi_k,\widehat a(\phi_k))\}}\indi_{\Omega_\e}\,
 \PP\big(M^\e\in\overline B(\phi_k,\widetilde a(\phi_k))\,\big|\,\zeta\big)\Big]\\
 \aad\leq\exp\big\{-h^2(\e)\big[\chi_\eta\big(\CI_\zeta(\ga-\phi_k)\big)+\CI_w(\phi_k)-\eta\big]\big\}.
\eea
Hence
\beqa{local-one-ball}
 \ad \limsup_{\e\to0}h^{-2}(\e)\log
 \PP\big(R^\e\in\overline B(\ga,a), M^\e\in B(\phi_k,a_k),\Omega_\e\big) \\
 \aad \leq-\chi_\eta\big(\CI_\zeta(\ga-\phi_k)\big)
-\CI_w(\phi_k)+\eta.
\eeqa
The upper bound of the LDP of $M^\e$ implies $
 \limsup_{\e\to0}h^{-2}(\e)
 \log\PP(M^\e\notin \CO_{N_\kappa}^\kappa)\leq-\kappa$.
Therefore, using the finite covering and
$h^{-2}(\e)\log\PP(\Omega_\e^c)\to-\infty$, we obtain
\bea\ad
\limsup_{\e\to0}h^{-2}(\e)
 \log\PP(R^\e\in\overline B(\ga,a)) \\
 \aad \leq\max\big\{-\kappa ,-\min_{1\leq k\leq N_\kappa}
 \big[\chi_\eta\big(\CI_\zeta(\ga-\phi_k)\big)+\CI_w(\phi_k)-\eta\big]\big\}.
\eea
For $x,y\in[0,\infty]$, $ \chi_\eta(x)+y-\eta
 \geq(x+y-2\eta)\wedge(2\eta)^{-1}$.
Since $\CI_\zeta(\ga-\phi_k)+\CI_w(\phi_k)\geq \CI(\ga)$,
\bea\ad
 \limsup_{a\to0}\limsup_{\e\to0}h^{-2}(\e)
 \log\PP\big(R^\e\in\overline B(\ga,a)\big)\leq\max\{-\kappa ,
 -[(\CI(\ga)-2\eta)\wedge(2\eta)^{-1}]\}.
\eea
Letting $\kappa \to\infty$ and then $\eta\to0$ proves the local upper bound. The preceding two bounds prove the local LDP. Together with the
exponential tightness above, they give the full LDP.

It remains to characterize $\CI$ as $\CI_C$. Let $x\in\SD_A$ and set
$v:=\dot x-Ax$.  Since
$\CI_\zeta=\CI_{C_F}$ and $\CI_w=\CI_{C_G}$, decompositions
$x=x_F+x_G$ of finite cost correspond to decompositions
$v=v_F+v_G$ in $L^2([0,1];H)$.  Hence \eqref{LC1+C2}, applied
pointwise, gives $
\CI(x)\geq\int_0^1 L_C(t,v(t))dt=\CI_C(x)$.
For the reverse inequality, suppose that $\CI_C(x)<\infty$ and
define $
 T_t(w_F,w_G)
 :=C_F(t)^{1/2}w_F+C_G(t)^{1/2}w_G,
$ and $T_t T_t^*=C(t)$.
The measurability argument in Lemma~\ref{lem:MP-energy}, applied to
$T_t:H\oplus H\to H$, gives a strongly measurable minimal-norm
preimage $(w_F,w_G)\in L^2([0,1];H\oplus H)$ of $v$.  Set
$v_F=C_F^{1/2}w_F$ and $v_G=C_G^{1/2}w_G$.  Then
$v=v_F+v_G$, and \corref{cor:add} gives $
 \int_0^1[ L_{C_F}(t,v_F(t))+L_{C_G}(t,v_G(t))] dt
 =\int_0^1L_C(t,v(t))dt$.
The convolutions generated by $v_F$ and $v_G$ belong to $\SD_A$ and
sum to $x$, which proves $\CI(x)\leq\CI_C(x)$.  If $\CI(x)<\infty$,
both paths in a finite-cost decomposition belong to $\SD_A$. Hence
$x\in\SD_A$.  Thus $\CI(x)=+\infty$ outside $\SD_A$.
\lemref{lem:quad-lift} shows that $\CI_C$ is good.
\end{proof}

\subsection{\texorpdfstring{Exponential equivalence of $\eta^\e$ and $\widehat{\eta}^\e$}{Exponential equivalence of eta epsilon and eta-hat epsilon}}\label{sec:equiv-etah}

\begin{prop}\label{prop:eta-equivalence}
Under the assumptions of Theorem~\ref{thm:main}, $\eta^\e$ and
$\widehat\eta^\e$ are exponentially equivalent in $C^0([0,1];H)$ at
speed $h^{-2}(\e)$. Equivalently, \eqref{eta-equiv} holds.
\end{prop}
The proof is postponed to \appref{app:exp-equivalence}.

\subsection[Proof of the main result]{Proof of Theorem~\ref{thm:main}}
\begin{proof}
Lemma~\ref{lem:local-R} gives the LDP for $R^\e$
with rate $\CI_C=\CI_{C_F+C_G}$.
Proposition~\ref{prop:cp} then gives the rate of
$\widehat\eta^\e=\mathcal T(R^\e)$ as
$S(x)=\CI_C(\mathcal T^{-1}x)$.
The map $\mathcal T$ preserves $\SD_A$.  Indeed, the
convolution with
$\overline{\nabla F}(\bar X(\cdot))x(\cdot)\in L^2([0,1];H)$ belongs
to $\SD_A$ by the maximal $L^2$-regularity of the negative self-adjoint operator $A$, and
$
 \frac{d}{dt}(\mathcal T^{-1}x)(t)-A(\mathcal T^{-1}x)(t)
 =\dot x(t)-Ax(t)-\overline{\nabla F}(\bar X(t))x(t).
$
The same argument applied to $\mathcal T$ gives the converse
inclusion.  Hence
\bea
S(x)
\ad \!\! =\CI_{C_F+C_G}\!\big(
 x-\int_0^\cdot e^{(\cdot-s)A}\overline{\nabla F}(\bar X(s))x(s)ds\big) \\
\ad \!\! =\int_0^1L_{C_F+C_G}
 (t,\dot x(t)-Ax(t)-\overline{\nabla F}(\bar X(t))x(t))dt, \quad \forall \, x\in\SD_A,
\eea
with value $+\infty$ otherwise.
Lemma~\ref{lem:MP-energy} converts this identity into
\eqref{main-MP-rate}.  Proposition~\ref{prop:eta-equivalence} transfers the LDP and its good
rate from $\widehat\eta^\e$ to $\eta^\e$, which establishes
the theorem.
\end{proof}

\section{State-dependent regime-switching process on a finite state space}
\label{sec:feedback}
In this section, we consider the random environment to be a continuous-state-dependent (or, in short, state-dependent)
regime-switching process $\al^\e$ having a finite state space $\CM=\{1,2,\dots, m_0\}$. Let $\Upsilon\in\CL(H,\rr^{d_\Upsilon})$ and let
$Q(t,\Upsilon x):=(q_{ij}(t,\Upsilon x))_{i,j\leq m_0}$ be  the generator of $\al^\e$ for $t\in [0,1],x\in H$. The transition rates of $\al^\e$ are
\beqa{al-tran}
 \PP\big(\alpha^\e(t+\dl)=j\mid \al^\e(t)=i, X^\e(t)=x\big)
\!=\!\left\{\barray
 \e^{-1}q_{ij}(t,\Upsilon x)\dl+o(\dl),
 \; & j\ne i,  \\
 1+\e^{-1}q_{ii}(t, \Upsilon x)\dl+ o(\dl), \; &\text{otherwise.}
\earray\right.
\eeqa
Throughout this section, the slow component satisfies $X^\e(0)=x_0\in H$ and
\begin{equation}\label{srde-al}
 dX^\e(t)=[AX^\e(t)+F(X^\e(t),\alpha^\e(t))]dt
 +\sqrt\e G(X^\e(t),\alpha^\e(t))dW^{\mathsf{Q}}(t).
\end{equation}

For large and moderate deviations for finite-dimensional diffusions in fast-varying Markovian environments, we refer to Huang et al. \cite{HMS16} for homogeneous Markov chains and to \cite{HYZ11,HY14,QY22} for non-homogeneous Markov chains. For diffusions with fast-varying state-dependent switching processes, see \cite{BDG18,Qia25}. Comprehensive treatments of hybrid diffusions can be found in \cite{YZ10,YZ98} and the references therein.

\begin{ass}\label{ass:feedback}
In addition to Assumptions~\ref{ass:A} and~\ref{ass:fg} with $\SE=\CM$, we assume that for $i\ne j$, $q_{ij}$ is continuously differentiable in $t$ and twice continuously differentiable in $y$,
with all the corresponding derivatives bounded and uniformly
continuous.  There are constants $0<\underline{q}\leq \bar q<\infty$ such that
$\underline{q}\leq q_{ij}(t,y)\leq \bar q$ for $i\ne j$ and $q_{ii}=-\sum_{j\ne i}q_{ij}$.
\end{ass}
The uniform lower bound $\underline{q}>0$ on the off-diagonal implies uniform irreducibility. Let $\pi(t,y)=(\pi_1(t,y),\ldots,\pi_{m_0}(t,y))$ be the invariant distribution of $Q(t,y)$, written as a row vector, so that $\pi(t,y)Q(t,y)=0$. Let $\La(t,y)$ be the matrix obtained from $Q(t,y)^\top$
by replacing its last row with $(1,\ldots,1)$. Then
\beq{La-pi}
 \La(t,y)\pi(t,y)^\top=(0,\ldots,0,1)^\top.
\eeq
By the uniform bounds of $q_{ij}(\cdot,\cdot)$, the family $\{\La(t,y):t\in [0,1],y\in \rr^{d_\Upsilon}\}$ has compact closure. Every matrix in this closure corresponds to an irreducible generator and is therefore nonsingular.
Hence $\inf_{t\in[0,1],y\in\rr^{d_\Upsilon}}| \det \La(t,y)|>0$ and $ \sup_{t\in [0,1] ,y\in\rr^{d_\Upsilon}} \|\La^{-1}(t,y)\|<\infty$.
Differentiating \eqref{La-pi} with respect to $t$ and $y$ and applying
\assmref{ass:feedback} therefore yields
\beq{bd-pi}
\sup_{t\in [0,1],y\in\rr^{d_\Upsilon}} \big(\|\partial_t\pi(t,y)\|
+\|\nabla_y\pi(t,y)\|+\|\nabla_y^2\pi(t,y)\|\big)<\infty.
\eeq

\begin{lem}\label{lem:sd-wp}
For any fixed $\e>0$, equation \eqref{srde-al} has a unique mild solution.
\end{lem}
\begin{proof}
Let $\mathfrak p(dt,dz)$ be a Poisson random measure on
$[0,1]\times\mathbb R$ with intensity $dt\,dz$, independent of
$W^{\mathsf Q}$. For $i\ne j$, let  $\Delta_{ij}^\e(t,x)$ be the consecutive left-closed, right-open intervals of the real line, each having length
$\e^{-1}q_{ij}(t,\Upsilon x)$. Define
\[
\mathfrak h^\e(t,x,i,z)
=\sum_{j\ne i}(j-i)\indi_{\{z\in\Delta_{ij}^\e(t,x)\}}.
\]
Then with the partition $\{\Dl^\e_{ij}(t,x): i,j\in \CM\}$, \eqref{al-tran} is equivalent to
\bea\ad
d\al^\e(t)=\int_{\mathbb R}
\mathfrak h^\e(t,X^\e(t-),\alpha^\e(t-),z)\,\mathfrak p(dt,dz).
\eea
All intervals are in
$[0,m_0(m_0-1)\bar q/\e)$. Hence the number of candidate
jump times on $[0,1]$ is almost surely finite. Between
successive jump times, $\alpha^\e$ is constant, and
\thmref{thm:XbarX} gives a unique mild solution for the
corresponding fixed-regime stochastic evolution equation. At each jump time the regime is updated by $\mathfrak h^\e$ while the
continuous component is unchanged. Successive construction yields
a pathwise unique mild solution of \eqref{srde-al}; see \cite{Xi09}.
\end{proof}
Define the averaged drift
\beqa{sd-averaged-drift}\ad
 \bar F_{\rm sd}(t,x)=\sum_{i=1}^{m_0} \pi_i(t,\Upsilon x)F(x,i),
\eeqa
where the subscript ``sd'' denotes state-dependent,
and let $\bar X$ be the solution of
\beq{avg-al}
 \dot{\bar X}(t)=A\bar X(t)+\bar F_{\rm sd}(t,\bar X(t)),\quad \bar X(0)=x_0\in H.
\eeq
\begin{lem}\label{lem:avg-al}
Equation \eqref{avg-al} has a unique mild solution $\bar X\in C([0,1];H)$. Define $\bar y_t:=\Upsilon \bar X(t)$ and $\nu(t):=\pi(t,\bar y_t)$. Then there exists a constant $K>0$ such that $\|\bar X\|_\infty <\infty$ and
\bea
|\partial_t \bar X(t)|+\|\partial_t \nu(t)\|\leq K(1+t^{-1}),\quad 0<t\leq 1.
\eea
\end{lem}

\begin{proof}
Condition~(H1) and \assmref{ass:feedback} give $|\bar F_{\rm sd}(t,x)|_H\leq K_0+K_F'|x|_H$ and
\bea
|\bar F_{\rm sd}(t,x_1)-\bar F_{\rm sd}(t,x_2)|_H \leq [K_F'+K \|\Upsilon\|(K_0+K_F'|x_2|_H)] |x_1-x_2|.
\eea
Then a fixed point argument gives a unique mild solution of \eqref{avg-al}. Hence we have
\bea
|\bar X(t)|_H \leq |x_0|_H +K_0 t + K_F' \int_0^t |\bar X(s)|_H ds.
\eea
Applying Gronwall's inequality yields $\|\bar X\|_\infty\leq (|x_0|+K_0)e^{K_F'}<\infty$. Define $f_\al(t):=\bar F_{\rm sd}(t,\bar X(t))$. It then yields $\|f_\al\|_\infty <\infty$. Moreover,
\bea
|(-A)^{1/2}\bar X(t)|_H \ad\leq |(-A)^{1/2}e^{tA}x_0|_H + \int_0^t |(-A)^{1/2}e^{(t-s)A} f_\al(s)|_H ds  \\
\aad \leq K t^{-1/2}|x_0|_H + K \|f_\al\|_\infty \int_0^t (t-s)^{-1/2}ds \leq K t^{-1/2}, \quad 0<t \leq 1.
\eea
Hence $\|\bar X(t)\|_{1/2} \leq K t^{-1/2}$ for $t>0$. For $t/2<s<t$, the mild solution form gives
\bea
|\bar X(t)-\bar X(s)|_H \ad \leq |(e^{(t-s)A}-I)\bar X(s)|_H + |\int_{s}^t e^{(t-r)A}f_\al(r)dr|_H \\
\aad \leq K(t-s)^{1/2}\|\bar X(s)\|_{1/2}+ \|f_\al\|_\infty(t-s)  \leq Kt^{-1/2}(t-s)^{1/2}.
\eea
Similarly, noting that $\bar y_t=\Upsilon \bar X(t)$, we have for $t/2\leq s<t$,
\bea
|f_\al(t)-f_\al(s)|_H \leq K |\bar X(t)-\bar X(s)|_H + K(|t-s|+\|\bar y_t-\bar y_s\|) \leq K t^{-1/2}(t-s)^{1/2}.
\eea
Thus, $f_\al$ is locally $1/2$-H\"{o}lder continuous on $(0,1]$. For any fixed $t>0$, we also have
\bea
A \bar X(t) = Ae^{tA} x_0 + (e^{tA}-I)f_\al(t)+\int_0^t A e^{(t-s)A} [f_\al(s)-f_\al(t)]ds,
\eea
which implies
\bea
|A \bar X(t)|_H \ad \leq Kt^{-1}|x_0|_H + 2\|f_\al\|_\infty +K \int_0^{t/2} \frac{1}{t-s}ds + K t^{-1/2} \int_{t/2}^t (t-s)^{-1/2} ds \\
\ad  \leq K(1+t^{-1}).
\eea
The local H\"{o}lder continuity of $f_\al$ and the analytic semigroup now imply that $\bar X(t)\in \CD(A) $ and $\partial_t \bar X(t)=A\bar X(t)+f_\al(t)$. Consequently,  $|\partial_t \bar X(t)|_H \leq K(1+t^{-1})$. Finally, $\partial_t \nu(t)=\partial_t \pi(t,\bar y_t)+\nabla_y \pi(t,\bar y_t)[\Upsilon \partial_t \bar X(t)]$. Hence $\|\partial_t \nu(t)\|\leq K (1+t^{-1})$ for $t\in (0,1]$.
\end{proof}
The dependence of $\pi$ on the state enters the linearization through the bounded operator:
\beqa{def-Bpi}\ad
\Bpi(t)v :=\sum_i (\nabla_y\pi_i(t,\bar y_t)[\Upsilon v]) F(\bar X(t),i),\quad \forall\, v\in H.
\eeqa
The derivative of the averaged drift is therefore
\beqa{nab-f}\ad
 \nabla\bar F_{\rm sd}(t,\bar X(t))v
 =\sum_i\nu_i(t)\nabla F(\bar X(t),i)v+\Bpi(t)v,
 \eeqa
and let us define the ``centered'' drift \bea
b_\al(t,i):=F(\bar X(t),i) - \bar F_{\rm sd}(t,\bar X(t)).
\eea

Next, we follow steps in  Figure~\ref{fig:roadmap}, where the modified steps are marked by dashed boxes.

\subsection{\texorpdfstring{Decomposition of $\eta^\e$ and main results}{Decomposition of eta epsilon and main results}}\label{sec:sd-decomp} For fixed $R>0$, define the stopping time \bea
\tau_R^\e:=\inf\{t\in[0,1]: |\eta^\e(t)| >R\}\wedge 1.
\eea
For any $s\leq \tau_R^\e$, we then have $\|\Upsilon X^\e(s)- \bar y_s \|_{\rr^{d_\Upsilon}}\leq \|\Upsilon\|\sqe h(\e)R$. Define
\bea
\la_\al^\e(t)\ad :=\frac{1}{\sqe h(\e)} \int_0^t e^{(t-s)A} b_\al(s,\al^\e(s))ds \\
M_\al^\e(t) \ad :=h^{-1}(\e)\int_0^t e^{(t-s)A} G(\bar X(s),\al^\e(s))dW^{\mathsf{Q}}(s) \\
\CR_1^\e(t)\ad := \frac{1}{\sqe h(\e)}\int_0^t e^{(t-s)A} [F(X^\e(s),\al^\e(s))-F(\bar X(s),\al^\e(s))]ds \\
\CR_2^\e(t)\ad := h^{-1}(\e) \int_0^t e^{(t-s)A} [G(X^\e(s),\al^\e(s))-G(\bar X(s),\al^\e(s))] dW^{\mathsf{Q}}(s)
\eea
Similar to \eqref{exp-eta}, we can decompose $\eta^\e$ as
\beqa{decomp-eta}
\eta^\e(t)=\la_\al^\e(t)+M_\al^\e(t)+\CR_1^\e(t)+\CR_2^\e(t).
\eeqa
To proceed, we apply the Poisson equation technique to rewrite $\la_\al^\e$ as a stochastic convolution of a martingale plus some remainder terms.

 Let $U(t)=\{U(t,i)\}_{i\in \CM}$ solve the following Poisson equation:
\begin{equation}\label{Poi-al}
 -Q(t, \Upsilon \bar X(t)) U(t)=b_\al(t),\quad
 \sum_i\nu_i(t)U_i(t)=0.
\end{equation}
Here $b_\al(t)=(b_\al(t,1),\dots b_\al(t,m_0))^\top$ and $U(t)=(U(t,1),\dots, U(t,m_0))^\top$ are elements of $H^{m_0}$.
Let $\indi \nu(t)$ denote the $m_0\times m_0$ matrix whose rows all equal $\nu(t)$, where $\indi:=(1,\ldots,1)^\top$. The matrix $-Q(t,\bar y_t)+\indi \nu(t)$ is invertible and its inverse is uniformly bounded. Since $\nu(t)Q(t,\bar y_t)=0$ and $\nu(t)\indi =1$, we have $\nu(t)[-Q(t,\bar y_t)+ \indi \nu(t)]=\nu(t)$. Together with $\nu(t)b_\al(t)=0$, this shows that the unique centered solution of \eqref{Poi-al} is given by
\beq{sol-Poi}
U(t)=(-Q(t,\bar y_t)+\indi \nu(t))^{-1}b_\al(t).
\eeq
We write $U(t,i)$ or $U_i(t)$, depending on convenience. By \lemref{lem:avg-al} and (H1), we have
\bea\ad
\sup_{t\in [0,1]}\max_{i\in \CM}|b_\al(t,i)|_H<\infty,\quad \text{ and }\quad
\sup_{t\in [0,1]}\max_{i\in \CM}|U(t,i)|_H<\infty.
\eea
By the chain rule, we have
\bea\ad
\partial_t U(t)=(-Q(t,\bar y_t)+\indi \nu(t))^{-1}\Big[\partial_t b_\al(t)+ \frac{d}{dt}Q(t,\bar y_t)U(t)-\indi\partial_t\nu(t) U(t)\Big].
\eea
By \lemref{lem:avg-al}, the following estimate holds:
\beq{sd-U-regularity}
\max_{i\in\CM}|\partial_tU(t,i)|_H\le K(1+t^{-1}),
\quad 0<t\le1.
\eeq

Let $\wdt{\mathfrak p}(ds,dz):=\mathfrak p(ds,dz)-ds\,dz$ be the compensated Poisson random measure on $[0,1]\times \rr$. Define the $H$-valued jump martingale as
\beqa{sd-NU}\ad
N_\al^\e(t)=\frac{\sqe}{h(\e)}\int_0^t\!\!\int_{\rr}
\bigl[U(s,\al^\e(s-)+\mathfrak h^\e(s,X^\e(s-),\al^{\e}(s-),z))-U(s,\al^\e(s-))\bigr]
\,\widetilde{\mathfrak p}(ds,dz).
\eeqa
For $t\ge\e$, applying the generalized It\^o formula
\cite[Section~2.8]{YZ10} and \eqref{Poi-al}, we have
\beqa{Ito-U}\ad
 U(t,\al^\e(t))-U(\e,\al^\e(\e))
 =\int_\e^t\partial_sU(s,\al^\e(s))\,ds
 +\frac{h(\e)}{\sqe}[N_\al^\e(t)-N_\al^\e(\e)]
\\
\aad -\frac1\e\int_\e^t b_\al(s,\al^\e(s))\,ds+\frac1\e\int_\e^t\sum_{j=1}^{m_0}
 [q_{\al^\e(s)j}(s,\Upsilon X^\e(s))
 -q_{\al^\e(s)j}(s,\bar y_s)]U(s,j)\,ds.
\eeqa
Consequently, for $t\in[0,1]$, we have the following decomposition
\beq{dec-b}
\frac1{\sqe h(\e)}\int_0^t b_\al(s,\al^\e(s))\,ds
 =N_\al^\e(t)-N_\al^\e(t\wedge \e)+\CR_3^\e(t)+\CR_4^\e(t),
\eeq
where
\[
\CR_3^\e(t)=\frac1{\sqe h(\e)}\int_\e^{t\vee\e}
\sum_{j=1}^{m_0}
[q_{\al^\e(s)j}(s,\Upsilon X^\e(s))
-q_{\al^\e(s)j}(s,\bar y_s)]U(s,j)\,ds,
\]
and $\CR_4^\e(t)=(\sqe h(\e))^{-1}\int_0^t b_\al(s,\al^\e(s))ds$ if $0\leq t <\e$; otherwise, it is defined by
\bea
\CR_4^\e(t)\ad= \frac1{\sqe h(\e)}\int_0^\e b_\al(s,\al^\e(s))\,ds\\
\aad\quad +\frac{\sqe}{h(\e)}\Big[
\int_\e^t\partial_sU(s,\al^\e(s))\,ds -U(t,\al^\e(t))+U(\e,\al^\e(\e))\Big].
\eea

For a c\`adl\`ag process $Y$, denote its jump at time $t>0$ by
$\Delta Y(t):=Y(t)-Y(t-)$. Then, at every switching time $t>\e$, \eqref{sd-NU} and the definition of $\CR_4^\e$ give
\beq{Dl-N}
\Dl \CR_4^\e(t)=-\sqe h^{-1}(\e)[U(t,\al^\e(t))-U(t,\al^\e(t-))]=-\Dl N_\al^\e(t).
\eeq
Hence $N_\al^\e-N_\al^\e(\cdot\wedge \e)+\CR_4^\e$ is continuous. Define \bea
\la_{\al,0}^\e(t):=\int_0^te^{(t-s)A}\,d[N_\al^\e(s)-N_\al^\e(s\wedge \e)]
\eea
and let $\wdt\la_\al^\e(t):=\la_{\al,0}^\e(t)+
\int_0^te^{(t-s)A}\,d\CR_4^\e(s)$. Then we have
\beqa{wdt-la}
\la_\al^\e(t)=\wdt\la_\al^\e(t)
+\int_0^t e^{(t-s)A}\,d\CR_3^\e(s).
\eeqa
For $k\ge1$ and $t\ge\e$, integration by parts yields
\bea\ad
\Big|\int_\e^t e^{-\alpha_k(t-s)}\,d\langle\CR_4^\e(s),e_k\rangle_H\Big|
\le 2\sup_{\e\le s\le t}
\Big|\langle\CR_4^\e(s)-\CR_4^\e(\e),e_k\rangle_H\Big|.
\eea
For each $k\geq 1$, the definition of $\CR_4^\e$ gives
\bea\ad
\sup_{\e\leq s\leq 1}|\qv{\CR_4^\e(s)-\CR_4^\e(\e),e_k}_H| \leq \frac{2\sqe}{h(\e)}\sum_{i\in \CM}\Big(|\qv{U(\e,i),e_k}_H|+\int_\e^1 |\qv{\partial_s U(s,i),e_k}_H| ds\Big)
\eea
Consequently, Minkowski's inequality and the finiteness of $\CM$ yield
\bea\ad
\Big(\sum_{k\ge1}\sup_{\e \leq s \leq 1}|\qv{\CR_4^\e(s)-\CR_4^\e(\e),e_k}_H|^2 \Big)^{1/2}  \leq \frac{K\sqe}{h(\e)}\Big(1+\int_\e^1 \max_{i\in \CM} |\partial_s U(s,i)|_H ds \Big)
\eea
Together with the preceding integration by parts estimate and
\eqref{sd-U-regularity}, we obtain
\beqa{sd-R4-conv}\disp
\left\|\int_0^\cdot e^{(\cdot-s)A}\,d\CR_4^\e(s)\right\|_\infty
\ad \le\frac{K\e}{\sqe h(\e)}
+\frac{K\sqe}{h(\e)}\left(1+
\int_\e^1\max_i|\partial_sU(s,i)|_H\,ds\right)\\
\aad
\le\frac{K\sqe}{h(\e)}(1+|\log\e|)\to0.
\eeqa

We now identify the covariance fields.  For every $n\in \NN_+$, define $N_{\al,n}^\e:=\Pi_n N_\al^\e$ and $Z_{\al,n}^\e(t):=h^{-1}(\e)\int_0^t\Pi_nG(\bar X(s),\al^\e(s))dW^{\mathsf{Q}}(s)$. Their predictable quadratic variations are
\beqa{brak}
h^2(\e)\qv{N_{\al,n}^\e}_t \ad =\int_0^{t} \sum_{j=1}^{m_0} q_{\al^\e(s) j}(s,\Upsilon X^\e(s))[\Pi_n U(s,j)-\Pi_n U(s,\al^\e(s))]^{\otimes 2} ds, \\
h^2(\e)\qv{Z_{\al,n}^\e}_t \ad = \int_0^t \Pi_n G(\bar X(s),\al^\e(s)) \mathsf{Q}\mathsf{Q}^* G(\bar X(s),\al^\e(s))^*\Pi_n ds.
\eeqa
Averaging the above two integrands with respect to the invariant measure $\nu(t)$ and replacing $X^\e$ by $\bar X$ lead to the covariance operators
\begin{align}
 C_{F,{\rm sd}}(t)
 &=\sum_{i,j}\nu_i(t)q_{ij}(t,\Upsilon \bar X(t))
 [U_j(t)-U_i(t)]\otimes[U_j(t)-U_i(t)],               \label{sd-CF}\\
 C_{G,{\rm sd}}(t)
 &=\sum_i\nu_i(t)G(\bar X(t),i)\mathsf{Q} \mathsf{Q}^*G(\bar X(t),i)^*.             \label{sd-CG}
\end{align}
By the bounds on $U$ and $q_{ij}$, \eqref{G-op}, and \lemref{lem:avg-al}, both $C_{F,{\rm sd}}(t)$ and $C_{G,{\rm sd}}(t)$ are nonnegative, self-adjoint, and bounded in $\CL(H)$ uniformly in $t\in[0,1]$.
Define $R_\al^\e:=\wdt \la_\al^\e+ M_\al^\e$ and let $\wdh\eta_\al^\e$ be the solution of
\beqa{wdh-eta-al}\ad
\wdh \eta_\al^\e(t)=R_\al^\e(t)+\int_0^t e^{(t-s)A} \nabla \bar F_{\rm sd}(s,\bar X(s)) \wdh \eta_\al^\e(s)ds
\eeqa
The proof reduces to proving the LDP for $R_\al^\e$ and the exponential equivalence of $\eta^\e$ and $\wdh \eta_\al^\e$.

\begin{thm}\label{thm:feedback-main}
Under Assumptions~\ref{ass:A},~\ref{ass:fg}, and~\ref{ass:feedback}, let $h(\e)$ satisfy \eqref{dev-s}.  Then
$\eta^\e$ satisfies an LDP in $C^0([0,1];H)$ with speed
$h^{-2}(\e)$ and good rate
\begin{equation}\label{sd-main-rate}
 S_{\rm sd}(x)=
 \begin{cases}
 \int_0^1L_{C_{F,{\rm sd}}+C_{G,{\rm sd}}}
 (t,\dot x(t)-Ax(t)-\nabla\bar F_{\rm sd}(t,\bar X(t))x(t))dt,
       &x\in\SD_A,\\
 +\infty,&\text{otherwise}.
 \end{cases}
\end{equation}
Equivalently, with $C_{\rm sd}=C_{F,{\rm sd}}+C_{G,{\rm sd}}$ and
$v_x(t)=\dot x(t)-Ax(t)-\nabla\bar F_{\rm sd}(t,\bar X(t))x(t)$,
\beqa{sd-MP-rate}
 S_{\rm sd}(x)=\frac12\int_0^1
 \big|\big[C_{\rm sd}(t)^{1/2}\big]^\dagger
 v_x(t)\big|_H^2dt
\eeqa
when $x\in\SD_A$, $v_x(t)\in
\Ran\bigl(C_{\rm sd}(t)^{1/2}\bigr)$ a.e., and the
integral is finite. Otherwise it is infinite.
\end{thm}

\begin{cor}\label{cor:feedback-independent}
Suppose that $Q(t,y)=Q(t)$ is independent of $y$.  Then
$\Bpi(t)=0$ and
$\nabla\bar F_{\rm sd}(t,\bar X(t))=\sum_i\nu_i(t)\nabla F(\bar X(t),i)$.
Consequently, the independent time-inhomogeneous finite-state chain
satisfies the MDP of Theorem~\ref{thm:feedback-main}, with rate
\eqref{sd-main-rate} and covariances
$C_{F,{\rm sd}},C_{G,{\rm sd}}$ evaluated using the invariant vector
$\nu(t)$ of $Q(t)$.
\end{cor}

\subsection{The LDP for the sum of two fluctuations} We first record an exponential inequality for jump martingales driven by $\al^\e$.

\begin{lem}\label{lem:jump-exp}
Let $\fg$ be
deterministic with values in a separable Hilbert space, with $\fg(s,i,i)=0$ and $|\fg|\le \kappa_{\fg}$, and define
\[
\fM_{\fg}(t)=\int_0^t \int_{\rr}
\fg\bigl(s,\al^\e(s-),
\al^\e(s-)+\mathfrak h^\e(s,X^\e(s-),\al^\e(s-),z)\bigr)
\,\widetilde{\mathfrak p}(ds,dz).
\]
If
$
V_{\fg}(t):=\frac1\e\int_0^t\sum_{j=1}^{m_0}
q_{\al^\e(s-)j}(s,\Upsilon X^\e(s-))
|\fg(s,\al^\e(s-),j)|^2\,ds\le \kappa_V
$ for $t\le1$,
then
\begin{equation}\label{fb-jump-exp}
\PP\Bigl(\sup_{t\le1}|\fM_{\fg}(t)|>u\Bigr)
\le2\exp\Bigl(-\frac{u^2}{2e\kappa_V+2\kappa_{\fg} u}\Bigr),\qquad u>0.
\end{equation}
\end{lem}
\begin{proof}
Fix $0<\varsigma\le1/\kappa_{\fg}$ and set
$\aleph(x)=\cosh(\varsigma|x|)$. As in
\cite[Theorem~3.2]{Pin94} and \cite{Pin99},
\bea\ad
\aleph(x+y)-\aleph(x)-\langle\nabla \aleph(x),y\rangle
\le\frac e2\varsigma^2|y|^2 \aleph(x),\qquad |y|\le \kappa_{\fg}.
\eea
The process $\fM_{\fg}$ has finitely many jumps on $[0,1]$ and is
absolutely continuous between them.  The jump measure of $\al^\e$ has
compensator
$\e^{-1}q_{\al^\e(s-)j}(s,\Upsilon X^\e(s-))\,ds$ on
$\{j\ne\al^\e(s-)\}$; see \cite[Theorem~II.1.8]{JS03}.  Therefore
$\exp\{-\frac e2\varsigma^2 V_{\fg}(t)\}\aleph(\fM_{\fg}(t))$ is a nonnegative
supermartingale starting from $1$.  Stopping it at the first time that
$|\fM_{\fg}|$ reaches $u$, as at the end of the proof of
\cite[Theorem~3.2]{Pin94} (see also \cite{Pin99}), gives
\[
\PP\Bigl(\sup_{t\le1}|\fM_{\fg}(t)|>u\Bigr)
\le2\exp\Bigl(-\varsigma u+\frac e2\varsigma^2\kappa_V\Bigr).
\]
Taking $\varsigma=\min\{u/(e\kappa_V),1/\kappa_{\fg}\}$ proves \eqref{fb-jump-exp}.
\end{proof}

The next lemma plays the same role as \lemref{lem:averaging} for the state-dependent case, where no independence of the environment and driving noise is assumed.

\begin{lem}\label{lem:sd-average}
Let $c:[0,1]\times \CM \to \rr$ be a bounded function such that $c(\cdot,i)$ is continuous for every $i\in \CM$ and $\sum_i \nu_i(t)c(t,i)=0$.
For every $R,r>0$,
\beqa{sd-average-stopped}\ad
 \limsup_{\e\to0}h^{-2}(\e)\log\PP\Big(
 \sup_{t\leq1}\Big|\int_0^{t\wedge\tau_R^\e}
 c(s,\al^\e(s))ds\Big|>r\Big)=-\infty.
\eeqa
The same assertion holds coordinatewise for bounded
$H_n$-valued functions.  Moreover, let
$\ell_v:[0,1]\times\CM\to H$, $v\in C^0([0,1];H)$, be centered in the
sense that $\sum_i\nu_i(t)\ell_v(t,i)=0$,
and suppose that on every compact $\CO \subset C^0([0,1];H)$, these
functions are uniformly bounded, uniformly continuous in $t$, and
Lipschitz in $v$ with respect to $\|\cdot\|_\infty$.  Then, for every
compact $\CO \subset C^0([0,1];H)$ and $R,r>0$,
\beq{sd-random-average}
 \limsup_{\e\to0}h^{-2}(\e)\log\PP\Big(
 \eta^\e\in \CO,\
 \sup_{t\leq1}\Big|\int_0^{t\wedge\tau_R^\e}
 e^{(t-s)A}\ell_{\eta^\e}(s,\alpha^\e(s))ds\Big|_H>r\Big)=-\infty.
\eeq
\end{lem}
The proof is given in \appref{app:lem43}.

\begin{prop}\label{prop:sd-driving}
The process $R_\al^\e=\wdt\la_\al^\e+M_\al^\e$ satisfies an LDP in
$C^0([0,1];H)$ with the uniform topology, speed $h^{-2}(\e)$, and good
rate function $\CI_{C_{F,{\rm sd}}+C_{G,{\rm sd}}}$.
\end{prop}

\begin{proof}
      \textbf{Step 1. Projected predictable quadratic variations}.
For $k,\iota\leq n, i\in\CM$, define
\begin{align*}
 \ell_{\iota k}^{1,n}(s,i)
 &:=\sum_{j=1}^{m_0} q_{ij}(s,\bar y_s)
 \langle U_j(s)-U_i(s),e_k\rangle_H
 \langle U_j(s)-U_i(s),e_\iota\rangle_H
 -\langle C_{F,{\rm sd}}(s)e_k,e_\iota\rangle_H \\
  \ell_{\iota k}^{2,n}(s,i)
 &:=\langle \mathsf{Q}^*G(\bar X(s),i)^*e_k,
              \mathsf{Q}^*G(\bar X(s),i)^*e_\iota \rangle_H
 -\langle C_{G,{\rm sd}}(s)e_k,e_\iota \rangle_H.
\end{align*}
Both functions are continuous and centered under $\nu$. By \eqref{brak}, for each fixed $R>0$, we have
\bea\ad
\sup_{t\leq \tau_R^\e}\Big\|h^2(\e)\qv{N_{\al,n}^\e}_t - \int_0^t \Pi_n C_{F,\rm sd}(s)\Pi_n ds \Big\| \\
\aad\quad  \leq K_n \Big(\sqe h(\e)R + \e+\max_{k,\iota \leq n,t\leq 1} \Big|\int_0^{t\wedge \tau_R^\e} \ell_{\iota k}^{1,n}(s,\al^\e(s)) ds \Big|\Big) \\
\ad
\sup_{t\leq \tau_R^\e} \Big\|h^2(\e)\qv{Z_{\al,n}^\e}_t-\int_0^t \Pi_n C_{G,\rm sd}(s)\Pi_n ds \Big\| \leq K_n \max_{\iota,k\leq n}\sup_{t\leq 1} \Big|\int_0^{t\wedge \tau_R^\e} \ell_{\iota k}^{2,n}(s,\al^\e(s))ds \Big|.
\eea
Indeed, before $\tau_R^\e$, the difference between $Q(s,\Upsilon X^\e(s))$ and $Q(s,\bar y_s)=Q(s,\Upsilon \bar X(s))$ is bounded by $K \sqe h(\e)R$. Lemma~\ref{lem:sd-average} then gives, for every $r>0$,
\begin{align}\limsup_{\e\to0}h^{-2}(\e)\log\PP\Big(
 \sup_{t\leq\tau_R^\e}\big\|h^2(\e)
 \langle N_{\al,n}^\e\rangle_t
 -\int_0^t\Pi_n C_{F,{\rm sd}}(s)\Pi_nds\big\|>r\Big)=-\infty  \label{jump-CF}\\
\limsup_{\e\to0}h^{-2}(\e)\log\PP\Big(
 \sup_{t\leq\tau_R^\e}\big\|h^2(\e)
 \langle Z_{\al,n}^\e\rangle_t
 -\int_0^t\Pi_n C_{G,{\rm sd}}(s)\Pi_nds \big\|>r\Big)
 =-\infty.                                            \label{sd-W-qv-limit}
\end{align}
\smallskip

\noindent\textbf{Step 2. Removal of the stopping time}.
Define $K_U:=1+2\sup_{t\leq 1,i\in \CM} |U(t,i)|$ and $\kappa_q:=m_0(m_0-1) \bar q$. By \eqref{Dl-N}, we have
\beq{est-Dl-N}
|\Dl N_\al^\e(t)|_H \leq \frac{K\sqe}{h(\e)}, \quad \frac{1}{\e}\int_0^1 \sum_{j\neq \al^\e(s)} q_{\al^\e(s) j} (s,\Upsilon X^\e(s))ds \leq \frac{\kappa_q}{\e}.
\eeq
By \assmref{ass:A}, $\|e^{tA}\|_{\CL(H)}\leq 1$ for every $t\geq 0$.
Hence by \cite{HS08}, it  has a unitary dilation $e^{tA}=\mathfrak{J}^* \mathfrak{U}(t)\mathfrak{J}$, where $\mathfrak{J}$ is an isometry and $\mathfrak{U}$ is a unitary group on a larger Hilbert space. Therefore, for every $H$-valued martingale $N$,
\beq{dila}
\sup_{t\leq 1} \Big|\int_0^t e^{(t-s)A} dN(s)\Big|_H \leq \sup_{t\leq 1} \Big|\int_0^t \mathfrak{U}(-s) \mathfrak{J} dN(s)\Big|.
\eeq
For $N\cd=N_\al^\e-N_\al^\e(\cdot\wedge \e)$, the martingale on the right of \eqref{dila} is
\beqa{sd-jump-mart}\ad
\fM(t)=\int_0^t \int_{\rr}
\fg_{>\e}(s,\al^\e(s-),\al^\e(s-)+\mathfrak h^\e(s,X^\e(s-),\al^\e(s-),z))\,
\widetilde{\mathfrak p}(ds,dz),
\eeqa
where $\fg_{>\e}(s,i,j):=\indi_{\{s>\e\}}\sqe h^{-1}(\e)\fU(-s)\mathfrak{J}[U(s,j)-U(s,i)]$ for $i,j\in \CM$.

For $\fg_{>\e}$, the isometry of $\fU(-s)\mathfrak{J}$ gives
$\kappa_{\fg}=K_U\sqe/h(\e)$ and $V_{\fg}(1)\le \kappa_V=K_U^2\kappa_q/h^2(\e)$.
Thus \eqref{dila} and \lemref{lem:jump-exp} yield
\begin{equation}\label{fb-ZF-tail}
\PP(\|\la_{\al,0}^\e\|_\infty>r)
\le2\exp\left\{-\frac{r^2h^2(\e)}
{2eK_U^2\kappa_q+2K_Ur\sqe h(\e)}\right\}.
\end{equation}
For the integral on $[0,\e]$, apply \eqref{fb-jump-exp} with
$
\fg_{\le\e}(s,i,j)=\indi_{\{s\le\e\}} \sqe h^{-1}(\e)
[U(s,j)-U(s,i)].
$
Here $\kappa_{\fg}=K_U\sqe/h(\e)$ and
$V_{\fg}(1)\le \e K_U^2\kappa_q/h^2(\e)$, so we have
\begin{equation}\label{fb-initial-tail}
\PP\bigl(\|N_\al^\e(\cdot\wedge\e)\|_\infty>r\bigr)
\le2\exp\left\{-\frac{r^2h^2(\e)}
{2eK_U^2\kappa_q\e+2K_Ur\sqe h(\e)}\right\}.
\end{equation}
By \eqref{fb-ZF-tail} and \eqref{sd-R4-conv},
\[
 \limsup_{\e\to0}h^{-2}(\e)
 \log\PP(\|\wdt\la_\al^\e\|_\infty>r)
 \le-\frac{r^2}{8eK_U^2\kappa_q},
 \qquad r>0.
\]
Moreover, applying Lemma~\ref{lem:conv-exp} and
\remref{rem:G-bounds} to
$\Phi^\e(s)=h^{-1}(\e)G(\bar X(s),\alpha^\e(s))\mathsf{Q}$ gives
\begin{equation}\label{sd-YG-tail}
 \limsup_{\e\to0}h^{-2}(\e)
 \log\PP(\|M_\al^\e\|_\infty>r)
 \le-\kappa_G r^2.
\end{equation}
Here $\kappa_G:=K_2^{-1}(K_G+K_G'\|\bar X\|_\infty)^{-2}$.

To proceed, define the stopped error for diffusion:
\bea\ad
\CR_{2,R}^\e(t):=h^{-1}(\e) \int_0^t e^{(t-s)A} \indi_{\{s\leq \tau_R^\e\}} [G(X^\e(s),\al^\e(s))-G(\bar X(s), \al^\e(s))]dW^{\mathsf{Q}}(s).
\eea
On the event $\{s\leq \tau_R^\e\}$, (H2) implies that the squared Hilbert--Schmidt norm of the above integrand is bounded by $K\e R^2(t-s)^{-\vartheta}$. Hence, \lemref{lem:conv-exp} gives, for every $r>0$,
\bea
\PP(\|\CR_{2,R}^\e\|_\infty > r) \leq K_1 \exp(-r^2/(K_2K \e R^2)).
\eea

Define $u_R^\e(t):=\sup_{s\leq t\wedge \tau_R^\e} |\eta^\e(s)|_H$.
By \eqref{decomp-eta}, \eqref{wdt-la}, the Lipschitz continuity of $F$ in $x$ and of $q_{ij}$ in $y$, and the boundedness of $U$, the terms $\CR_1^\e$ and $\int_0^\cdot e^{(\cdot-s)A}d\CR_3^\e(s)$ are bounded by $K\int_0^{t\wedge \tau_R^\e} u_R^\e(s)ds$. Hence
\bea\ad
u_R^\e(t)\leq \|\wdt \la_\al^\e\|_\infty + \|M_\al^\e\|_\infty +  \|\CR_{2,R}^\e\|_\infty
 +K \int_0^{t\wedge \tau_R^\e} u_R^\e(s)ds.
\eea
Thus, Gronwall's inequality implies that for sufficiently small $\e$, there exists a constant $\kappa>0$ (independent of $\e$ and $R$) such that
\beq{stop-sub}
\{\tau_R^\e <1\} \subset \{\|\wdt \la_\al^\e\|_\infty > \kappa R\} \cup \{\|M_\al^\e\|_\infty > \kappa R\} \cup \{\|\CR_{2,R}^\e\|_\infty > \kappa R\}.
\eeq
Consequently, combining \eqref{fb-ZF-tail}--\eqref{stop-sub}, we have
\beq{lim-tau}
\lim_{R\to \infty} \limsup_{\e\to 0} h^{-2}(\e) \log\PP(\tau_R^\e<1)=-\infty.
\eeq
\smallskip

\noindent\textbf{Step 3. The projected LDP and passage to $H$.} Fix $n\ge1$ and set $V_n^\e:=N_{\al,n}^\e+Z_{\al,n}^\e$ and
$C_{{\rm sd},n}=\Pi_n C_{\rm sd}\Pi_n|_{H_n}$, where
$C_{\rm sd}:=C_{F,{\rm sd}}+C_{G,{\rm sd}}$.
The predictable quadratic covariation $\langle N_{\al,n}^\e,Z_{\al,n}^\e\rangle$ vanishes because $N_{\al,n}^\e$ is purely discontinuous
and $Z_{\al,n}^\e$ is continuous. Hence
$\qv{V_n^\e}=\qv{N_{\al,n}^\e}+\qv{Z_{\al,n}^\e}$.

For every $r>0$, the triangle inequality and a union bound show that
the probability below is at most $\PP(\tau_R^\e<1)$ plus the
probabilities in \eqref{jump-CF} and \eqref{sd-W-qv-limit} with
threshold $r/2$. Letting $\e\to0$ and then $R\to\infty$, and using
\eqref{lim-tau}, gives
\bea\ad
\limsup_{\e\to0}h^{-2}(\e)\log\PP\Big(
\sup_{t\le1}\Big\|h^2(\e)\qv{V_n^\e}_t
-\int_0^t C_{{\rm sd},n}(s)\,ds\Big\|>r\Big)=-\infty.
\eea
Since $Z_{\al,n}^\e$ is continuous,
$
\Delta V_n^\e(t)=\Pi_n\Delta N_\al^\e(t).
$
To verify the scaling in \cite[Corollary~6.7]{Puh94}, take an
arbitrary sequence $\e_k\to0$ and set
$m_k:=\lfloor h^2(\e_k)\rfloor$. The compensator of $\mathfrak p(ds,dz)$ is $ds\,dz$. Hence the
predictable compensator of the jump measure of $V_n^{\e_k}$ is
\[
\mathfrak n_k(ds,dx)
=\frac{ds}{\e_k}
\sum_{j\ne\al^{\e_k}(s-)}\!\!\!
q_{\al^{\e_k}(s-)j}
\bigl(s,\Upsilon X^{\e_k}(s-)\bigr)\times
\dl_{
\sqrt{\e_k}h^{-1}(\e_k)
\Pi_n[
U(s,j)-U(s,\al^{\e_k}(s-))]
}(dx),
\]
where $\delta_y$ denotes the Dirac measure at $y\in H_n$.

Set $
a_k:=K_U\sqrt{\e_k}h^{-1}(\e_k)$.
Since $\Pi_n$ is a contraction and
$|U(s,j)-U(s,i)|_H\le K_U$, we have
$
|x|_H\le a_k$  for $\mathfrak n_k$-a.e. $(s,x)$.
Moreover,
$
m_ka_k
\le K_U\sqrt{\e_k}\,h(\e_k)\to0.
$
Consequently, for every $\delta>0$ and all sufficiently large $k$, $
m_k\int_0^1\int_{H_n}
|x|_H^2
\indi_{\{m_k|x|_H>\delta\}}
\,\mathfrak n_k(ds,dx)=0$.
The explicit expression for $\mathfrak n_k$ also gives
\[
\begin{aligned}
\int_0^1\int_{H_n}
|x|_H^2\,\mathfrak n_k(ds,dx)
&=
\frac1{\e_k}\int_0^1
\sum_{j\ne\al^{\e_k}(s-)}
q_{\al^{\e_k}(s-)j}
\bigl(s,\Upsilon X^{\e_k}(s-)\bigr)\\
&\quad \times
\Big|
\frac{\sqrt{\e_k}}{h(\e_k)}
\Pi_n\bigl[
U(s,j)-U(s,\al^{\e_k}(s-))
\bigr]
\Big|_H^2ds \le
\frac{K_U^2\kappa_q}{h^2(\e_k)}.
\end{aligned}
\]
Therefore, we obtain
$
m_k\int_0^1\int_{H_n}
|x|_H^2\,\mathfrak n_k(ds,dx)
\le m_k h^{-2}(\e_k) K_U^2\kappa_q
\le K_U^2\kappa_q$.
For $y\in\mathbb R$, Taylor's formula gives
$
|e^y-1-y-\frac12y^2|
\le\frac16e^{|y|}|y|^3.
$
Fix $R'<\infty$. If $\psi\in H_n$, $|\psi|_H\le R'$, and
$x\in\operatorname{supp}\mathfrak n_k(s,\cdot)$, then
$
|m_k\langle\psi,x\rangle_H |
\le R' m_ka_k\to0.
$
Hence, for some $K_{R'}>0$,
\bea
&\disp\sup_{\substack{t\le1\\
\psi\in H_n,\ |\psi|_H\le R'}}
\left|
\frac1{m_k}\int_0^t\int_{H_n}
\big[
e^{m_k\langle\psi,x\rangle_H}
-1-m_k\langle\psi,x\rangle_H
-\frac{m_k^2}{2}\langle\psi,x\rangle_H^2
\big]\mathfrak n_k(ds,dx)
\right|\\
&\qquad\le
K_{R'}\,m_k a_k
\big(
m_k\int_0^1\int_{H_n}
|x|_H^2\,\mathfrak n_k(ds,dx)
\big)
\to0.
\eea

Together with the convergence above, these estimates verify
the Cram\'er and exponential Lindeberg conditions in
\cite[Corollary~6.7]{Puh94}. It follows that $V_n^\e$ satisfies an LDP in $D^0([0,1];H_n)$, also for the uniform topology by \cite[Remark~2.4]{Puh94}, with good rate function
$
\int_0^1L_{C_{{\rm sd},n}}(t,\dot z(t))\,dt
$
for $z\in W^{1,2}([0,1];H_n)$ with $z(0)=0$, and with value
$+\infty$ otherwise. The $W^{1,2}$ condition for finite cost follows
from \eqref{low-LC}. Since the argument applies to every sequence
$\e_k\to0$, the conclusion holds for the full family as $\e\to0$.

By \eqref{Dl-N},
$
N_\al^\e-N_\al^\e(\cdot\wedge\e)+\CR_4^\e
$
is continuous. Since $Z_{\al,n}^\e$ is continuous, it follows that
$
V_n^\e+\Pi_n\CR_4^\e
-\Pi_nN_\al^\e(\cdot\wedge\e)
\in C^0([0,1];H_n).
$
By the definitions of $\wdt\la_\al^\e$, $M_\al^\e$, and
$R_\al^\e$, we have
$
\Pi_nR_\al^\e(t)
=
\int_0^t e^{(t-s)A_n}\,
d\big[
V_n^\e(s)+\Pi_n\CR_4^\e(s)
-\Pi_n N_\al^\e(s\wedge\e)
\big]$. Let $\Phi_n$ be the map in \eqref{def-Phi-n}. Integration by parts gives
\bea
\Pi_nR_\al^\e(t)=V_n^\e(t)+\Pi_n\CR_4^\e(t)
-\Pi_nN_\al^\e(t\wedge\e)
+\int_0^tA_n\Pi_nR_\al^\e(s)\,ds
\eea
and
\beq{fb-Volterra-inverse}
\Phi_n^{-1}(\Pi_n R_\al^\e)
=V_n^\e+\Pi_n\CR_4^\e-\Pi_nN_\al^\e(\cdot\wedge\e).
\eeq
The definition of $\CR_4^\e$ gives
$
\|\Pi_n\CR_4^\e\|_\infty
\le K\sqe h^{-1}(\e)(1+|\log\e|)\longrightarrow0.
$
Combining this deterministic estimate with \eqref{fb-initial-tail}
shows that, for every $r>0$,
\beq{fb-Volterra-error}
\limsup_{\e\to0}h^{-2}(\e)\log\PP(
\|\Pi_n\CR_4^\e-\Pi_nN_\al^\e(\cdot\wedge\e)\|_\infty>r
)=-\infty.
\eeq

By \eqref{fb-Volterra-inverse} and \eqref{fb-Volterra-error},
$\Phi_n^{-1}(\Pi_nR_\al^\e)$ and $V_n^\e$ are exponentially
equivalent in $D^0([0,1];H_n)$ with the uniform metric.
Hence, \cite[Theorem~4.2.13]{DZ09} transfers the LDP of $V_n^\e$
to $\Phi_n^{-1}(\Pi_nR_\al^\e)$. The resulting rate function is
finite only on $W^{1,2}([0,1];H_n)\subset C^0([0,1];H_n)$, while
$\Phi_n^{-1}(\Pi_nR_\al^\e)$ has continuous paths. The inclusion
principle \cite[Lemma~4.1.5 (b)]{DZ09} therefore gives the same LDP
in $C^0([0,1];H_n)$ with the uniform topology. Since $\Phi_n$ is continuous, the contraction
principle \cite[Theorem~4.2.1]{DZ09} gives the LDP for
$\Pi_nR_\al^\e$. By the properties of $\Phi_n$ after \eqref{def-Phi-n}, the rate function
is $\CI_n$ in \lemref{lem:quad-lift} with covariance
$C_{{\rm sd},n}$.

To control the tail, set $K_{U,n}=2\sup_{t\le1,i\in\CM}|(I-\Pi_n)U(t,i)|_H$.
Each $U(\cdot,i)$ is continuous on $[0,1]$, so $\{U(t,i):t\le1,\,i\in\CM\}$ is compact in $H$ and $K_{U,n}\to0$ as $n\to\infty$.
Apply \eqref{dila} to the increment
$N_\al^\e-N_\al^\e(\cdot\wedge\e)$ and use
\eqref{fb-jump-exp} with
$(I-\Pi_n)[U(s,j)-U(s,i)]$ in place of $U(s,j)-U(s,i)$.
Taking $\kappa_{\fg}=K_{U,n}\sqe/h(\e)$ and
$\kappa_V=K_{U,n}^2\kappa_q/h^2(\e)$ gives
\[
\PP\bigl(\|(I-\Pi_n)\la_{\al,0}^\e\|_\infty>r\bigr)
\le2\exp\Big(-\frac{r^2h^2(\e)}
{2eK_{U,n}^2\kappa_q+2K_{U,n}r\sqe h(\e)}\Big).
\]
Together with \eqref{sd-R4-conv}, this gives, for every $r>0$,
\begin{equation}\label{I-Pi-la}
 \lim_{n\to\infty}\limsup_{\e\to0}h^{-2}(\e)
 \log\PP\big(\|(I-\Pi_n)\wdt\la_\al^\e\|_\infty>r\big)
 =-\infty.
\end{equation}
For the Wiener part, fix $0<\delta<1-\vartheta$. By
\eqref{abstract-Gtail} and $2^\vartheta e^{-x}\leq K_\delta x^{-\delta}$, we have
\bea
 \|e^{uA}(I-\Pi_n)G(\bar X(s),\al^\e(s))\mathsf{Q}\|_{\CL_2(H)}^2
 \leq K_\delta (K_G+K_G'\|\bar X\|_\infty)^2\alpha_{n+1}^{-\delta}
 u^{-(\vartheta+\delta)},\, 0<u\leq1.
\eea
Applying \lemref{lem:conv-exp} to the predictable integrand
$h^{-1}(\e)(I-\Pi_n)G(\bar X(s),\al^\e(s))\mathsf{Q}$ gives
\beq{sd-YG-Galerkin}
 \limsup_{\e\to0}h^{-2}(\e)
 \log\PP(\|(I-\Pi_n)M_\al^\e\|_\infty>r)
 \leq -\kappa_\delta r^2\alpha_{n+1}^{\delta},\qquad r>0.
\eeq
 Here $\kappa_\delta:=[K_2K_\delta(K_G+K_G'\|\bar X\|_\infty)^2]^{-1}$.
A union bound using \eqref{I-Pi-la} and \eqref{sd-YG-Galerkin} at radius
$r/2$ gives \eqref{abstract-lift-tail} for $Y^\e=R_\al^\e$. The projected
LDP proved above verifies the other hypothesis of
\corref{thm:abstract-lift}. The corollary gives the claimed LDP with
rate $\CI_{C_{F,{\rm sd}}+C_{G,{\rm sd}}}$.
\end{proof}

\subsection{\texorpdfstring{Exponential equivalence of $\eta^\e$ and $\wdh \eta_\al^\e$}{Exponential equivalence of eta epsilon and eta-hat alpha epsilon}}
\begin{prop}\label{prop:equi-eta-al}
Under the assumptions of \thmref{thm:feedback-main}, for every $r>0$,
\beq{equi-eta-al}
\limsup_{\e\to 0} h^{-2}(\e)\log\PP(\|\eta^\e-\wdh \eta_\al^\e\|_\infty>r)=-\infty.
\eeq
\end{prop}
The proof is postponed to \appref{app:prop-equi-eta-al}.
\subsection[Proof of the main result]{Proof of \thmref{thm:feedback-main}}
\begin{proof}
\propref{prop:sd-driving} gives the LDP for $R_\al^\e$ with speed $h^{-2}(\e)$ and rate $\CI_{C_{F,\rm{sd}}+C_{G,\rm{sd}}}$. Applying \lemref{lem:Tc} with $\nabla \bar F_{\rm sd}(t,\bar X(t))$ in place of $\overline{\nabla F}(\bar X(t))$ and the contraction principle to \eqref{wdh-eta-al}, we obtain the LDP of $\wdh \eta_\al^\e$. The corresponding Volterra map preserves $\SD_A$. Indeed, for $x\in \SD_A$, define $z(t):=x(t)-\int_0^t e^{(t-s)A} \nabla \bar F_{\rm sd}(s,\bar X(s)) x(s)ds$. Then $z\in \SD_A$ and $\dot{z}(t)-A z(t)=\dot{x}(t)-A x(t)-\nabla \bar F_{\rm sd}(t,\bar X(t)) x(t)$. By applying the same argument to the inverse Volterra map, we obtain the converse. Hence, the rate in \eqref{sd-main-rate} follows. By \propref{prop:equi-eta-al} and \lemref{lem:MP-energy}, the proof is complete.
\end{proof}

\section{Stochastic reaction--diffusion realizations}\label{sec:rd-applications}
The following conditions verify the abstract assumptions for the stochastic reaction--diffusion equation.
\begin{prop}\label{prop:rd-realization}
Let $\SO\subset\rr^d$ be bounded and regular,
$H=L^2(\SO)$, and $A$ be a negative self-adjoint uniformly
elliptic realization with Dirichlet or Neumann
boundary conditions and compact resolvent.  Suppose its eigenfunctions
satisfy $e_k\in L^\infty(\SO)$ and that $\mathsf{Q} e_k=\lambda_k^{\mathsf Q}e_k$.
If $d=1$, assume for some $\beta\in(0,1)$
\begin{equation}\label{kap-zeta}
 \sup_k|\lambda_k^{\mathsf Q}|<\infty,\qquad
 \sum_{k\geq1}\alpha_k^{-\beta}\|e_k\|_\infty^2<\infty.
\end{equation}
If $d\geq2$, assume for some $p_0>2$ and $\beta>0$ with
$\beta(p_0-2)/p_0<1$ that
\begin{equation}\label{kap-zeta2}
 \sum_{k\geq1}(\lambda_k^{\mathsf Q})^{p_0}\, \|e_k\|_\infty^2<\infty,
 \qquad
 \sum_{k\geq1}\alpha_k^{-\beta}\,\|e_k\|_\infty^2<\infty.
\end{equation}
Moreover, we assume that uniformly in $(\xi,z)$, the following conditions hold:
\begin{align*}
 |f_r(\xi,r,z)|&\leq K_F',& |f_{rr}(\xi,r,z)|&\leq K_F'',&
 \sup_z|f(\cdot,0,z)|_{L^2}&<\infty,\\
 |g(\xi,r,z)-g(\xi,r',z)|&\leq K_g'|r-r'|,&
 |g(\xi,r,z)|&\leq K_g.
\end{align*}
Define $F(x,z)(\xi):=f(\xi,x(\xi),z)$ and
$[G(x,z)y](\xi):=g(\xi,x(\xi),z)y(\xi)$. Then, Assumptions~\ref{ass:A} and~\ref{ass:fg} hold with $ \vartheta=\beta$ if $d=1$ and $\vartheta=\beta(p_0-2)/p_0$ if $d\geq 2$.
In particular, space-time white noise $\mathsf{Q}=I$ is included whenever
$d=1$ and \eqref{kap-zeta} holds.
\end{prop}
The proof is given in \appref{app:rd-verification}.
Proposition~\ref{prop:rd-realization} and Theorem~\ref{thm:main} give
the following reaction--diffusion examples.  They cover a bounded
additive environment and a nonuniformly ergodic environment.

\begin{exm}\label{ex:additive-environment}
Let $f(\xi,r,z)=f_0(r)+a(z)\psi(\xi)$, where
$f_0\in C^2(\rr)$ has bounded first and second derivatives, $a$ is
bounded, and $\psi\in H$. Since $f_0$ cancels in
\eqref{F-tz}, we have
$F_c(t,z)=(a(z)-\pi(a))\psi$ and $U(t,z)=u(z)\psi$, where
$u(z):=\int_0^\infty P_r(a-\pi(a))(z)\,dr$. Moreover,
\eqref{U-P-bound} gives
$|u(z)|\leq4\|a\|_\infty\mixer(z)/(-\log\mixrate)$, and
$C_F=2\pi\big((a-\pi(a))u\big)\,\psi\otimes\psi$.
Thus $u$ is bounded under uniform ergodicity, although its boundedness
is not required. Since $F_c$ is independent of time, the temporal regularity in \eqref{Fc-H-bounds} is automatic. The assumption $x_0\in\CD((-A)^\theta)$ is used only for \eqref{Fc-H-bounds}, so the conclusion of \thmref{thm:main} holds for every $x_0\in H$. The same coefficient is admissible in the switching model of \secref{sec:feedback}.
\end{exm}

\begin{exm}\label{ex:OU-environment}
A bounded solution of the Poisson equation does not require uniform ergodicity. Let
$d\zeta_t=-\zeta_tdt+dB_t$, where $B_t$ is independent
of $W_0$, with
$\pi=N(0,1/2)$, and set $u(z)=\arctan z$ and
$a(z)=z/(1+z^2)+z/(1+z^2)^2$. For
$\mathscr L=\frac12\partial_{zz}-z\partial_z$, we have
$-\mathscr Lu=a$, and both $a$ and $u$ are centered under $\pi$.
Dynkin's formula gives
$\int_0^tP_ra(z)\,dr=u(z)-P_tu(z)$. Since
$P_tu(z)\to\pi(u)=0$, $u$ is the bounded Poisson potential.
For $f(\xi,r,z)=f_0(r)+a(z)\psi(\xi)$,
$F_c(t,z)=a(z)\psi$ and $U(t,z)=u(z)\psi$, so the bounds \eqref{U-H-bounds}
hold for every $x_0\in H$. Gaussian integration by parts gives
$C_F=\sigma^2(a)\psi\otimes\psi$, where
$\sigma^2(a)=2\pi(au)=\pi\big((u')^2\big)
=\int_\rr(1+z^2)^{-2}\pi(dz)>0$.
Thus $C_F$ is nonzero and rank one when $\psi\neq0$.
\end{exm}

\begin{exm}\label{ex:feedback-observation}
Let $\upsilon_1,\ldots,\upsilon_{d_\Upsilon}\in H$ and define
$\Upsilon x=(\langle x,\upsilon_1\rangle_H,\ldots,
\langle x,\upsilon_{d_\Upsilon}\rangle_H)$. If
$\widetilde q_{ij}$ satisfies \assmref{ass:feedback}, the rates
$\widetilde q_{ij}(t,\Upsilon x)$ have the form required in \secref{sec:feedback}.
For $H=L^2(\SO)$, smooth mollifiers or normalized indicators give
bounded local-average observations. Point evaluation
$x\mapsto x(\xi_k)$ instead requires a state space continuously
embedded in $C(\overline{\SO})$.
\end{exm}

\section*{List of notation}
\leavevmode\par\nobreak\vskip4pt
{\small\parindent0pt\parskip0pt\lineskip0pt
\hrule\nobreak\vskip1.5pt
\notrow{Symbol}{Meaning}{Reference}\nobreak\vskip1.5pt\hrule\nobreak\vskip1.5pt
\notrow{$H$, $\CL(H)$, $\CL_2(H)$}{state space, bounded linear and Hilbert--Schmidt operators}{\secref{sec:assump}}
\notrow{$\mathsf{Q}$, $Q(t,y)$}{Wiener covariance and regime switching generator}{\secref{sec:assump}, \eqref{al-tran}}
\notrow{$W_0$, $W^{\mathsf{Q}}$}{cylindrical noises, $W^{\mathsf{Q}}=\mathsf{Q}W_0$}{\secref{sec:assump}}
\notrow{$C^0([0,1];H)$}{continuous paths starting at $0$, sup-norm $\|\cdot\|_\infty$}{\secref{sec:assump}}
\notrow{$D^0([0,1];H)$}{c\`adl\`ag paths from $0$, uniform unless specified}{\secref{sec:assump}}
\notrow{$A$, $e_k$, $\alpha_k$}{generator, eigenvectors, eigenvalues}{\assmref{ass:A}}
\notrow{$\mathcal H_\delta$, $\Pi_n$, $H_n$, $A_n$}{fractional Sobolev spaces, Galerkin projection, subspace and projected generator}{\secref{sec:assump}}
\notrow{$F$, $G$, $\vartheta$}{drift, noise coefficient, smoothing exponent}{\assmref{ass:fg}}
\notrow{$h(\e)$, $J_N^{\rm hom}$}{deviation scale and rate function of $\Ga_\e^{\rm hom}$}{\eqref{dev-s}, \eqref{J-N}}
\notrow{$L_C$, $\CI_C$, $\SD_A$}{local rate function, rate function and path space}{\eqref{LC}, \eqref{quad-L}, \secref{sec:intro}}
\notrow{$[C^{1/2}]^\dagger$}{minimum-norm preimage on $\Ran(C^{1/2})$}{\eqref{LC-t}}
\notrow{$\Scv$}{convolution operator $v\mapsto\int_0^te^{(t-s)A}v(s)\,ds$}{\appref{app:quad-lift}}
\notrow{$\Phi_n$, $\Psi^N$}{spectral Volterra and time-grid contraction maps}{\eqref{def-Phi-n}, \eqref{Phi-map}}
\notrow{$\CI_n$, $\CI_{G,n}$}{quadratic rates for spectral projections}{\lemref{lem:quad-lift}, \eqref{new:IN}}
\notrow{$\zeta$, $P_t$, $\pi$, $\mixrate$, $\mixer$}{environment, semigroup, invariant law, mixing rate and weight}{\eqref{exp-erg}}
\notrow{$C_\ST$}{Green--Kubo covariance of a centered map $\ST$}{\eqref{CF-GK}}
\notrow{$F_c$, $F_k$, $F^N$}{centered drift, its values at $t_k$ and the vector $F^N$}{\eqref{F-tz}, \lemref{lem:LDP-L}}
\notrow{$U(t,z)$, $U(t,i)$, $\fU$}{solutions of the Poisson equations and unitary dilation group}{\eqref{def-U}, \eqref{Poi-al}, \eqref{dila}}
\notrow{$C_F$, $C_G$, $C$}{environment, Wiener and total covariances}{\eqref{CF}, \eqref{CG}}
\notrow{$\la^\e$, $M^\e$, $R^\e$}{environmental and Wiener terms and their sum}{\eqref{la}, \eqref{M}, \eqref{def-R}}
\notrow{$\overline{\nabla F}$, $\CT$}{averaged derivative and Volterra map}{\secref{sec:intro}, \eqref{Tu}}
\notrow{$\wdh\eta^\e$, $\wdh\eta_\al^\e$}{auxiliary linearized processes}{\eqref{hat-eta}, \eqref{wdh-eta-al}}
\notrow{$D_\e^N$, $L_\e^N$, $A_\e^N$}{time-discretization of $\la^\e$}{\eqref{decomp-la}}
\notrow{$\mathscr S_N$, $\CI_N^{\mathrm{hom}}$}{piecewise semigroup paths and time-grid rate}{\lemref{lem:LDP-L}, \eqref{IN}}
\notrow{$\zeta^{(\ell)}$, $\wdh P$}{unit-time path chain and its kernel}{\corref{cor:zeta-path}}
\notrow{$\varrho(i)$, $N_{\mathfrak a}(n)$}{regeneration times and number of visits to the atom}{\appref{app:additive-holder}}
\notrow{$\CI_\zeta$, $\CI_w$, $\CI$}{rates of $\la^\e$, $M^\e$ and $R^\e$}{\eqref{I}, \eqref{new:IW}, \eqref{new:sumrate}}
\notrow{$\al^\e$, $m_0$, $\Upsilon$}{switching process, number of states and observation}{\eqref{al-tran}}
\notrow{$\bar y_t$, $\nu(t)$}{observed limit $\Upsilon\bar X(t)$ and invariant vector $\pi(t,\bar y_t)$}{\lemref{lem:avg-al}}
\notrow{$\Bpi$, $\nabla\bar F_{\rm sd}$}{extra drift operator and drift derivative}{\eqref{def-Bpi}, \eqref{nab-f}}
\notrow{$C_{F,{\rm sd}}$, $C_{G,{\rm sd}}$}{covariances in the switching case}{\eqref{sd-CF}, \eqref{sd-CG}}
\vskip1.5pt\hrule\par}

\appendix
\section{Common analytic and Galerkin estimates}\label{app:common-tools}

\subsection{Proof of \lemref{lem:MP-energy}}\label{app:lem-MP-energy}
Since $\qv{Cy,y}=\|T^* y\|^2$, we have $L_C(v)=\sup_{y}\{\qv{v,y}-\half \|T^* y\|^2\}$. If $T u=v$, then
$
\qv{v,y}-\half \|T^* y\|^2 = \qv{u,T^* y}-\half \|T^* y\|^2 \leq \half \|u\|^2.
$
Hence $L_C(v)\leq \half \|u\|^2$ for every preimage $u$ of $v$. In particular, $L_C(v)<\infty$ if $v\in\Ran T$.

Conversely, suppose $J:=L_C(v)<\infty$. Replacing $y$ by $\tau y$ with $\tau \in \rr$ gives $\tau \qv{v,y}-\half \tau^2 \|T^* y\|^2 \leq J$ for all $\tau\in \rr$. If $\|T^* y\|^2 >0$, maximizing over $\tau\in \rr$ gives $\qv{v,y}^2 \leq 2J \|T^* y\|^2$. If $T^* y=0$, then $\tau \qv{v,y}\leq J, \forall\,\tau\in \rr$, hence $\qv{v,y}=0$. In both cases
\beq{vy-leq}
|\qv{v,y}|\leq \sqrt{2J}\|T^* y\|, \quad y\in H.
\eeq
Thus, $T^* y \mapsto \qv{v,y}$ is a well-defined linear functional on $\Ran T^*$ of norm at most $\sqrt{2J}$. It extends to the closure $\lbar{\Ran T^*}=(\ker T)^\perp$, and the Riesz representation theorem gives a $u\in (\ker T)^\perp$ such that $\qv{u,T^* y}=\qv{v,y}$ for all $y\in H$. Therefore, $\qv{T u-v, y}=0$ for all $y$, i.e., $T u = v$ and \eqref{vy-leq} implies $\|u\|^2\leq 2J$. Every other solution of $T u'=v$ has the form $u'=u+u^o$ with $u^o\in \ker T$. Thus, $\|u'\|^2 = \|u\|^2 + \|u^o\|^2 \geq \|u\|^2$ since $u\perp \ker T$. Hence $u$ is the unique solution of minimum norm and then $J=\half \|u\|^2$. This proves \eqref{LC}, including the value $+\infty$ for $v\notin \Ran T$. For each fixed $t\in [0,1]$, let $T=C(t)^{1/2}$, the minimum-norm solution of $C(t)^{1/2}u=v$ is by definition $[C(t)^{1/2}]^\dagger v$. It then gives \eqref{LC-t}.

For measurability, let $\dl\downarrow 0$ and define $
u_\dl(t):=C(t)^{1/2}[C(t)+\dl I]^{-1} v(t)$.
Strong measurability follows from the continuous functional calculus and the uniform
boundedness of $C$. Indeed, the continuous scalar function $\la\mapsto \la^{1/2} (\la+\dl)^{-1}$ and $\la \mapsto \la^{1/2}$ can be uniformly approximated by polynomials on the common spectral interval $[0,\sup_{t}\|C(t)\|]$, and polynomials in a strongly measurable bounded operator field are strongly measurable. At a point where $v(t)=C(t)^{1/2} u(t)$ with $u(t)\perp \ker[C(t)^{1/2}]=\ker[C(t)]$, the spectral theorem gives
\bea
u_\dl(t)=C(t)[C(t)+\dl I]^{-1} u(t) \to u(t) \text{ and } |u_\dl(t)|\leq |u(t)|,
\eea
because $\la/(\la+\dl) \uparrow 1$ for $\la>0$ and the spectral measure of $u(t)$ does not charge $\{0\}$. The finite energy set $\{t:L_C(t,v(t)) <\infty\}$ is measurable, since \eqref{LC} is a supremum over a countable dense subset of $H$, and $t\mapsto L_C(t,v(t))$ is measurable for the same reason. Hence $u=\lim_{\dl\to0} u_{\dl} \indi_{\{L_C(\cdot,v)<\infty\}}$ is strongly measurable, where we extend $u$ by zero on the complement of the finite energy set. Eq. \eqref{LC-t} gives $\|u\|_{L^2}^2 =2\int_0^1 L_C(t,v(t))dt<\infty$. Since $|u_\dl(t)|\leq |u(t)|$ and $u\in L^2([0,1];H)$, dominated convergence implies that $u_\dl \to u$ in $L^2([0,1];H)$. The proof is complete. \qed

\subsection{Proof of Corollary~\ref{cor:add}}\label{app:cor-add}
By \lemref{lem:MP-energy} with $T=C_1^{1/2}$, every $v_1$ of finite energy is of the form $v_1=C_1^{1/2}u$ with $L_{C_1}(v_1)=\half |u|^2$ for the minimal $u$, and every $u\in H$ gives $L_{C_1}(C_1^{1/2}u)\leq \half |u|^2$. Therefore, $
\inf\{L_{C_1}(v_1): \CK v_1 =w\}=\half \inf\{|u|^2: \CK C_1^{1/2} u=w\}$.
Now apply \eqref{LC} to $T=\CK C_1^{1/2}\in \CL(H,H')$, for which $T T^*=\CK C_1 \CK^*$. Its right-hand side equals $L_{\CK C_1 \CK^*}(w)$ and the infimum over $u$ is a minimum. This proves \eqref{LC1}. To prove \eqref{LC1+C2}, take $\CU=H\oplus H$ and $T(u_1,u_2)=C_1^{1/2}u_1 + C_2^{1/2} u_2$, so that $T T^* = C_1 + C_2$ and $\|(u_1,u_2)\|^2 =|u_1|^2 +|u_2|^2$. Then, similarly to the argument above,
\bea\ad
\inf_{v_1+v_2=v}\{L_{C_1}(v_1)+L_{C_2}(v_2)\}\!=\!\half \inf\{|u_1|^2+|u_2|^2: C_1^{1/2}u_1 + C_2^{1/2} u_2 = v\} \!=\! L_{C_1+C_2}(v).
\eea
By \eqref{LC}, the minimum is attained at the minimal-norm solution. \qed

\subsection{Proof of Lemma~\ref{lem:quad-lift}}\label{app:quad-lift}

Choose $K\geq1$ such that
$\sup_{t\in[0,1]}\|C(t)\|_{\CL(H)}\leq K$. Then,
\beq{low-LC}
L_C(t,v)
=
\sup_{u\in H}
\big\{
\qv{v,u}_H-\frac12\qv{C(t)u,u}_H
\big\}
\geq |v|_H^2/2K,\quad \forall\,v\in H.
\eeq
The same holds for $L_{C_n}$ on $H_n$. Let
$
J(x):=\sup_{n\geq1}\CI_n(\Pi_nx),
x\in C^0([0,1];H).
$

We first prove that $J=\CI_C$. Let $m\leq n$. If $\CI_n(\Pi_nx)=+\infty$, then
$\CI_m(\Pi_mx)\leq \CI_n(\Pi_nx)$ is immediate. Suppose that
$\CI_n(\Pi_nx)<\infty$. Set
$
x_n:=\Pi_nx$, $
v_n:=\dot x_n-A_nx_n$.
Then $\Pi_mx=\Pi_mx_n$ and, since $\Pi_m$ commutes with $A_n$ on
$H_n$,
$
v_m:=(\Pi_mx)^{\cdot}-A_m\Pi_mx=\Pi_mv_n.
$
For every $u\in H_m$, viewed also as an element of $H_n$,
$
\qv{v_m(t),u}_H=\qv{v_n(t),u}_H$ and $
\qv{C_m(t)u,u}_H
=
\qv{C_n(t)u,u}_H
=
\qv{C(t)u,u}_H.
$
Hence, for almost every $t$,
\bea
L_{C_m}(t,v_m(t))
\ad=
\sup_{u\in H_m}
\big\{
\qv{v_n(t),u}_H
-\frac12\qv{C(t)u,u}_H
\big\} \\
\ad \leq
\sup_{u\in H_n}
\big\{
\qv{v_n(t),u}_H
-\frac12\qv{C(t)u,u}_H
\big\}
=
L_{C_n}(t,v_n(t)).
\eea
Integrating over $[0,1]$ gives
$
\CI_m(\Pi_mx)\leq \CI_n(\Pi_nx).
$
Thus $\CI_n(\Pi_nx)\uparrow J(x)$ as $n\to\infty$.

We next show that $J(x)<\infty$ implies $x\in\SD_A$. Suppose that
$J(x)\leq R$. Then $\CI_n(\Pi_nx)\leq R$ for every $n$. With
$x_n=\Pi_nx$ and $v_n=\dot x_n-A_nx_n$, \eqref{low-LC} gives
$
\int_0^1|v_n(t)|_H^2\,dt
\leq
2K  \CI_n(\Pi_nx)
\leq
2KR.
$
Since $A_n$ is negative and self-adjoint and $x_n(0)=0$, we have
\bea
\int_0^1|v_n(t)|_H^2\,dt
&=
\int_0^1
\bigl(
|\dot x_n(t)|_H^2+|A_nx_n(t)|_H^2
\bigr)\,dt
+|(-A_n)^{1/2}x_n(1)|_H^2.
\eea
Hence both $\{\dot x_n\}$ and $\{A_nx_n\}$ are bounded in
$L^2([0,1];H)$. Since $x\in C^0([0,1];H)$, the set $x([0,1])$ is compact in $H$. The strong
convergence $\Pi_n\to I$ is thus uniform on $x([0,1])$, and
$
\|\Pi_nx-x\|_\infty\to0.$
After passing to a subsequence,
$
\dot x_n\rightharpoonup y,$ and $
A_nx_n\rightharpoonup z
$
weakly in $L^2([0,1];H)$. The convergence
$x_n\to x$ in $L^2([0,1];H)$ and the weak convergence of
$\dot x_n$ imply that $x\in W^{1,2}([0,1];H)$ and $\dot x=y$.
Moreover, $A_nx_n=A\Pi_nx$. The operator induced by $A$ on
$L^2([0,1];H)$ has a closed linear graph and hence a weakly
closed graph. It yields
$
x\in L^2([0,1];\CD(A))$, $
Ax=z.$ As $x_n(0)=0$ and $x_n\to x$ uniformly, $x(0)=0$. Thus
$x\in\SD_A$.

It remains to identify $J(x)$ for $x\in\SD_A$. Set
$
v:=\dot x-Ax\in L^2([0,1];H).
$
Then
$
(\Pi_nx)'-A_n\Pi_nx=\Pi_nv.
$
Therefore,
$
 \CI_n(\Pi_nx)
=
\int_0^1
L_{C_n}(t,\Pi_nv(t))\,dt.
$
For every $t\in[0,1]$,
\bea
L_{C_n}(t,\Pi_nv(t))
\ad=
\sup_{u\in H_n}
\big\{
\qv{\Pi_nv(t),u}_H
-\frac12\qv{C_n(t)u,u}_H
\big\} \\
\ad=
\sup_{u\in H_n}
\big\{
\qv{v(t),u}_H
-\frac12\qv{C(t)u,u}_H
\big\}.
\eea
The spaces $H_n$ are increasing and
$\cup_{n\geq1}H_n$ is dense in $H$. For fixed $t$, the map
$
u\mapsto
\qv{v(t),u}_H-\frac12\qv{C(t)u,u}_H
$
is continuous on $H$. Hence
\bea\ad
L_{C_n}(t,\Pi_nv(t))
\uparrow
\sup_{u\in H}
\big\{
\qv{v(t),u}_H
-\frac12\qv{C(t)u,u}_H
\big\}
=
L_C(t,v(t)).
\eea
Since these functions are nonnegative,
the monotone convergence theorem yields
\bea
J(x)\ad=
\lim_{n\to\infty} \CI_n(\Pi_nx) =
\int_0^1
\lim_{n\to\infty}
L_{C_n}(t,\Pi_nv(t))\,dt =
\int_0^1 L_C(t,v(t))\,dt
=
 \CI_C(x).
\eea
Together with the implication $J(x)<\infty\Rightarrow x\in\SD_A$,
this proves $J=\CI_C$ on $C^0([0,1];H)$.

We next prove that $ \CI_C$ is a good rate function. Let $(\Scv v)(t):=\int_0^te^{(t-s)A}v(s)\,ds$ for $v\in L^2([0,1];H)$ and, for $R<\infty$, set
\[
\SK_R:=\big\{\Scv v:\  \|v\|_{L^2([0,1];H)}^2\leq2KR\big\}.
\]%
Let
$ \CI_C(x)\leq R$ and $v:=\dot x-Ax$. By \eqref{low-LC},
$
\|v\|_{L^2([0,1];H)}^2\leq2KR.
$
Since $x(0)=0$,
$
x(t)=\int_0^t e^{(t-s)A}v(s)\,ds.
$
Thus $\{x\in C^0([0,1];H):\CI_C(x)\leq R\}\subset\SK_R$. We show that $\SK_R$ is relatively compact in $C^0([0,1];H)$. Let $x\in\SK_R$ with the corresponding $v$.
For every $n$,
\bea
|(I-\Pi_n)x(t)|_H^2
\ad=
\sum_{k\geq n+1}
\Big|
\int_0^t e^{-\alpha_k(t-s)}
\qv{v(s),e_k}_H\,ds
\Big|^2 \leq
\frac{1}{2\alpha_{n+1}}
\int_0^1|v(s)|_H^2\,ds.
\eea
Consequently,
\beq{IC-tail}
\|(I-\Pi_n)x\|_\infty
\leq
\left(\frac{KR}{\alpha_{n+1}}\right)^{1/2}.
\eeq
The right-hand side tends to zero uniformly over
$\SK_R$. For fixed $m$, we have $ \Pi_mx(t) = \int_0^t e^{(t-s)A_m}\Pi_mv(s)\,ds. $ The bound on $v$ gives
$
\|\Pi_mx\|_{W^{1,2}([0,1];H_m)}
\leq
(2+\alpha_m)\|v\|_{L^2([0,1];H)}
\leq
(2+\alpha_m)(2KR)^{1/2}.
$
Since $H_m$ is finite dimensional,
$\{\Pi_mx: x\in\SK_R\}$ is relatively compact in
$C^0([0,1];H_m)$. Combining this finite-dimensional compactness
with \eqref{IC-tail} shows that
$\SK_R$, and hence $\{x\in C^0([0,1];H): \CI_C(x)\leq R\}$, is relatively compact in $C^0([0,1];H)$. Since $R$ is arbitrary, $\Scv$ is compact from $L^2([0,1];H)$ to $C^0([0,1];H)$.

For every $n$, $\CI_n$ is lower semicontinuous on
$C^0([0,1];H_n)$, and $\Pi_n:C^0([0,1];H)\to C^0([0,1];H_n)$ is continuous.
Since
$
 \CI_C(x)=\sup_{n\geq1} \CI_n(\Pi_nx),
$
$ \CI_C$ is lower semicontinuous on $C^0([0,1];H)$. Its level sets are
closed and relatively compact, so
$ \CI_C$ is good.

Finally, let
$
z\in C^0([0,1];H_n),\,
 \CI_n(z)\leq R,
$
and set $v=\dot z-A_nz$. By \eqref{low-LC},
$
\|v\|_{L^2([0,1];H)}^2\leq2KR,
$
where $v$ is viewed as an $H$-valued function. Moreover,
$
z(t)=\int_0^t e^{(t-s)A}v(s)\,ds.
$
Hence $z\in\SK_R$ for every $n$, so
$
\cup_{n\geq1}
\{
z\in C^0([0,1];H_n):\CI_n(z)\leq R
\}\subset\SK_R$
is relatively compact in $C^0([0,1];H)$. This proves the last assertion. \qed

\section{Exponentially ergodic Markov environment}\label{app:independent-environment}

\subsection{Proof of Lemma~\ref{lem:A1}}\label{app:A1}
We first recall the result of \cite[Theorem~1.1]{Che97}.
\begin{thm}[Chen \cite{Che97}]\label{thm:che97}
Let $\{a_k\}_{k\geq 1}$ satisfy the assumptions preceding Lemma~\ref{lem:A1}. Define $S_n=\sum_{k=1}^n a_k$. Suppose that
$S_n/b_n \to 0$ in probability. Then $\{S_n/b_n\}_{n\geq 1}$
satisfies an MDP with speed
$n/b_n^2$ and rate function $
 \widetilde\CI(z)=\sup_{V\in E^*}
 \{V(z)-\tfrac12\sigma_V^2\}, z\in E$, where $
 \sigma_V^2=\EE V(a_1)^2+
 2\sum_{j=1}^m\EE[V(a_1)V(a_{1+j})],\, V\in E^*$.
\end{thm}
We need the following lemma to prove our result.
\begin{lem}\label{lem:st}
Let $0\leq s<t\leq 1$ be fixed. If the conditions in \thmref{thm:che97} hold, then $\{\UU_n(t)-\UU_n(s)\}_{n\geq 1}$ satisfies an MDP in $E$ with speed $n/b_n^2$ and rate function
\begin{equation}\label{Ist}
 \widetilde\CI_{s,t}(z)=(t-s)\sup_{V\in E^*}
 \Big\{V\Big(\frac z{t-s}\Big)-\frac12\sigma_V^2\Big\},\qquad z\in E.
\end{equation}
\end{lem}

\begin{proof}
The function $\wdt\CI_{s,t}$ defined in \eqref{Ist} has compact level sets due to \thmref{thm:che97}. Given $0\leq s<t \leq 1$, let $n_k:=\lf kt \rf- \lf ks \rf
$. Define sequences $\wdt a_1=a_{\lf k s\rf
+1}, \wdt a_2=a_{\lf ks \rf
+2},\dots, \wdt a_{n_k}=a_{\lf kt \rf
}$ and $\wdt S_{n_k}=\sum_{i=1}^{n_k}\wdt a_i$. Hence, since $b_{n_k}\leq b_k$, $\wdt S_{n_k}/ b_k \to 0$ in probability. Note that $\frac{k}{b_k^2}\cdot \frac{b_k^2}{n_k}\to \frac{1}{t-s}$ as $k\to \infty$. Thus, \thmref{thm:che97} implies our result.
\end{proof}

\begin{proof}[Proof of Lemma~\ref{lem:A1}]
Given a partition $\CP: 0\leq t_0< t_1<\dots < t_N =1$, denote $\varpi_n^\CP=(\UU_n(t_1),\UU_n(t_2)-\UU_n(t_1),\dots, \UU_n(t_N)-\UU_n(t_{N-1}))$. Since $\{a_k\}$ is $m$-dependent, removing the first $m$ terms of each block makes the $N$ block sums independent, and the removed terms are negligible at speed $n/b_n^2$ by the exponential integrability of $\|a_1\|_E$. By \lemref{lem:st} and \cite{DZ09}, the sequence $\{\varpi_n^\CP\}_{n\geq 1}$ satisfies an MDP in $E^N$ with speed $n/b_n^2$ and rate function $\sum_{i=1}^N\wdt\CI_{t_{i-1},t_i}(z_i)$, $z=(z_1,\dots,z_N)$, with $\wdt\CI_{s,t}$ in \eqref{Ist}.
Then by the contraction principle for the map $(z_1,\dots,z_N)\mapsto (z_1,z_1+z_2, \dots, \sum_{i=1}^N z_i)$ on $E^N$, the proof is completed.
\end{proof}

\subsection{Proof of \lemref{lem:B-mdp}}\label{app:Bmdp}
Let $\{X_k\}_{k\geq0}$ be an irreducible Markov chain with
transition kernel $P$. It is well-known that it has some small set, i.e., $\exists\; 0<{}\iota <1, \exists\; \CS \in \SB(\SE)$ and $\exists\; \nu$ a probability measure on $\SE$ such that $\nu(\CS)>0$ and  $
P^m \geq \iota \indi_{\CS}\otimes \nu$.

The split chain method is an embedding of the Markov chain $\{X_k\}_{k\geq 0}$ into a larger probability space on which a $\{0,1\}$-valued random sequence $\{Z_k\}_{k\geq 0}$ is also defined. This embedding is done in such a way that the sequences
\beq{spl1}
\{\wdt \Phi_n\}_{n\geq 0} =\{(X_{nm}, Z_n)\}_{n\geq 0}, \quad \{\wdt \Psi_n\}_{n\geq 1}=\{(X_{(n-1)m+1}, \dots, X_{nm},Z_n)\}_{n\geq 1}
\eeq
\beq{spl2}
\{(X_0,Z_0), (X_{(n-1)m+1}, \dots, X_{nm}, Z_n)\}_{n\geq 1}
\eeq
are Markov chains inheriting properties of the initial chain $\{X_k\}_{k\geq 0}$.  The chain
$(\widetilde\Phi_n)$ has an atom $\mathfrak a$, that is,
$\widetilde\pi(\mathfrak a)>0$ and there exists a probability measure
$\widetilde\nu$ such that
\bea
\bar P(x,\cdot)=\wdt\nu\cd,\quad \forall\, x\in \mathfrak a.
\eea
This feature enables us to decompose \eqref{spl2} into blocks of return to this atom $\mathfrak a$. These blocks are independent or $1$-dependent random sequences. Let $\II=\{0,1\}$ and denote by $2^{\II}$ the power set of $\II$. For any given probability measure $\mu$ on $(\SE,\SB(\SE))$, we could define the measure $\wdt \mu$ on $(\SE\times \II, \SB(\SE)\times 2^{\II}) $ by
\bea
\wdt \mu = \mu\big((1-\iota\indi_{\CS})\indi_{(\cdot)}\big)\otimes \dl_0 +\mu\big(\iota \indi_{\CS}\indi_{(\cdot)}\big)\otimes \dl_1,
\eea
where $\dl_i\cd$ are Dirac measures for $i=0,1$. Moreover, we could also define the Markov kernels $\wdt P$ on $(\SE\times \II) \times (\SB^m (\SE) \otimes 2^{\II})$ and $\bar P$ on $(\SE \times \II)\times (\SB(\SE) \times 2^{\II})$. Denote by $\wdt{\PP}_{\mu^*}$ the law of $\{\wdt{\Psi}_n\}$ with transition $\wdt P$ and initial measure $\mu^*=\wdt{\PP}_{\mu^*}((X_{-m+1},\dots,X_0, Z_0)\in \cdot)$ satisfying $
\wdt{\PP}_{\mu^*}((X_0,Z_0)\in \cdot)=\wdt\mu\cd.$
Denote
\bea
\wdt{\PP}_{\mu^*}\Big|_{(\SB(\SE)\otimes2^{\II})\otimes_{n=1}^\infty(\SB^m(\SE)\otimes 2^{\II})}=\wdt{\PP}_{\wdt \mu}.
\eea
Then we have the following results, which can be found in Guillin \cite[Lemma 1]{Gui00}.
\begin{lem}\label{lem:reg} We have
\begin{itemize}
\item[1.] Under the law $\wdt{\PP}_{\wdt \mu}$, the sequence $\{\wdt\Phi_n\}_{n\geq 0}$ with general state space $(\SE\times \II, \SB(\SE)\otimes 2^{\II})$, transition probability $\bar P$ and initial distribution $\wdt\mu$, is an aperiodic Harris recurrent Markov chain.
\item[2.] If $\pi$ is invariant for $\{X_n\}$, then $\wdt \pi$ is an invariant measure for the chain $\{\wdt\Phi_n\}$.
\item[3.] $\mathfrak a=\CS\times \{1\}$ is an atom for the chain $\{\wdt\Phi_n\}$.
\item[4.] The following equality in law holds: $
\operatorname{Law}_{\wdt{\PP}_{\wdt\mu}}(\{X_n\})=\operatorname{Law}_{\PP_{\mu}}(\{X_n\})$.
\end{itemize}
\end{lem}
Define the successive return times to the atom $\mathfrak a$ by
\beq{return}\left\{\barray
\tau(0)=\inf\{n\geq 0:  \wdt \Phi_n =(X_{nm}, Z_n)\in \mathfrak a\}\\
\tau(k+1) =\inf\{n> \tau(k): \wdt \Phi_n =(X_{nm}, Z_n)\in \mathfrak a\}
\earray\right.\eeq
If $\widetilde{\mathcal F}_n=\sigma(\widetilde\Phi_0,\ldots,
\widetilde\Phi_n)$, then
$
 \{\tau(0)\leq n\}
 =\bigcup_{j=0}^n\{\widetilde\Phi_j\in\mathfrak a\}
 \in\widetilde{\mathcal F}_n$,
and, inductively,
\bea
 \{\tau(k+1)\leq n\}
 =\bigcup_{j=1}^n
 \{\tau(k)<j,\ \widetilde\Phi_j\in\mathfrak a\}
 \in\widetilde{\mathcal F}_n.
\eea
Thus all $\tau(k)$ are stopping times.  Harris recurrence and
$\widetilde\pi(\mathfrak a)>0$ imply
$\widetilde\PP_{\widetilde\mu}(\tau(k)<\infty)=1$.  Kac's return-time
formula at the atom gives
$\EE_{\widetilde\nu}[\tau(0)+1] =\widetilde\pi(\mathfrak a)^{-1}$.
The following lemma can be found in Guillin \cite[Lemma~2]{Gui00}.
\begin{lem}\label{lem:reg-geo} We have
\begin{itemize}
\item[1.] Under the law $\wdt\PP_{\wdt \mu}$, the blocks $\{X_{m\tau(k-1)+m},\dots, X_{m\tau(k)+m-1}\}$ are $1$-dependent random variables with common law
$
\operatorname{Law}_{\wdt\PP_{\wdt \nu}}\big(\{X_0, \dots, X_{m\tau(0)+m-1}\}\big)$.
Moreover, $\{\tau(k)-\tau(k-1)\}_{k\geq1}$ are i.i.d.\ random variables with common law $\operatorname{Law}_{\wdt \PP_{\wdt \nu}}(\tau(0)+1)$.
\item[2.] If $\{X_n\}$ is geometrically ergodic, then so is $\{\wdt\Phi_n\}$, and there is $\kappa>0$ such that $\EE_{\wdt\nu}e^{\kappa\tau(0)}<\infty$ and $\EE_{\wdt\mu}e^{\kappa\tau(0)}<\infty$ for every $\mu$ satisfying \eqref{chain-cond-mu}, see \cite[Chapter~IV, Lemma~2.2 and (2.10)--(2.12)]{Che99}.
\end{itemize}
\end{lem}

\noindent{{\bf Step 1}}. We decompose the Markov chain $\{X_k\}$ into blocks of return to the atom $\mathfrak a$ as follows:
\beqa{decomp}\disp
\sum_{k=0}^M \fT(X_k)
\ad = \sum_{k=0}^{(m\tau(0)+m-1)\wedge M} \fT(X_k) +\sum_{k=1}^{\mathfrak v(M)-1} \xi_k(\fT) + \sum_{m(\ell(M)+1) \leq k \leq M} \fT(X_k) \\
\ad =: S_M^{\rm ini} + \sum_{k=1}^{\mathfrak v(M)-1}\xi_k(\fT) + S_M^{\rm ter},
\eeqa
where $
\xi_k(\fT)=\sum_{j=m(\tau(k-1)+1)}
^{m\tau(k)+m-1} \fT(X_j)$.  The $\mathfrak v(M)$ is the number of times of return to the atom $\mathfrak a$ before $M$: $
\mathfrak v(M) =\sum_{k=0}^{\lfloor M/m\rfloor-1} \indi_{\mathfrak a}(\wdt \Phi_k)$,
and $\ell(M)=\tau([\mathfrak v(M)-1]\vee 0)$ is the last moment that $\{\wdt \Phi_k\}$ visits $\mathfrak a$. By \lemref{lem:reg-geo}, the $\{\xi_k\}$ is a sequence of $1$-dependent random variables in the Banach space $E$, having common law $
\operatorname{Law}_{\wdt \PP_{\wdt \nu}}(\sum_{k=0}^{m\tau(0)+m-1} \fT(X_k))$.
Introduce $e(M)=\lfloor\wdt \pi(\mathfrak a)M/m\rfloor$, then the decomposition \eqref{decomp} can be rewritten as
\beq{decomp-2}\barray\disp
\sum_{k=0}^{M} \fT(X_k)
\ad = \sum_{k=1}^{e(M)} \xi_k + S_M^{\rm ini} + \Big(\sum_{k=1}^{\mathfrak v(M)-1} \xi_k(\fT) - \sum_{k=1}^{e(M)} \xi_k(\fT)\Big) + S_M^{\rm ter} \\
\ad =: \sum_{k=1}^{e(M)} \xi_k + S_M^{\rm ini} +S_M^{\rm rem} + S_M^{\rm ter}.
\earray\eeq
Let $M=\lf nt_i\rf,
i=1,\dots, N$. It can be shown that only the first term of \eqref{decomp-2} will contribute to the moderate deviations. The hypothesis $\UU_n(1)\to0$ of \lemref{lem:A1} holds for $\{\xi_k\}$ by the assumption on $\MM_n(1)(\fT)$; see \cite[Chapter~IV, (2.21)]{Che99}. Leveraging \lemref{lem:A1}, the exponential integrability of $\{\xi_k\}$ and the identification of the variance in \eqref{iden}, it follows that
\bea\ad
\frac{1}{b_n} \Big(\sum_{k=1}^{e(\lf nt_1\rf
)} \xi_k(\fT), \dots, \sum_{k=1}^{e(\lf nt_N\rf
)} \xi_k(\fT)\Big)
\eea
satisfies an MDP with speed $n/b_n^2$ having the good rate function
\beq{JTP}
J_\fT^\CP = \sum_{i=1}^N (t_i - t_{i-1}) J_\fT \Big(\frac{z_i - z_{i-1}}{ t_i - t_{i-1}} \Big),
\eeq
where $J_\fT$ is given by \eqref{JT}, $z=(z_1,\dots, z_N)$, and $z_0\mathrel{:=}0$.
To obtain the desired MDP,  it suffices to establish the following negligibility of terms $S_M^{\rm ini}$, $S_M^{\rm ter}$, and $S_M^{\rm rem}$. That is, for each $i \in \{1,\dots, N\}$ and every $r>0$,
\beq{neg-1}
\limsup_{n\to \infty} \frac{n}{b_n^2} \log\wdt \PP_{\wdt \mu}\Big(\Big\|\sum_{k=0}^{(m\tau(0)+m-1)\wedge \lf nt_i\rf
} \fT(X_k)\Big\|_E  > r b_n \Big) =-\infty.
\eeq
\beq{neg-2}
\limsup_{n\to \infty} \frac{n}{b_n^2} \log\wdt \PP_{\wdt \mu} \Big(\Big\|\sum_{m(\ell(\lf nt_i\rf
)+1)\leq k\leq \lf nt_i\rf
} \fT(X_k)\Big\|_E  > r b_n \Big) =-\infty.
\eeq
\beq{neg-3}
\limsup_{n\to \infty} \frac{n}{b_n^2} \log \wdt \PP_{\wdt \mu} \Big(\Big\| \sum_{k=1}^{\mathfrak v(\lf nt_i\rf
)-1} \xi_k(\fT) -\sum_{k=1}^{e(\lf nt_i\rf
)} \xi_k(\fT)\Big\|_E  > r b_n \Big) =-\infty.
\eeq
These are (2.14)--(2.16) in the proof of Chen \cite[Chapter~IV, Theorem~2.1]{Che99}, with $n$ replaced by $\lfloor nt_i\rfloor$. That proof uses the boundedness of $\fT$, the decomposition \eqref{decomp}, the exponential moments in \lemref{lem:reg-geo}, and (2.21) there for \eqref{neg-3}. Thus, for every finite partition $\CP$ of $[0,1]$, $\PP_\mu(\MM_n^\CP\in\cdot)$ satisfies the MDP with speed $n/b_n^2$ and rate function $J_\fT^\CP$ in \eqref{JTP}, by using $\PP_\mu (\MM_n^\CP \in\cdot) = \wdt \PP_{\wdt \mu}(\MM_n^\CP \in\cdot )$, where $\MM_n^\CP:=(\MM_n(t_1)(\fT),\dots,\MM_n(t_N)(\fT))$. By Djellout--Guillin--Wu \cite[Lemma A.1]{DGW99}, we also need \eqref{A2}, which is proved in the next step.
\smallskip

\def\nns{{\lf ns\rf}}
\def\nnt{{\lf nt \rf}}
\noindent{{\bf Step 2}}. We claim that for any $r>0$,
\beq{A2}
\lim_{\dl\to 0} \sup_{0\leq s\leq 1} \limsup_{n\to \infty} \frac{n}{b_n^2} \log \PP_\mu \Big( \sup_{s\leq t\leq (s+\dl)\wedge 1} \|\MM_n(t)(\fT)-\MM_n(s)(\fT)\|_E  >r \Big) = -\infty.
\eeq
In fact, the Markov property implies
\beq{M-diff}\barray
\ad\PP_\mu\Big(\sup_{s\leq t\leq(s+\dl)\wedge1}
 \|\MM_n(t)(\fT)-\MM_n(s)(\fT)\|_E>r\Big)\\
\ad=\PP_\mu\Big(\sup_{s\leq t\leq(s+\dl)\wedge1}
 \Big\|\sum_{k=\nns}^{\nnt-1}\fT(X_k)\Big\|_E>rb_n\Big)\\
\ad=\PP_{\mu P^{\nns}}\Big(
 \max_{0\leq j\leq \lf n((s+\dl)\wedge1)\rf-\nns-1}
 \Big\|\sum_{k=0}^{j}\fT(X_k)\Big\|_E>rb_n\Big)\\
\ad\leq\|\mu P^{\nns}-\pi\|_{\rm var}
 +\PP_\pi\Big(\max_{0\leq j\leq\lf n\dl\rf}
 \Big\|\sum_{k=0}^{j}\fT(X_k)\Big\|_E>rb_n\Big).
\earray\eeq
The condition \eqref{chain-cond-mu} implies that for sufficiently large $n$, $\|\mu P^n -\pi\|_{\text{var}} \leq \varsigma^{-n}$. It yields
\beq{log-geo}
\limsup_{n\to \infty}\frac{n}{b_n^2}\log\|\mu P^{\nns}-\pi\|_{\text{var}}=-\infty.
\eeq
for every $s\in(0,1]$. When $s=0$, we do not compare with $\pi$, and we apply the decomposition below under $\wdt\PP_{\wdt\mu}$ in place of $\wdt\PP_{\wdt\pi}$. The three terms are then treated in the same way.
Focusing on the second term of the last line in \eqref{M-diff}, we have
\bea\ad
\PP_\pi \Big(\max_{0\leq j\leq
\lf n\dl \rf} \Big\|\sum_{k=0}^j \fT(X_k)\Big\|_E >r b_n \Big) =\wdt \PP_{\wdt \pi} \Big(\max_{0\leq j\leq
\lf n\dl \rf} \Big\|\sum_{k=0}^j \fT(X_k)\Big\|_E > r b_n \Big) \\
\aad = \wdt \PP_{\wdt \pi} \Big( \max_{0\leq j \leq
\lf n\dl \rf} \Big\| S_j^{\rm ini} + \sum_{k=1}^{\mathfrak v(j)-1} \xi_{k}(\fT)+ S_j^{\rm ter} \Big\|_E  >r  b_n \Big)  \\
\aad \leq \wdt \PP_{\wdt \pi} \Big(\sum_{k=0}^{m\tau(0)+m-1} \|\fT(X_k)\|_E >\frac{r b_n }{3}\Big) +\wdt \PP_{\wdt \pi} \Big(\max_{0\leq j \leq
\lf n\dl \rf} \Big\| \sum_{k=1}^{\mathfrak v(j)-1} \xi_k(\fT) \Big\|_E >\frac{r b_n}{3}\Big) \\
\aad \quad +\wdt \PP_{\wdt \pi}\Big( \max_{0\leq j \leq
\lf n \dl \rf}  \Big\| \sum_{m(\ell(j)+1) \leq k \leq j} \fT(X_k) \Big\|_E  > \frac{r b_n}{3}\Big)
\eea
The first and third terms can be dealt with similarly to \eqref{neg-1} and \eqref{neg-2}. For the second term, $\mathfrak v(j)-1\leq j$, and (2.25)--(2.28) in the proof of \cite[Chapter~IV, Theorem~2.1]{Che99}, which combine the Ottaviani-type inequality \cite[Chapter~IV, Lemma~2.3]{Che99} for $1$-dependent sequences with the MDP upper bound for $\{\xi_k\}$, give
\bea\ad
\limsup_{n\to\infty}\frac{n}{b_n^2}\log\wdt\PP_{\wdt\pi}\Big(\max_{0\leq j\leq \lf n\dl \rf
}\Big\|\sum_{k=1}^{j}\xi_k(\fT)\Big\|_E>\frac{rb_n}{3}\Big)\leq-\frac{\wdt\pi(\mathfrak a)}{m\dl}\inf_{\|z\|_E\geq r/9}J_\fT(z).
\eea
Since $\inf_{\|z\|_E\geq r/9}J_\fT(z)>0$, the right-hand side tends to $-\infty$ as $\dl\to0$, uniformly in $s$. Thus, by \cite[Lemma~A.1]{DGW99} and \cite[Section~4.6]{DZ09}, the LDP holds with rate $\sup_{\CP}J_\fT^{\CP} = \wdt J_\fT$, which is our desired rate function. The proof is completed. \qed

\subsection{Proofs of Lemmas~\ref{lem:neg-D} and~\ref{lem:neg-A}}\label{app:additive-holder}
This subsection proves Lemmas~\ref{lem:neg-D} and~\ref{lem:neg-A}. We first explain why the semigroup has to be separated from the environment before the regeneration argument is applied. For $s\in[t_k,t_{k+1})$, the integrand of $D_\e^N$ in \eqref{decomp-la} splits as
\beqa{D-integrand-split}
 \ad e^{(t-s)A}F_c(s,\zeta_{s/\e})-e^{(t-t_k)A}F_c(t_k,\zeta_{s/\e})\\
 \aad =\big[e^{(t-s)A}-e^{(t-t_k)A}\big]F_c(s,\zeta_{s/\e})
 +e^{(t-t_k)A}\big[F_c(s,\zeta_{s/\e})-F_c(t_k,\zeta_{s/\e})\big].
\eeqa
The second difference is small because $F_c$ is $\theta$-H\"older continuous in time by \eqref{Fc-H-bounds}. The first difference is small only when it acts on a smooth vector: \eqref{e-I} gives $|(e^{(t-s)A}-e^{(t-t_k)A})x|_H\leq K_\chi(s-t_k)^\chi\|x\|_\chi$, and this requires $\sup_{s,z}\|F_c(s,z)\|_\chi<\infty$ for some $\chi>0$. Such spatial regularity of $F_c$ is not assumed. We therefore never apply a semigroup increment to $F_c$ itself. Instead we integrate by parts in time, so that the semigroup acts on increments of the additive functional
\beq{Ga-e}
 \Ga^\e(t):=\feh\int_0^tF_c(s,\zeta_{s/\e})\,ds,\qquad t\in[0,1].
\eeq
The price is the singular factor $\|Ae^{rA}\|_{\CL(H)}\leq K_Ar^{-1}$, where $K_A$ denotes the constant in \eqref{frac-A} with $\delta=1$. This factor is integrable against $|\Ga^\e(s)-\Ga^\e(t)|_H$ as soon as $\Ga^\e$ is H\"older continuous. For $\sigma\in(0,1)$ and $g:[0,1]\to H$, let
\[
 [g]_{\sigma}:=\sup_{0\le s<t\le1}\frac{|g(t)-g(s)|_H}{(t-s)^{\sigma}}
\]
be the $\sigma$-H\"older seminorm. The proof has three ingredients:
\begin{itemize}
\item[(a)] A deterministic integration by parts bounds the semigroup errors by $[\Ga^\e]_{\sigma}$ (\lemref{lem:conv-holder}).
\item[(b)] The regeneration argument gives an exponential bound for a single increment of an additive functional (\lemref{lem:holder-increment}). Applied to $F_c(s,\cdot)-F_c(t_k,\cdot)$ on $[t_k,t_{k+1}]$, it controls the second term of \eqref{D-integrand-split}, which contains no semigroup inside the time integral.
\item[(c)] A dyadic chaining argument turns (b) into a Gaussian tail bound for $[\Ga^\e]_{\sigma}$ (\lemref{lem:additive-holder}).
\end{itemize}

\begin{lem}\label{lem:conv-holder}
Let $\sigma\in(0,1)$, and let $g:[0,1]\to H$ be Lipschitz continuous.
\begin{itemize}
\item[(a)] For $0\le u\le t\le1$,
\beqa{holder-V}\ad
 \int_u^te^{(t-s)A}\dot g(s)\,ds
 =e^{(t-u)A}[g(t)-g(u)]+\int_u^tAe^{(t-s)A}[g(s)-g(t)]\,ds,
\eeqa
and consequently
$\big|\int_u^te^{(t-s)A}\dot g(s)\,ds\big|_H\le(1+K_A/\sigma)(t-u)^{\sigma}[g]_{\sigma}$.
\item[(b)] Let $N\geq2$, $\Dl=N^{-1}$ and $t_k=k\Dl$. For $t\in[t_j,t_{j+1})$,
\beqa{R-N-IBP}\ad
 \sum_{k=0}^{j-1}\int_{t_k}^{t_{k+1}}\big[e^{(t-s)A}-e^{(t-t_k)A}\big]\dot g(s)\,ds
 =\sum_{k=0}^{j-1}\int_{t_k}^{t_{k+1}}Ae^{(t-s)A}[g(s)-g(t_{k+1})]\,ds,
\eeqa
and the norm of this sum is at most $K_A(1+\sigma^{-1}+\log N)N^{-\sigma}[g]_{\sigma}$.
\end{itemize}
\end{lem}

\begin{proof}
(a) Since $e^{tA}$ is analytic, $s\mapsto e^{(t-s)A}$ is differentiable in $\CL(H)$ on $[u,t)$ with derivative $-Ae^{(t-s)A}$, and $g$ is absolutely continuous. Hence, for almost every $s\in[u,t)$,
\[
 \frac{d}{ds}\big(e^{(t-s)A}[g(s)-g(t)]\big)
 =-Ae^{(t-s)A}[g(s)-g(t)]+e^{(t-s)A}\dot g(s).
\]
Integrate over $[u,t']$ and let $t'\uparrow t$. The boundary term at $t'$ is bounded by $|g(t')-g(t)|_H\to0$, and
$|Ae^{(t-s)A}[g(s)-g(t)]|_H\le K_A[g]_{\sigma}(t-s)^{\sigma-1}$ is integrable on $[u,t)$. This proves \eqref{holder-V}. Since $e^{tA}$ is a contraction, the first term of \eqref{holder-V} is bounded by $(t-u)^{\sigma}[g]_{\sigma}$, and the second by $K_A[g]_{\sigma}\int_u^t(t-s)^{\sigma-1}ds=K_A\sigma^{-1}(t-u)^{\sigma}[g]_{\sigma}$.

(b) For $k\le j-1$ and $s\in[t_k,t_{k+1}]$, let $S_k(s):=e^{(t-s)A}-e^{(t-t_k)A}$. Then $S_k(t_k)=0$ and $\frac{d}{ds}S_k(s)=-Ae^{(t-s)A}$ for $s<t$. The product rule, applied to $S_k(s)[g(s)-g(t_{k+1})]$ as in (a), gives \eqref{R-N-IBP}; when $k=j-1$ and $t=t_j$ the limit $t'\uparrow t$ is taken as in (a). Since $t\ge t_j$, the term $k=j-1$ is bounded by
$K_A[g]_{\sigma}\int_{t_{j-1}}^{t_j}(t_j-s)^{\sigma-1}ds=K_A\sigma^{-1}\Dl^{\sigma}[g]_{\sigma}$. For $k\le j-2$, we have $t-s\ge(j-k-1)\Dl$ on $[t_k,t_{k+1}]$, so the corresponding term is bounded by
\[
 \frac{K_A[g]_{\sigma}}{(j-k-1)\Dl}\int_{t_k}^{t_{k+1}}(t_{k+1}-s)^{\sigma} ds
 \le\frac{K_A\Dl^{\sigma}[g]_{\sigma}}{j-k-1}.
\]
Summing over $k$ and using $\sum_{l=1}^{N}l^{-1}\le1+\log N$ proves the bound.
\end{proof}

\subsubsection{Regeneration split chain structure}

The probabilistic estimates use the regeneration split chain argument of Lemma~\ref{lem:B-mdp}, following Guillin \cite[Theorem~2, pp.~298--303]{Gui01}. Let $\zeta^{(\ell)}=\{\zeta_{\ell+s}\}_{s\in[0,1]}$, $\ell\ge0$, be the unit-time path chain of $\zeta$ from the proof of \corref{cor:zeta-path}, with state space $\wdh\SE$ and transition probability $\wdh P$. By Down et al.\
\cite[Theorem~5.3]{DMT95}, the chain $\{\zeta^{(\ell)}\}_{\ell \geq 0}$ is
geometrically ergodic. For simplicity, we assume that the Markov chain $\{\zeta^{(\ell)}\}_{\ell \geq 0}$ has an atom: there exist a measurable set $\mathfrak a$ and a probability measure $\nu_{\mathfrak a}$ such that $\wdh P(z,\cdot)=\nu_{\mathfrak a}$ for all $z\in \mathfrak a$. In the general case, the chain has a small set, from which one builds the associated split chain with an atom; see \appref{app:Bmdp} and Guillin \cite[Remark~2.9]{Gui01}.

Define the regeneration times
\bea
\varrho(0)
\ad =\inf\big\{n\geq 0: \zeta^{(n)}\in \mathfrak a \big\},\quad
\varrho(i+1) =\inf\big\{n > \varrho(i):  \zeta^{(n)}\in \mathfrak a \big\},
\eea
the number of visits to $\mathfrak a$ before time $n$,
$N_{\mathfrak a}(n):=\sum_{\ell=0}^{n-1}\indi_{\{\zeta^{(\ell)}\in \mathfrak a\}}$. For $i\geq0$, let $\CG_i$ be the $\sigma$-field of the chain stopped at $\varrho(i)$, that is, $\CB \in\CG_i$ if and only if $\CB\cap\{\varrho(i)=n\}\in\sg(\zeta^{(0)},\dots,\zeta^{(n)})$ for every $n\ge0$. The initial law of $\zeta$ will be $\mu$ or $\pi$; in both cases the chain has the following properties, see \lemref{lem:reg-geo} and \cite{Num84,AC98,Che99}.
\begin{itemize}
\item[(i)] $\{\varrho(i)-\varrho(i-1)\}_{i\ge1}$ are i.i.d.\ random variables with law $\operatorname{Law}_{\nu_{\mathfrak a}}(\varrho(0)+1)$. Moreover, given $\CG_{i-1}$, the chain $\{\zeta^{(\varrho(i-1)+\ell+1)}\}_{\ell \geq 0}$ has law $\PP_{\nu_{\mathfrak a}}$; in particular, $\varrho(i)-\varrho(i-1)$ is independent of $\CG_{i-1}$.
\item[(ii)] There is $\kappa_\varrho>0$ such that
\beq{exp-varrho-0}
K_\varrho^0:=\max\big\{\EE_\mu e^{\kappa_\varrho \varrho(0)},\,\EE_\pi e^{\kappa_\varrho \varrho(0)}\big\}<\infty,
\qquad
K_\nu:=\EE_{\nu_{\mathfrak a}}e^{\kappa_\varrho(\varrho(0)+1)}<\infty .
\eeq
For $\mu$ this follows from \eqref{cond-mu} and geometric ergodicity, and for $\nu_{\mathfrak a}$ because $\varrho(0)+1$ under $\PP_{\nu_{\mathfrak a}}$ is the return time to the atom. With $K_\varrho:=\kappa_\varrho^{-1} K_\nu^{1/2}$ we have $\EE_{\nu_{\mathfrak a}}[\varrho(0)+1]^p \leq p!\, K_\varrho^p$ for $p\geq 2$, and $\PP(\varrho(i)-\varrho(i-1)>v)\le K_\nu e^{-\kappa_\varrho v}$ for $v\ge0$.
\item[(iii)] Kac's formula \cite[Theorem 10.0.1]{MT12}: if $f$ is bounded and $\EE_\pi f(\zeta^{(0)})=0$, then
\beq{EE-f}
\EE_{\nu_{\mathfrak a}}\Big[\sum_{\ell=0}^{\varrho(0)} f(\zeta^{(\ell)})\Big]=0.
\eeq
\end{itemize}
Constants depending only on $\kappa_\varrho$, $K_\nu$ and $K_\varrho^0$ are called regeneration constants.

\subsubsection{One increment of an additive functional}

Let $b:[0,1]\times\SE\to H$ be jointly measurable with $\pi(b(s,\cdot))=0$ for every $s$. Fix an interval $[a,t]\subset[0,1]$, let $\delta:=t-a$, and assume that for some constants $\bar K_b$, $K_b$ and $\theta\in(1/2,1]$ and any $ s,s'\in[a,t],\ z\in\SE$, we have
\beqa{b-piecewise-holder}
 |b(s,z)|_H\le \bar K_b,\qquad |b(s,z)-b(s',z)|_H\le K_b|s-s'|^\theta.
\eeqa

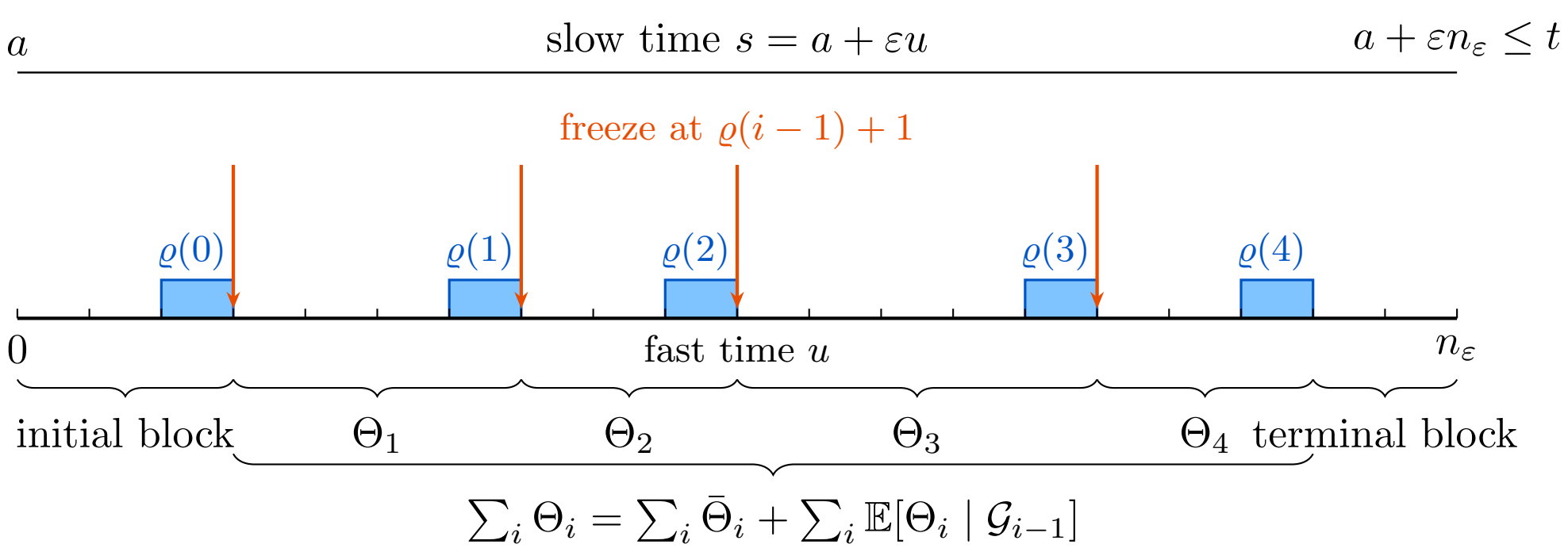
\begin{figure}[ht]
\centering
\definecolor{RegenAtom}{RGB}{0,84,196}
\definecolor{RegenAtomFill}{RGB}{128,196,255}
\definecolor{RegenFreeze}{RGB}{232,76,0}
\begin{tikzpicture}[x=0.6cm,y=1cm,font=\footnotesize]
 \draw[black,line width=0.45pt] (0,2.05) -- (20,2.05);
 \node[above] at (0,2.05) {$a$};
 \node[above] at (10,2.05) {slow time $s=a+\e u$};
 \node[above] at (20,2.05) {$a+\e n_\e\le t$};

 \foreach \r/\lab in {2/{\varrho(0)},6/{\varrho(1)},9/{\varrho(2)},
                     14/{\varrho(3)},17/{\varrho(4)}} {
   \path[fill=RegenAtomFill,draw=RegenAtom,line width=0.6pt]
     (\r,0) rectangle (\r+1,0.32);
   \node[above,text=RegenAtom,font=\scriptsize,inner sep=1pt]
     at (\r+0.42,0.36) {$\lab$};
 }
 \draw[black,thick] (0,0) -- (20,0);
 \foreach \l in {0,...,20} {\draw[black] (\l,0) -- (\l,0.08);}
 \node[below] at (0,0) {$0$};
 \node[below,font=\scriptsize] at (10,0) {fast time $u$};
 \node[below] at (20,0) {$n_\e$};

 \node[text=RegenFreeze,font=\scriptsize] at (10,1.57)
   {freeze at $\varrho(i-1)+1$};
 \foreach \s in {3,7,10,15} {
   \draw[RegenFreeze,line width=0.8pt,-{Stealth[length=4pt,width=3.2pt]}]
     (\s,1.28) -- (\s,0.08);
 }

 \draw[decorate,decoration={brace,mirror,amplitude=4pt}]
   (0,-0.51) -- (3,-0.51) node[midway,below=5pt] {initial block};
 \draw[decorate,decoration={brace,mirror,amplitude=4pt}]
   (3,-0.51) -- (7,-0.51) node[midway,below=5pt] {$\Theta_1$};
 \draw[decorate,decoration={brace,mirror,amplitude=4pt}]
   (7,-0.51) -- (10,-0.51) node[midway,below=5pt] {$\Theta_2$};
 \draw[decorate,decoration={brace,mirror,amplitude=4pt}]
   (10,-0.51) -- (15,-0.51) node[midway,below=5pt] {$\Theta_3$};
 \draw[decorate,decoration={brace,mirror,amplitude=4pt}]
   (15,-0.51) -- (18,-0.51) node[midway,below=5pt] {$\Theta_4$};
 \draw[decorate,decoration={brace,mirror,amplitude=4pt}]
   (18,-0.51) -- (20,-0.51) node[midway,below=5pt] {terminal block};
 \draw[decorate,decoration={brace,mirror,amplitude=5pt}]
   (3,-1.13) -- (18,-1.13)
   node[midway,below=6pt]
   {$\sum_i\Theta_i=\sum_i\bar\Theta_i+\sum_i\EE[\Theta_i\mid\CG_{i-1}]$};
\end{tikzpicture}
\caption{Regeneration blocks in \lemref{lem:holder-increment} in the atom case.
A shaded cell marks a visit to the atom at index $\varrho(i)$ and is the
last cell of its block. The complete block $\Theta_i$ starts at
$\varrho(i-1)+1$, where the arrow marks the freezing point.
The initial, complete, and terminal sums correspond to the tail
probabilities $\mathrm{I}_1$, $\mathrm{I}_2$, and $\mathrm{I}_3$, respectively.}
\label{fig:blocks}
\end{figure}

\begin{lem}[Exponential bound for one increment]\label{lem:holder-increment}
There are regeneration constants $\kappa_{\rm reg}>0$ and $K_{\rm reg}\geq1$ such that the following holds. Let $0<\e\leq1$, $r>0$, and assume
\beqa{holder-bias-small}\ad
 \frac{\sqe}{h(\e)}\big(\bar K_b+2K_bK_\varrho^2\,\delta\,\e^{\theta-1}\big)\le\frac r8 .
\eeqa
Then we have
\beqa{holder-increment-bound}
\ad \PP_\mu\Big(\Big|\int_a^tb(s,\zeta_{s/\e})\,ds\Big|_H>r\sqe h(\e)\Big) \\
 \ad\le2\exp\Bigl(-\frac{\kappa_{\rm reg}\,r^2h^2(\e)}{\bar K_b^2\delta+\bar K_br\sqe h(\e)}\Bigr)
 +\frac{K_{\rm reg}}{\e}\exp\Bigl(-\frac{\kappa_\varrho rh(\e)}{4\bar K_b\sqe}\Bigr) +\indi_{\{a>0\}} K_\phi\mixrate^{a/\e}.
\eeqa
\end{lem}
The proof adapts the regeneration decomposition in \lemref{lem:B-mdp} to a time-dependent observable. Figure~\ref{fig:blocks} identifies the initial, complete, and terminal blocks.

\begin{proof}
Set $r_\e:=rh(\e)/\sqe$. In Steps~1--5 the dependence on $\e$ and on the interval $[a,t]$ is suppressed in the notation. The substitution $s=a+\e u$ gives
\[
 \int_a^tb(s,\zeta_{s/\e})\,ds=\e\int_0^{\delta/\e}\Xi(u,\zeta_{a/\e+u})\,du,
 \qquad \Xi(u,z):=b(a+\e u,z),
\]
so the event in \eqref{holder-increment-bound} is $\{|\int_0^{\delta/\e}\Xi(u,\zeta_{a/\e+u})du|_H>r_\e\}$. It is determined by the shifted path $(\zeta_{a/\e+u})_{u\ge0}$. If $a>0$, the Markov property at time $a/\e$, \eqref{exp-erg} and \eqref{cond-mu} show that the law of the shifted path under $\PP_\mu$ differs in total variation from the law of $(\zeta_u)_{u\ge0}$ under $\PP_\pi$ by at most $\|\mu P_{a/\e}-\pi\|_{\rm var}\le K_\phi\mixrate^{a/\e}$. Hence
\beq{holder-fresh}
 \PP_\mu\Big(\Big|\int_a^tb(s,\zeta_{s/\e})\,ds\Big|_H>r\sqe h(\e)\Big)
 \!\le\!\PP_{\nu_0}\Big(\Big|\int_0^{\delta/\e}\Xi(u,\zeta_u)\,du\Big|_H>r_\e\Big)\!+\!\indi_{\{a>0\}}K_\phi\mixrate^{a/\e},
\eeq
where $\nu_0:=\pi$ if $a>0$ and $\nu_0:=\mu$ if $a=0$. From now on $\zeta$ has initial law $\nu_0$, and probabilities and expectations are taken under $\PP_{\nu_0}$.

\emph{Step 1: reduction to the path chain.} Let $n_\e:=\lfloor\delta/\e\rfloor$. The integral over $[n_\e,\delta/\e]$ has norm at most $\bar K_b\le r_\e/8$ by \eqref{holder-bias-small}. If $n_\e=0$, the proof is complete. Otherwise, for $\omega\in\wdh\SE$ we set
\beqa{frak-f}
f_\ell(\omega):=\int_0^1 \Xi\big(\ell\wedge(n_\e-1)+w,\,\omega(w)\big)\,dw, \quad \ell\geq0,
\eeqa
so that $\int_0^{n_\e}\Xi(u,\zeta_u)du=\sum_{\ell=0}^{n_\e-1}f_\ell(\zeta^{(\ell)})$. For $\ell\le n_\e-1$ and $w\in[0,1)$ the time $a+\e(\ell+w)$ lies in $[a,t]$. Therefore \eqref{b-piecewise-holder} gives
\beq{diff-ff}
 \|f_\ell\|_\infty\le \bar K_b,\qquad
 |f_\ell(\omega)-f_j(\omega)|_H \leq K_b(\e|\ell-j|)^{\theta}, \quad \ell,j \geq 0,
\eeq
and, by Fubini's theorem and the stationarity of $\zeta$ under $\PP_\pi$, we have \bea \EE_\pi f_\ell(\zeta^{(0)})=\int_0^1\pi(\Xi((\ell\wedge(n_\e-1))+w,\cdot))\,dw=0.\eea

\emph{Step 2: regeneration decomposition.} Split the sum at the first and last indices $\varrho(i)<n_\e$:
\[
 \sum_{\ell=0}^{n_\e-1}f_\ell(\zeta^{(\ell)})
 =\sum_{\ell=0}^{\varrho(0)\wedge (n_\e-1)}f_\ell(\zeta^{(\ell)})
 +\sum_{i=1}^{N_{\mathfrak a}(n_\e)-1}\Theta_i
 +\sum_{\varrho(N_{\mathfrak a}(n_\e)-1)<\ell\leq n_\e-1}f_\ell(\zeta^{(\ell)}).
\]
where $\Theta_i:=\sum_{\ell=\varrho(i-1)+1}^{\varrho(i)}f_\ell(\zeta^{(\ell)})$.
If $N_{\mathfrak a}(n_\e)\le1$, the middle sum is empty, and if $N_{\mathfrak a}(n_\e)=0$, the last sum is empty as well. Since the remainder in Step~1 is at most $r_\e/8$,
\[
 \PP_{\nu_0}\Big(\Big|\int_0^{\delta/\e}\Xi(u,\zeta_u)\,du\Big|_H>r_\e\Big)
 \le\mathrm{I}_1+\mathrm{I}_2+\mathrm{I}_3,
\]
where $\mathrm{I}_1$ and $\mathrm{I}_3$ are the probabilities that the first and the last sum exceed $r_\e/4$ in norm, and $\mathrm{I}_2$ is the probability that the middle sum exceeds $3r_\e/8$ in norm.

The first sum is bounded in norm by $(\varrho(0)+1)\bar K_b$. By (ii) and Chebyshev's inequality,
\[
 \mathrm{I}_1\le\PP_{\nu_0}\Big(\varrho(0)+1>\frac{r_\e}{4\bar K_b}\Big)\le e^{\kappa_\varrho}K_\varrho^0\exp\Big(-\frac{\kappa_\varrho r_\e}{4\bar K_b}\Big).
\]
If the last sum is not empty, it has fewer than $\varrho(N_{\mathfrak a}(n_\e))-\varrho(N_{\mathfrak a}(n_\e)-1)$ terms, and $N_{\mathfrak a}(n_\e)\le n_\e$. By (i) and (ii),
\[
 \mathrm{I}_3\le\sum_{i=1}^{n_\e}\PP_{\nu_0}\Big(\varrho(i)-\varrho(i-1)>\frac{r_\e}{4\bar K_b}\Big)\le n_\e K_\nu\exp\Big(-\frac{\kappa_\varrho r_\e}{4\bar K_b}\Big).
\]
For $\mathrm{I}_2$, the blocks $\Theta_i$ are defined for every $i\ge1$, and $\Theta_i$ is $\CG_i$-measurable. Let $\bar\Theta_i:=\Theta_i-\EE_{\nu_0}[\Theta_i\mid\CG_{i-1}]$. Then $\{\bar\Theta_i\}$ are martingale differences with respect to $\{\CG_i\}$, and, since $N_{\mathfrak a}(n_\e)-1\le n_\e$, we have
\beq{sum-dga}
\Big|\sum_{i=1}^{N_{\mathfrak a}(n_\e)-1}\Theta_i\Big|_H
\le\max_{0\le j\le n_\e}\Big|\sum_{i=1}^j\bar\Theta_i\Big|_H
+\sum_{i=1}^{n_\e}\big|\EE_{\nu_0}[\Theta_i\mid\CG_{i-1}]\big|_H .
\eeq

\emph{Step 3: the conditional bias.} By (i), we have
$
 \EE_{\nu_0}[\Theta_i\mid\CG_{i-1}]=\mathfrak b(\varrho(i-1))$ and
 \bea\ad
 \mathfrak b(q):=\EE_{\nu_{\mathfrak a}}\Big[\sum_{\ell=0}^{\varrho(0)}f_{q+1+\ell}(\zeta^{(\ell)})\Big],\quad q\in\NN_+\cup\{0\}.
\eea
Write $f_{q+1+\ell}=f_{q+1}+(f_{q+1+\ell}-f_{q+1})$, that is, freeze the observable at the beginning $q+1$ of the block. The frozen part has zero mean by \eqref{EE-f} and $\EE_\pi f_{q+1}(\zeta^{(0)})=0$. Equation~\eqref{diff-ff} and the regeneration moment bound (ii) control the remainder. As $\theta\le1$,
\beq{cond-The}
 |\mathfrak b(q)|_H\le\EE_{\nu_{\mathfrak a}}\Big[\sum_{\ell=0}^{\varrho(0)}K_b(\e \ell)^\theta\Big]
 \le K_b\e^\theta\,\EE_{\nu_{\mathfrak a}}(\varrho(0)+1)^2\le2K_bK_\varrho^2\e^\theta .
\eeq
Summing over $i\le n_\e\le\delta/\e$ and using \eqref{holder-bias-small}, we obtain
$\sum_{i=1}^{n_\e}|\EE_{\nu_0}[\Theta_i\mid\CG_{i-1}]|_H\le2K_bK_\varrho^2\delta\e^{\theta-1}\le r_\e/8$. By \eqref{sum-dga},
\beq{leq-CP4}
 \mathrm{I}_2\le\PP_{\nu_0}\Big(\max_{0\le j\le n_\e}\Big|\sum_{i=1}^j\bar\Theta_i\Big|_H>\frac{r_\e}4\Big).
\eeq

\emph{Step 4: the martingale part.} Since $|\Theta_i|_H\le \bar K_b(\varrho(i)-\varrho(i-1))$, the conditional Jensen inequality, (i) and (ii) give, for $p\ge2$,
\beq{cond-mom-The}
\EE_{\nu_0}\big[|\bar\Theta_i|_H^p\bigm|\CG_{i-1}\big]
\le2^p\bar K_b^p\,\EE_{\nu_{\mathfrak a}}(\varrho(0)+1)^p
\le\frac{p!}2\,2(2K_\varrho \bar K_b)^2(2K_\varrho \bar K_b)^{p-2}.
\eeq
Apply \eqref{Hil-Ber-tail} with $v_i^\beta=8K_\varrho^2\bar K_b^2$, $\kappa_0=2K_\varrho \bar K_b$, and $\kappa_v=8K_\varrho^2\bar K_b^2n_\e$. Setting $\kappa_{\rm reg}=1/(256(K_\varrho^2\vee K_\varrho))$ gives
\[
 \mathrm{I}_2\le2\exp\Big(-\frac{r_\e^2/16}{16K_\varrho^2\bar K_b^2n_\e+K_\varrho \bar K_br_\e}\Big)
 \le2\exp\Big(-\frac{\kappa_{\rm reg}\,r_\e^2}{\bar K_b^2\delta/\e+\bar K_br_\e}\Big).
\]

\emph{Step 5: conclusion.} The identity
\[
\frac{r_\e^2}{\bar K_b^2\delta/\e+\bar K_br_\e}
=\frac{r^2h^2(\e)}{\bar K_b^2\delta+\bar K_br\sqe h(\e)}
\]
together with $r_\e/\bar K_b=rh(\e)/(\bar K_b\sqe)$ and $n_\e\le\delta/\e\le1/\e$, gives the two exponents in \eqref{holder-increment-bound}. The prefactor of the second term is bounded by $K_{\rm reg}/\e$ with $K_{\rm reg}=e^{\kappa_\varrho}K_\varrho^0+K_\nu$. Combine the estimates for $\mathrm{I}_1$, $\mathrm{I}_2$, and $\mathrm{I}_3$ with \eqref{holder-fresh}.
\end{proof}

\begin{lem}{\rm (H\"older tail of $\Ga^\e$)}\label{lem:additive-holder}
Let $K$ be as in \eqref{Fc-H-bounds}. For every $\sigma\in(0,1/2)$ there is $\kappa_\sigma>0$, depending only on $\sigma$ and the regeneration constants, such that $\Ga^\e$ in \eqref{Ga-e} satisfies
\beq{additive-holder-tail}
 \limsup_{\e\to0}h^{-2}(\e)\log\PP_\mu
 \big([\Ga^\e]_{\sigma}>R\big)
 \le-\frac{\kappa_\sigma R^2}{K^2},\qquad R>0.
\eeq
\end{lem}

\begin{proof}
Fix $\sigma\in(0,1/2)$ and $R>0$, and recall that $\sqe h(\e)\to0$ by \eqref{dev-s}.

\emph{Step 1: dyadic increments.} For a continuous $g:[0,1]\to H$ and $p\ge0$, let
$d_p(g):=\max_{0\le j<2^p}2^{p\sigma}|g((j+1)2^{-p})-g(j2^{-p})|_H$. We claim that
\beqa{holder-dyadic}\ad
 [g]_{\sigma}\le K_\sigma\sup_{p\ge0}d_p(g),\qquad K_\sigma:=2\big(1+(1-2^{-\sigma})^{-1}\big).
\eeqa
Let $\bar d:=\sup_pd_p(g)$, let $s<t$ be dyadic rationals, and choose $p\ge0$ with $2^{-p-1}<t-s\le2^{-p}$. For $k\ge p+1$, let $s_k:=\lceil2^ks\rceil2^{-k}$ and $t_k':=\lfloor2^kt\rfloor2^{-k}$. Then $s_k=s$ and $t_k'=t$ for large $k$, and $s_k-s_{k+1}$ and $t_{k+1}'-t_k'$ belong to $\{0,2^{-k-1}\}$. Moreover, $t_{p+1}'-s_{p+1}$ is a multiple of $2^{-p-1}$ lying in $(-2^{-p-1},2^{-p}]$, so it equals $0$, $2^{-p-1}$ or $2^{-p}$. Going from $s$ to $s_{p+1}$, then to $t_{p+1}'$, then to $t$ along these points,
\[
 |g(t)-g(s)|_H
 \le2\bar d\,2^{-(p+1)\sigma}
 +2\sum_{k\ge p+1}\bar d\,2^{-(k+1)\sigma}
 \le K_\sigma\bar d\,2^{-(p+1)\sigma}\le K_\sigma\bar d\,(t-s)^{\sigma} .
\]
By continuity, this holds for all $s<t$, which proves \eqref{holder-dyadic}.

\emph{Step 2: fine scales are automatic.} By \eqref{holder-dyadic}, $\{[\Ga^\e]_{\sigma}>R\}\subset\bigcup_{p,j}\{|\Ga^\e((j+1)2^{-p})-\Ga^\e(j2^{-p})|_H>r_p\}$, where $r_p:=R\,2^{-p\sigma}/K_\sigma$. By \eqref{Fc-H-bounds}, $\Ga^\e$ is Lipschitz with constant $K/(\sqe h(\e))$, so
$|\Ga^\e((j+1)2^{-p})-\Ga^\e(j2^{-p})|_H\le K\,2^{-p}/(\sqe h(\e))\le r_p$ whenever $2^{-p}\le\delta_\e$, where
$
 \delta_\e:=\big(R\sqe h(\e)/(K_\sigma K)\big)^{1/(1-\sigma)} .
$
Hence only the scales $2^{-p}>\delta_\e$ remain:
\beqa{holder-coarse-union}\ad
 \{[\Ga^\e]_{\sigma}>R\}
 \subset\bigcup_{p:\,2^{-p}>\delta_\e}\ \bigcup_{0\le j<2^p}
 \big\{|\Ga^\e((j+1)2^{-p})-\Ga^\e(j2^{-p})|_H>r_p\big\},
\eeqa
and the number of pairs $(p,j)$ in \eqref{holder-coarse-union} is less than $2/\delta_\e$.

\emph{Step 3: coarse scales.} Let $\beta:=(1-2\sigma)/(1-\sigma)\in(0,1)$. By the definition of $\delta_\e$,
\beqa{holder-fine-scales}\ad
 \frac{\delta_\e^{\sigma}}{\sqe h(\e)}=K_\sigma'(\sqe h(\e))^{-\beta},\qquad
 \frac{\delta_\e}{\e h^2(\e)}=K_\sigma''(\sqe h(\e))^{-\beta},
\eeqa
with $K_\sigma',K_\sigma''>0$ depending only on $\sigma$, $R$ and $K$, and $(\sqe h(\e))^{-\beta}\to\infty$. We apply \lemref{lem:holder-increment} with $b=F_c$, $\bar K_b=K_b=K$, $[a,t]=[j2^{-p},(j+1)2^{-p}]$, $\delta=2^{-p}$ and $r=r_p$, for $2^{-p}>\delta_\e$. Condition \eqref{holder-bias-small} holds for small $\e$, uniformly in $(p,j)$: indeed, $\sqe K/(h(\e)r_p)\le K_\sigma K\sqe/(h(\e)R\delta_\e^{\sigma})\to0$ by \eqref{holder-fine-scales}, and
\[
 \frac{\sqe}{h(\e)r_p}\,2KK_\varrho^2\,2^{-p}\,\e^{\theta-1}
 =\frac{2K_\sigma KK_\varrho^2}{R}\,2^{-p(1-\sigma)}\,\frac{\e^{\theta-1/2}}{h(\e)}
 \le\frac{2K_\sigma KK_\varrho^2}{R}\,\frac{\e^{\theta-1/2}}{h(\e)}\to0,
\]
because $\theta>1/2$ and $h(\e)\to\infty$. We now evaluate the three terms of \eqref{holder-increment-bound}.
\begin{itemize}
\item Since $2^{-p}>\delta_\e$, we have $r_p\sqe h(\e)\le K\,2^{-p}$. Hence
\[
 \frac{\kappa_{\rm reg}\,r_p^2h^2(\e)}{K^2\,2^{-p}+Kr_p\sqe h(\e)}\ge\frac{\kappa_{\rm reg}\,r_p^2h^2(\e)}{2K^2\,2^{-p}}=\Lambda_\e2^{p(1-2\sigma)},
 \qquad \Lambda_\e:=\frac{\kappa_{\rm reg}R^2h^2(\e)}{2K_\sigma^2K^2} .
\]
\item Since $r_p\ge R\delta_\e^{\sigma}/K_\sigma$, \eqref{holder-fine-scales} gives
\[
\frac{\kappa_\varrho r_ph(\e)}{4K\sqe}
\ge \kappa_\sigma'h^2(\e)(\sqe h(\e))^{-\beta},
\qquad \kappa_\sigma'=\frac{\kappa_\varrho RK_\sigma'}{4K_\sigma K}.
\]
\item For $j\ge1$, $j2^{-p}/\e\ge\delta_\e/\e=K_\sigma''h^2(\e)(\sqe h(\e))^{-\beta}$.
\end{itemize}
Summing \eqref{holder-increment-bound} over the pairs in \eqref{holder-coarse-union} yields
\beqa{holder-dyadic-tail}
\ad \PP_\mu\big([\Ga^\e]_{\sigma}>R\big) \le2\sum_{p\ge0}2^p\exp\big(-\Lambda_\e2^{p(1-2\sigma)}\big)\\
 \aad+\frac2{\delta_\e}\Big[\frac{K_{\rm reg}}{\e}\exp\big(-\kappa_\sigma'h^2(\e)(\sqe h(\e))^{-\beta}\big)
 +K_\phi\exp\big(K_\sigma''\log\mixrate\,h^2(\e)(\sqe h(\e))^{-\beta}\big)\Big].
\eeqa
For $\lambda\ge1$ and $u\ge1$, $\lambda u\ge\lambda+u-1$. Hence, once $\Lambda_\e\ge1$, the first line is at most $2\bar K_\sigma e^{-\Lambda_\e}$ with $\bar K_\sigma:=\sum_{p\ge0}2^p\exp(1-2^{p(1-2\sigma)})<\infty$, and its contribution to $h^{-2}(\e)\log\PP_\mu(\cdot)$ is at most $-\kappa_{\rm reg}R^2/(2K_\sigma^2K^2)+o(1)$. In the second line, $\delta_\e\ge(R\sqe/(K_\sigma K))^{1/(1-\sigma)}$ if $h(\e)\ge1$, so for a constant $\wdt K_\sigma$ depending only on $\sigma$, $R$, $K$ and $K_{\rm reg}$ we have
\[
 h^{-2}(\e)\log\frac{2K_{\rm reg}}{\delta_\e\e}\le \wdt K_\sigma\,\frac{1+\log(1/\e)}{h^2(\e)}=o\big((\sqe h(\e))^{-\beta}\big),
\]
because $(\sqe h(\e))^{-\beta}h^2(\e)=\e^{-\beta/2}h^{2-\beta}(\e)\ge\e^{-\beta/2}$,
whereas both exponents divided by $h^2(\e)$ are at most $-\min\{\kappa_\sigma',K_\sigma''|\log\mixrate|\}(\sqe h(\e))^{-\beta}\to-\infty$. Thus the second line contributes $-\infty$ to the logarithmic upper limit. This is \eqref{additive-holder-tail} with $\kappa_\sigma:=\kappa_{\rm reg}/(2K_\sigma^2)$.
\end{proof}

\subsubsection{Proof of \lemref{lem:neg-D}}\label{app:neg-D}
Fix $\sigma\in(0,1/2)$ and $N\ge2$, and let $\Dl=N^{-1}$ and $t_k=k\Dl$. For $t\in[t_j,t_{j+1})$, the splitting \eqref{D-integrand-split} gives $D_\e^N(t)=D_1^N(t)+D_2^N(t)$ with
\bea
 D_1^N(t)\ad:=\feh\sum_{k=0}^{j-1}\int_{t_k}^{t_{k+1}}\big[e^{(t-s)A}-e^{(t-t_k)A}\big]F_c(s,\zeta_{s/\e})\,ds,\\
 D_2^N(t)\ad:=\feh\sum_{k=0}^{j-1}e^{(t-t_k)A}\int_{t_k}^{t_{k+1}}\big[F_c(s,\zeta_{s/\e})-F_c(t_k,\zeta_{s/\e})\big]\,ds .
\eea
For $D_1^N$, since $\dot\Ga^\e(s)=(\sqe h(\e))^{-1} F_c(s,\zeta_{s/\e})$, \lemref{lem:conv-holder} (b) with $g=\Ga^\e$ gives
$\|D_1^N\|_\infty\le K_A(1+\sigma^{-1}+\log N)N^{-\sigma}[\Ga^\e]_{\sigma}$. By \lemref{lem:additive-holder},
\beq{D-semigroup-tail}
 \limsup_{\e\to0}h^{-2}(\e)\log\PP_\mu\Big(\|D_1^N\|_\infty>\frac r2\Big)
 \le-\frac{\kappa_\sigma r^2N^{2\sigma}}{4K_A^2K^2(1+\sigma^{-1}+\log N)^2}.
\eeq
For $D_2^N$, since $e^{tA}$ is a contraction,
\[
 \|D_2^N\|_\infty\le\feh\sum_{k=0}^{N-1}\Big|\int_{t_k}^{t_{k+1}}\big[F_c(s,\zeta_{s/\e})-F_c(t_k,\zeta_{s/\e})\big]\,ds\Big|_H ,
\]
and therefore
\[
 \PP_\mu\Big(\|D_2^N\|_\infty>\frac r2\Big)
 \le\sum_{k=0}^{N-1}\PP_\mu\Big(\Big|\int_{t_k}^{t_{k+1}}\big[F_c(s,\zeta_{s/\e})-F_c(t_k,\zeta_{s/\e})\big]\,ds\Big|_H>\frac{r}{2N}\sqe h(\e)\Big).
\]
Fix $k$ and apply \lemref{lem:holder-increment} with $b(s,z):=F_c(s,z)-F_c(t_k,z)$, $[a,t]=[t_k,t_{k+1}]$, $\delta=\Dl$, and $r/(2N)$ in place of $r$. By \lemref{lem:Fc-U}, $b$ is centered, and \eqref{b-piecewise-holder} holds with $\bar K_b=K N^{-\theta}$, $K_b=K$. Condition \eqref{holder-bias-small} reads
$\sqe h^{-1}(\e)\big(KN^{-\theta}+2KK_\varrho^2N^{-1}\e^{\theta-1}\big)\le r/(16N)$. It holds for small $\e$, because $\sqe/h(\e)\to0$ and $\e^{\theta-1/2}/h(\e)\to0$ for $\theta>1/2$. The three terms of \eqref{holder-increment-bound} behave as follows.
\begin{itemize}
\item The first exponent is
\[
 \frac{\kappa_{\rm reg}\,(r/2N)^2h^2(\e)}{K^2N^{-2\theta-1}+KN^{-\theta}(r/2N)\sqe h(\e)}
 =\frac{\kappa_{\rm reg}\,r^2N^{2\theta-1}}{4K^2}\,h^2(\e)\,(1+o(1)), \;\text{ as } \sqe h(\e)
\to 0.
\]
\item The second exponent divided by $h^2(\e)$ is $\kappa_\varrho rN^{\theta-1}/(8K\sqe h(\e))\to\infty$, while $h^{-2}(\e)\log(K_{\rm reg}/\e)\cdot\sqe h(\e)\to0$.
\item For $k\ge1$, $h^{-2}(\e)\log(K_\phi\mixrate^{t_k/\e})\le h^{-2}(\e)\log K_\phi+\log\mixrate/(N\e h^2(\e))\to-\infty$ due to the fact that $\e h^2(\e)\to0$.
\end{itemize}
Hence, after the finite union over $k$,
\beq{CP4-1}
 \limsup_{\e\to0}h^{-2}(\e)\log\PP_\mu\Big(\|D_2^N\|_\infty>\frac r2\Big)\le-\frac{\kappa_{\rm reg}\,r^2N^{2\theta-1}}{4K^2}.
\eeq

By \eqref{D-semigroup-tail} and \eqref{CP4-1},
\[
 \limsup_{\e\to0}h^{-2}(\e)\log
 \PP_\mu(\|D_\e^N\|_\infty>r)
 \le-\kappa'r^2\min\Bigl\{N^{2\theta-1},
 \frac{N^{2\sigma}}{(1+\sigma^{-1}+\log N)^2}\Bigr\}
\]
with $\kappa'>0$ independent of $N$. Since $\theta>1/2$ and $\sigma>0$, the right-hand side tends to $-\infty$ as $N\to\infty$. This proves \eqref{D}. \qed

\subsubsection{Proof of \lemref{lem:neg-A}}\label{app:lem-neg-A}
Fix $\sigma\in(0,1/2)$. For $t\in[t_j,t_{j+1})$, $A_\e^N(t)=\int_{t_j}^te^{(t-s)A}\dot\Ga^\e(s)\,ds$ with $\Ga^\e$ in \eqref{Ga-e}. By \lemref{lem:conv-holder}(a) with $u=t_j$,
$\|A_\e^N\|_\infty\le(1+K_A/\sigma)N^{-\sigma}[\Ga^\e]_{\sigma}$.
Here \eqref{holder-V} subtracts $\Ga^\e(t)$ inside the integral. Subtracting $\Ga^\e(t_j)$ instead leads to the integrand $Ae^{(t-s)A}[\Ga^\e(s)-\Ga^\e(t_j)]$, which is not integrable near $s=t$ without spatial regularity of $F_c$. By \lemref{lem:additive-holder},
\[
 \limsup_{\e\to0}h^{-2}(\e)\log
 \PP_\mu(\|A_\e^N\|_\infty>r)
 \le-\frac{\kappa_\sigma r^2N^{2\sigma}}{(1+K_A/\sigma)^2K^2}
 \xrightarrow[N\to\infty]{}-\infty,
 \qquad r>0.
\]
This proves \eqref{A}. \qed

\subsection{Proof of \lemref{lem:CI-CF}}\label{app:lem-CI-CF}
By \eqref{CF-time-holder}, $\|C_F(t)-C_F(s)\|_{\CL(H)}\le K|t-s|^\theta$; in particular, $C_F$ is continuous on $[0,1]$.
Let $K_C:=\max(1,\sup_{t\in[0,1]}\|C_F(t)\|_{\CL(H)})$.

We first prove \eqref{CF-leq}. It suffices to consider a subsequence for which $\sup_{N}\CI_N^{\rm hom}(\ga_N)<\infty$, otherwise it is trivial. Then for $\ga_N \in \mathscr{S}_N$,
\bea\ad
v^N_k= [\ga_N(t_{k+1})-e^{\Dl A}\ga_N(t_k)]/\Dl, \quad v^N(t)=v_k^N, \quad t_k\leq t <t_{k+1}.
\eea
As $e^{tA}$ is a contraction, $\|e^{\Dl A}C_F(t_k)e^{\Dl A^*}\|\le K_C$.
The definition of $L_C$ therefore gives
\[
L(e^{\Dl A}C_F(t_k)e^{\Dl A^*},v_k^N)
\ge\frac{|v_k^N|_H^2}{2K_C},
\qquad
\int_0^1|v^N(t)|_H^2\,dt
\le2K_C\CI_N^{\rm hom}(\ga_N).
\]
Thus, after passing to a subsequence, $v^N\rightharpoonup v$ in $L^2([0,1];H)$.

We now proceed to identify $v$. Let $\phi\in C^1([0,1]; H_m)$ for some $m$. We have
\bea\ad
\sum_{k=0}^{N-1}\Dl \qv{v^N_k, \phi(t_{k+1})}_H = \qv{\ga_N(1),\phi(1)}_H +\sum_{k=1}^{N-1} \qv{\ga_N(t_k), \phi(t_k)-e^{\Dl A^*} \phi(t_{k+1})}_H.
\eea
On $H_m$, we know that $A_m$ is bounded. Hence
\bea
[\phi(t_k)-e^{\Dl A^*} \phi(t_{k+1})]/\Dl \to -\dot{\phi}(t_k)-A^* \phi(t_k)
\eea
uniformly in $k$. Since $\ga_N \to \ga$ uniformly and $v^N \rightharpoonup v$, the preceding identity gives
\bea\ad
\int_0^1 \qv{v(t),\phi(t)}_H dt = \qv{\ga(1),\phi(1)}_H -\int_0^1 \qv{\ga(t), \dot{\phi}(t)+A^* \phi(t)}_H dt .
\eea
Working on the
coordinates, we have $\dot{\ga}_i+\al_i \ga_i = v_i$ and $\ga_i(0)=0$. Moreover,
\bea\ad
\int_0^1 |v_i|^2 dt =\int_0^1 (|\dot{\ga}_i|^2 +\al_i^2 |\ga_i|^2)dt+ \al_i |\ga_i(1)|^2.
\eea
Summing over $i$ implies that $\ga \in\SD_A$ and $v=\dot{\ga}-A\ga$.

Let $u\in C([0,1];H_m)$. The variational representation \eqref{IN} yields
\bea\ad
\CI_N^{\rm hom}(\ga_N) \geq \sum_{k=0}^{N-1} \Dl \big[\qv{v^N_k, u(t_k)}_H - \half \qv{e^{\Dl A} C_F(t_k) e^{\Dl A^*} u(t_k), u(t_k)}_H\big]
\eea
The first term on the right-hand side of the above inequality converges by the weak convergence of $v^N$. Moreover, by the continuity of $C_F$, the second term converges to $\half \int_0^1 \qv{C_F(t) u(t), u(t)}_H dt$. Therefore, taking the supremum over all such finite-dimensional continuous $u$ and using \lemref{lem:quad-lift}, we have \bea\ad
\liminf_{N\to \infty}\CI_N^{\rm hom}(\ga_N) \geq  \int_0^1 L_{C_F}(t,v(t))dt = \CI_{C_F}(\ga),
\eea
which proves \eqref{CF-leq}.

To prove the converse inequality \eqref{CF-geq}, we suppose that $\CI_{C_F}(\ga)<\infty$ and set $v=\dot{\ga}-A \ga$. By \lemref{lem:MP-energy}, there exists a measurable $g\in L^2([0,1];H)$ such that $v(t)=C_F(t)^{1/2}g(t)$ and $\CI_{C_F}(\ga)=\half \int_0^1 |g(t)|_H^2 dt$. For $t_k\leq t <t_{k+1}$, define  $
g_k^N :=\Dl^{-1} \int_{t_k}^{t_{k+1}} g(s)ds$ and $w_k^N= e^{\Dl A} C_F(t_k)^{1/2} g^N_k$. Starting from $\ga_N(0)=0$, we define $\ga_N(t_{k+1})= e^{\Dl A} \ga_N(t_k) + \Dl w_k^N$ and set $\ga_N(t):=e^{(t-t_k)A} \ga_N(t_k)$ for $t_k \leq t<t_{k+1}$. Then $\ga_N \in \mathscr{S}_N$. Because
$
e^{\Dl A} C_F(t_k)e^{\Dl A^*}= (e^{\Dl A} C_F(t_k)^{1/2})(e^{\Dl A} C_F(t_k)^{1/2})^*,
$
the quadratic variational formula gives
$
L(e^{\Dl A} C_F(t_k) e^{\Dl A^*}, w_k^N)\leq \half |g_k^N|_H^2$.
By Jensen's inequality, it yields
\beq{temporal-recovery-cost}
\CI_N^{\rm hom}(\ga_N)\leq \half\sum_{k=0}^{N-1} \Dl |g_k^N|_H^2 \leq \half \int_0^1 |g(t)|_H^2 dt = \CI_{C_F}(\ga).
\eeq
It remains to verify the uniform convergence. By the continuity of $C_F$, we have
\bea\ad
\sup_{|t-s|\leq \Dl}\|C_F(t)^{1/2}-C_F(s)^{1/2}\|_{\CL(H)} \to 0,
\eea
as $N\to \infty$. Moreover, $\sum_{k=0}^{N-1} g_k^N \indi_{[t_k,t_{k+1})} \to g$ in $L^2([0,1];H)$. Consequently, the following convergence holds in $L^2([0,1];H)$:
\bea\ad
\sum_{k=0}^{N-1} C_F(t_k)^{1/2}g_k^N \indi_{[t_k,t_{k+1})} \to C_F(\cdot)^{1/2}g(\cdot) =v(\cdot)
\eea
Iteration of the defining recursion yields, for $1\leq j\leq N$,
\bea\ad
\ga_N(t_j)=\sum_{k=0}^{j-1} \Dl e^{(t_j-t_k)A} C_F(t_k)^{1/2} g_k^N, \quad \ga(t_j)=\sum_{k=0}^{j-1} \int_{t_k}^{t_{k+1}} e^{(t_j-s)A} v(s)ds
\eea
We choose $g^\ell\in C([0,1];H_{m_\ell})$  such that
$\|g^\ell-g\|_{L^2([0,1];H)}\to0$.  Let
$v^\ell(t)=C_F(t)^{1/2}g^\ell(t)$, and let
$\gamma_N^\ell$ and $\gamma^\ell$ be obtained from $g^\ell$ by the
same formulas as $\gamma_N$ and $\gamma$.  If
$g_k^{\ell,N}=\Dl^{-1}\int_{t_k}^{t_{k+1}}g^\ell(s)ds$, then
\begin{align*}
 &\max_{k<N}\sup_{t_k\leq s<t_{k+1}}
 |C_F(t_k)^{1/2}g_k^{\ell,N}-v^\ell(s)|_H\\
 & \leq K_C^{1/2}
 \sup_{|t-s|\leq\Dl}|g^\ell(t)-g^\ell(s)|_H
 +\|g^\ell\|_\infty
 \sup_{|t-s|\leq\Dl}
 \|C_F(t)^{1/2}-C_F(s)^{1/2}\|_{\CL(H)}\to 0.
\end{align*}
Since $v^\ell([0,1])$ is compact in $H$,
$
\sup_{0\leq r\leq\Dl}\sup_{0\leq t\leq1}
 |(e^{rA}-I)v^\ell(t)|_H\to 0.
$
The preceding three displays and the contraction property imply
\begin{align*}
 \max_{j\leq N}|\gamma_N^\ell(t_j)-\gamma^\ell(t_j)|_H
 &\leq\max_{k<N}\sup_{t_k\leq s<t_{k+1}}
 |C_F(t_k)^{1/2}g_k^{\ell,N}-v^\ell(s)|_H\\
 &\, +\sup_{0\leq r\leq\Dl}\sup_{0\leq t\leq1}
 |(e^{rA}-I)v^\ell(t)|_H\to 0.
\end{align*}
For every $N$ and $\ell$,
 $\max_{j\leq N}|\gamma_N(t_j)-\gamma_N^\ell(t_j)|_H
 \leq K_C^{1/2}\|g-g^\ell\|_{L^1([0,1];H)}$ and
 $\|\gamma-\gamma^\ell\|_\infty
 \leq K_C^{1/2}\|g-g^\ell\|_{L^1([0,1];H)}$.
Consequently,
\begin{equation}\label{temporal-recovery}
 \limsup_{N\to\infty}\max_{j\leq N}
 |\gamma_N(t_j)-\gamma(t_j)|_H
 \leq2K_C^{1/2}\|g-g^\ell\|_{L^2([0,1];H)}\to 0.
\end{equation}
For $t_k\leq t<t_{k+1}$,
\begin{align*}
 |\gamma_N(t)-\gamma(t)|_H
 &\leq|\gamma_N(t_k)-\gamma(t_k)|_H
   +|\gamma(t)-\gamma(t_k)|_H +|(e^{(t-t_k)A}-I)\gamma(t_k)|_H.
\end{align*}
Letting $N\to\infty$ in the preceding estimate and using the
compactness of $\gamma([0,1])$ yield $
 \|\gamma_N-\gamma\|_\infty\to 0$.
The jumps satisfy $
 \max_{k<N}|\gamma_N(t_{k+1})-\gamma_N(t_{k+1}-)|_H
 =\max_{k<N}\Dl|w_k^N|_H
 \leq K_C^{1/2}\Dl^{1/2}\|g\|_{L^2([0,1];H)}\to 0$.
Together with \eqref{temporal-recovery-cost}, the uniform convergence proves
\eqref{CF-geq} when
$\CI_{C_F}(\gamma)<\infty$.  If
$\CI_{C_F}(\gamma)=+\infty$, the choice $\gamma_N=\gamma$
satisfies
$
 \|\gamma_N-\gamma\|_\infty=0$ and
 \bea
\ad \limsup_{N\to\infty}\CI_N^{\rm hom}(\gamma_N)\leq+\infty
 =\CI_{C_F}(\gamma).\eea

It remains to prove the relative compactness.  Let $N_j\to\infty$ and
$\CI_{N_j}^{\rm hom}(\gamma_{N_j})\leq \kappa$.  With $N=N_j$, it follows that $
 \int_0^1|v^N(t)|_H^2dt
 =\sum_{k=0}^{N-1}\Dl|v_k^N|_H^2
 \leq2K_C\kappa$.
Then
 $\gamma_N(t_j)=\sum_{k=0}^{j-1}\Dl
 e^{(j-1-k)\Dl A}v_k^N$.
Hence
$
 \max_{j\leq N}|\gamma_N(t_j)|_H
 \leq\sum_{k=0}^{N-1}\Dl|v_k^N|_H
 \leq(2K_C\kappa)^{1/2}.
$
Since $\Pi_mA=A\Pi_m$, the Cauchy--Schwarz inequality leads to
\begin{align*}
 |(I-\Pi_m)\gamma_N(t_j)|_H^2
 &\leq\Big\{\sum_{k=0}^{j-1}\Dl
 e^{-2\alpha_{m+1}(j-1-k)\Dl}\Big\}
 \sum_{k=0}^{j-1}\Dl|v_k^N|_H^2\\
 &\leq2K_C\kappa \frac{\Dl}
 {1-e^{-2\alpha_{m+1}\Dl}} \leq2K_C\kappa \Big(\Dl+\frac1{2\alpha_{m+1}}\Big).
\end{align*}
For $t_k\leq t<t_{k+1}$,
$(I-\Pi_m)\gamma_N(t)=e^{(t-t_k)A}(I-\Pi_m)\gamma_N(t_k)$,
thus
\bea\ad
 \lim_{m\to\infty}\limsup_{j\to\infty}
 \sup_{t\leq1}|(I-\Pi_m)\gamma_{N_j}(t)|_H=0.
 \eea
For grid indices $0\leq r<q\leq N$, iteration from $t_r$ to $t_q$ yields
\bea\ad
 \Pi_m\gamma_N(t_q)-\Pi_m\gamma_N(t_r)
 =(e^{(q-r)\Dl A_m}-I)\Pi_m\gamma_N(t_r) +\sum_{k=r}^{q-1}\Dl
 e^{(q-1-k)\Dl A_m}\Pi_m v_k^N.
\eea
Since $\|e^{tA_m}-I\|_{\CL(H_m)}\leq\alpha_m t$, the preceding
energy and uniform bounds imply
\begin{align*}
 |\Pi_m\gamma_N(t_q)-\Pi_m\gamma_N(t_r)|_H
 &\leq(2K_C\kappa)^{1/2}
 \{\alpha_m(q-r)\Dl+((q-r)\Dl)^{1/2}\}.
\end{align*}
The same bounds on each interval yield
\begin{align*}
 \sup_{|t-s|\leq\delta}
 |\Pi_m\gamma_N(t)-\Pi_m\gamma_N(s)|_H
 &\leq(2K_C\kappa)^{1/2}
 \{\alpha_m(\delta+4\Dl)
       +(\delta+2\Dl)^{1/2}\},\\
 \max_{k<N}\Dl|v_k^N|_H
 &\leq(2K_C\kappa \Dl)^{1/2}\to 0.
\end{align*}
Let $\overline\gamma_{N,m}$ be the continuous piecewise linear
interpolation of $\Pi_m\gamma_N(t_k)$.  Hence $
 \|\Pi_m\gamma_N-\overline\gamma_{N,m}\|_\infty
 \leq(2K_C\kappa)^{1/2}(2\alpha_m\Dl+\Dl^{1/2})
 \to 0$.
Its modulus of continuity satisfies
\begin{equation*}
 \lim_{\delta\downarrow0}\limsup_{j\to\infty}
 \sup_{|t-s|\leq\delta}
 |\overline\gamma_{N_j,m}(t)-
   \overline\gamma_{N_j,m}(s)|_H=0.
\end{equation*}
Since $H_m$ is finite dimensional, the uniform bound and the last
modulus estimate give a
subsequence on which $\overline\gamma_{N_j,m}$ converges uniformly.
A diagonal subsequence works for every $m$.  For two members of this
subsequence, the spatial-tail estimate gives
\begin{align*}
 \|\gamma_{N_p}-\gamma_{N_q}\|_\infty
 &\leq\|(I-\Pi_m)\gamma_{N_p}\|_\infty
      +\|(I-\Pi_m)\gamma_{N_q}\|_\infty +\|\Pi_m\gamma_{N_p}-\Pi_m\gamma_{N_q}\|_\infty
 \to 0
\end{align*}
by first letting $p,q\to\infty$ and then $m\to\infty$.
Thus $(\gamma_{N_j})$ is relatively compact in the uniform norm.
The vanishing jumps and the piecewise linear interpolations show that
every limit is continuous.  Since $\gamma_{N_j}(0)=0$, $
 \gamma\in C^0([0,1];H), \gamma(0)=0$, thus $\gamma\in C^0([0,1];H)$. \qed

\subsection{Proof of \lemref{new:averaging}}\label{app:new-averaging}
For $i,j\leq n$, define
\bea
c_{ij}(t,z):=\langle \mathsf{Q}^*G(\bar X(t),z)^*e_i,\mathsf{Q}^*G(\bar X(t),z)^*e_j\rangle_H
\eea
and
$\bar c_{ij}(t):=\int_\SE c_{ij}(t,z)\pi(dz)$.
\remref{rem:G-bounds} gives
\[|c_{ij}(t,z)|\leq |\mathsf{Q}^*G(\bar X(t),z)^* e_i|_H |\mathsf{Q}^*G(\bar X(t),z)^*e_j|_H \leq e^{\al_i+\al_j}(K_G+K_G'\|\bar X\|_\infty)^2.\]  Since the range of $\bar X$ is compact in
$H$, $\bar X$ is uniformly continuous, and (H2) gives
\bea\ad
\lim_{|t-s|\to0}\sup_{z\in\SE}
|c_{ij}(t,z)-c_{ij}(s,z)|=0.
\eea
Since $|\bar c_{ij}(t)-\bar c_{ij}(s)|\leq \sup_z |c_{ij}(t,z)-c_{ij}(s,z)|$, the function $b_{ij}:=c_{ij}-\bar c_{ij}$ is bounded by $2e^{\al_i+\al_j}(K_G+K_G'\|\bar X\|_\infty)^2$, centered under $\pi$ for every $t$, and satisfies $\sup_{z}|b_{ij}(t,z)-b_{ij}(s,z)| \leq 2 \sup_{z}|c_{ij}(t,z)-c_{ij}(s,z)| \to 0$. Thus, $b_{ij}$ satisfies the hypotheses of
Lemma~\ref{lem:averaging}.  Define
\beqa{Err}\ad
\text{Err}_n^\e(t)=\int_0^t\,
[\Pi_n G(\bar X(s),\zeta_{s/\e})\mathsf{Q}\mathsf{Q}^*G(\bar X(s),\zeta_{s/\e})^*\Pi_n-C_{G,n}(s)]\,ds.
\eeqa
In the basis $(e_1,\ldots,e_n)$, by the definition of $C_G$ in \eqref{CG},
\bea\ad
\qv{\text{Err}_n^\e(t)e_j,e_i}_H
=\int_0^t [c_{ij}(s,\zeta_{s/\e})-\bar c_{ij}(s)]\,ds=\int_0^t b_{ij}(s,\zeta_{s/\e})ds
.
\eea
Moreover, since the operator norm on $H_n$ is dominated by the Frobenius norm,
\bea\ad
\|\text{Err}_n^\e(t)\|_{\CL(H_n)} \leq \Big(\sum_{i,j=1}^n |\qv{\text{Err}_n^\e(t)e_j, e_i}_H|^2\Big)^{1/2} \leq n \max_{i,j\leq n} |\qv{\text{Err}_n^\e(t) e_j, e_i}_H|
\eea
Consequently,
\bea\ad
 \PP_\mu\Big(\sup_{t\leq1}\|\text{Err}_n^\e(t)\|_{\CL(H_n)}>r\Big)\leq\sum_{i,j=1}^n
 \PP_\mu\Big(\sup_{t\leq1}
 \Big|\int_0^t b_{ij}(s,\zeta_{s/\e})ds\Big|
 >\frac{r}{n}\Big).
\eea
Lemma~\ref{lem:averaging} with $r/n$ in place of $r$ then completes the proof. \qed

\subsection{Proof of \lemref{new:finite}}\label{app:new-finite}
Conditioned on $\zeta$, the process $Z_{G,n}^\e$ is a centered Gaussian
martingale in the $n$-dimensional space $H_n$, and its predictable quadratic variation satisfies
\bea\ad
 h^2(\e)\langle Z_{G,n}^\e\rangle_t =\int_0^t\Pi_nG(\bar X(s),\zeta_{s/\e})\mathsf{Q}\mathsf{Q}^*G(\bar X(s),\zeta_{s/\e})^*\Pi_n \,ds
\eea
and $\|\Pi_n G(\bar X(s),z)\mathsf{Q}\mathsf{Q}^*G(\bar X(s),z)^* \Pi_n\|_{\CL(H_n)} \leq K_n:=ne^{2\al_n}(K_G\!+\!K_G'\|\bar X\|_\infty)^2$ for all $s,z$.
Lemma~\ref{new:averaging} gives that
$h^2(\e)\langle Z_{G,n}^\e\rangle_t\to\int_0^tC_{G,n}(s)ds$ at speed
$h^{-2}(\e)$, uniformly in $t$.

For later
use, the same bound on the predictable quadratic variation gives, for fixed
$n$, $a\in(0,1]$, and $r,R>0$,
\begin{align}
 \PP\Big(\sup_{|t-s|\leq a}
 |Z_{G,n}^\e(t)-Z_{G,n}^\e(s)|_H>r\,\Big|\,\zeta\Big)
 &\leq4n(2+a^{-1})
 \exp\Big\{-\frac{h^2(\e)r^2}{32nK_na}\Big\},
 \label{new:path-osc}\\
 \PP\big(\|Z_{G,n}^\e\|_\infty>R\,\big|\,\zeta\big)
 &\leq2n\exp\Big\{-\frac{h^2(\e)R^2}{2nK_n}\Big\}.
 \label{new:compact-containment}
\end{align}
The convergence of the predictable quadratic variation from Lemma~\ref{new:averaging} is exactly the hypothesis of
\cite[Corollary~6.7]{Puh94} for the continuous martingale $Z_{G,n}^\e$ with $Z_{G,n}^\e(0)=0$, as in Step~3 of the proof of \propref{prop:sd-driving}.
It gives (a), that is, the path rate
$J_{G,n}(z)=\int_0^1L_{C_{G,n}}(s,\dot z(s))ds$.

For (b), stochastic Fubini's theorem (the integrand is bounded since $A_n$ is bounded and $H_n$ is finite dimensional) establishes
\bea\disp
 \int_0^t A_n \Pi_n M^\e(s)ds
\ad =\frac1{h(\e)}\int_0^t\int_r^tA_ne^{(s-r)A}\Pi_n G(\bar X(r),\zeta_{r/\e})\mathsf{Q}\,ds\,dW_0(r)
\\
\aad =\Pi_n M^\e(t)-Z_{G,n}^\e(t),
\eea
because $\int_r^tA_ne^{(s-r)A}ds=e^{(t-r)A}-I$ on $H_n$, so $\Pi_n M^\e(t)=Z_{G,n}^\e(t)+\int_0^tA_n \Pi_nM^\e(s)ds$.  Hence $\Pi_nM^\e=\Phi_n(Z_{G,n}^\e)$ with $\Phi_n$ in \eqref{def-Phi-n}. Since $\Phi_n$ is continuous, the contraction principle applied to the LDP in (a) proves (b). The contracted rate is \eqref{new:IN} because \mbox{$d(\Phi_n^{-1}x)/dt=\dot x-A_nx$}. \qed

\subsection{Proof of \lemref{lem:conditional-M}}\label{app:lem-conditional-M}
Recall that $\CI_w=\CI_{C_G}$ by \eqref{new:IW}. For a deterministic c\`adl\`ag path $\omega$ in $\SE$, let $\Pcond_\e(\omega,\cdot)$ be the law in $C^0([0,1];H)$ of the process
$h^{-1}(\e)\int_0^\cdot e^{(\cdot-s)A}G(\bar X(s),\omega(s/\e))\mathsf{Q}\,dW_0(s)$. Since $W_0$ is independent of $\zeta$, $\Pcond_\e(\zeta,\cdot)$ is a regular conditional law of $M^\e$ given $\zeta$: $\PP(M^\e\in V\mid\zeta)=\Pcond_\e(\zeta,V)$ almost surely, for every Borel set $V\subset C^0([0,1];H)$.

\emph{Step 1: the events $\Omega_\e$.} For a path $\omega$ as above and $\bar n\in\NN_+$, set
\[
 \text{Err}_{\bar n}^\e(\omega)(t):=\int_0^t\big[\Pi_{\bar n}G(\bar X(s),\omega(s/\e))\mathsf{Q}\mathsf{Q}^*G(\bar X(s),\omega(s/\e))^*\Pi_{\bar n}-C_{G,\bar n}(s)\big]\,ds ,
\]
so that $\text{Err}_{\bar n}^\e(\zeta)$ is the error in \eqref{Err}. The spectral projections are nested: $\Pi_n\Pi_{\bar n}=\Pi_{\bar n}\Pi_n=\Pi_n$ for $n\leq \bar n$. Since $C_{G,\bar n}=\Pi_{\bar n}C_G\Pi_{\bar n}|_{H_{\bar n}}$, we have $\text{Err}_n^\e(\omega)(t)=\Pi_n\text{Err}_{\bar n}^\e(\omega)(t)|_{H_n}$, so
\beq{Err-nested}
 \|\text{Err}_n^\e(\omega)\|_{C([0,1];\mathcal L(H_n))}
 \leq\|\text{Err}_{\bar n}^\e(\omega)\|_{C([0,1];\mathcal L(H_{\bar n}))},
 \qquad n\leq \bar n.
\eeq
The dimension $\bar n$ and the accuracy of the predictable quadratic variation are independent parameters. Fix sequences $\bar n_j\uparrow\infty$ in $\NN_+$ and $r_j\downarrow0$, and let $\mathscr{A}_{\e,j}$ be the set of paths $\omega$ with
$\|\text{Err}_{\bar n_j}^\e(\omega)\|_{C([0,1];\mathcal L(H_{\bar n_j}))}\leq r_j$.
For each $j$, \eqref{new:qv} with $n=\bar n_j$ and $r=r_j$ gives $\bar\e_j>0$ such that
$\PP(\zeta\notin\mathscr{A}_{\e,j})\leq e^{-j h^2(\e)}$ for $0<\e\leq\bar\e_j$; we may take $\bar\e_j\downarrow0$.
For $0<\e\leq\bar\e_1$ let $j_\e:=\max\{j:\e\leq\bar\e_j\}$, so that $j_\e\to\infty$ as $\e\to0$, and set
$\mathscr{A}_\e:=\mathscr{A}_{\e,j_\e}$ and $\Omega_\e:=\{\zeta\in\mathscr{A}_\e\}$.  Then $\Omega_\e$ is $\sigma(\zeta)$-measurable and
$h^{-2}(\e)\log\PP(\Omega_\e^c)\leq-j_\e \to -\infty$. Moreover, by \eqref{Err-nested}, for every fixed $n$, every $\e$ with $\bar n_{j_\e}\geq n$, and every $\omega\in\mathscr{A}_\e$,
\beq{Err-small}
 \|\text{Err}_n^\e(\omega)\|_{C([0,1];\mathcal L(H_n))}\leq r_{j_\e}.
\eeq

\emph{Step 2: a sequential conditional LDP.} Let $\e_k\downarrow 0$ and
$\omega_k\in \mathscr{A}_{\e_k}$. We claim that the laws
$\Pcond_{\e_k}(\omega_k,\cdot)$ satisfy an LDP in $C^0([0,1];H)$ with speed $h^{-2}(\e_k)$ and good rate $\CI_w$.
Indeed, fix $n$. For large $k$ we have $\bar n_{j_{\e_k}}\geq n$, hence, by \eqref{Err-small},
\bea\ad
 \sup_{t\leq1}\big\|\text{Err}_n^{\e_k}(\omega_k)(t)\big\|_{\mathcal L(H_n)}\leq r_{j_{\e_k}}\to0 .
\eea
Under $\Pcond_{\e_k}(\omega_k,\cdot)$, the process $\Phi_n^{-1}(\Pi_nm)$ is a continuous centered Gaussian martingale. Hence, for every $H_n$-valued step function $u$,
\bea
 \ad h^{-2}(\e_k)\log
 \int_{C^0([0,1];H)}\exp\big\{h^2(\e_k)
 \int_0^1\big\langle u(s),
 d\, \big[\Phi_n^{-1}(\Pi_nm)\big](s)\big\rangle_H\big\}
 \Pcond_{\e_k}(\omega_k,dm)\\
 \ad =\frac12\int_0^1
 \langle \Pi_n G(\bar X(s),\omega_k(s/\e_k))\mathsf{Q}\mathsf{Q}^*G(\bar X(s),\omega_k(s/\e_k))^*
 \Pi_n u(s),u(s)\rangle_Hds \\
 \aad \to
 \frac12\int_0^1\langle C_{G,n}(s)u(s),u(s)\rangle_Hds.
\eea
The conditional versions of \eqref{new:path-osc} and
\eqref{new:compact-containment} hold uniformly in $\omega_k$, and
\eqref{new:cond-tail} gives
\bea\ad
 \limsup_{k\to\infty}h^{-2}(\e_k)\log
 \Pcond_{\e_k}\!\big(\omega_k,
 \{m:\|(I-\Pi_n)m\|_\infty>r\}\big)
 \leq-\kappa_\delta r^2\alpha_{n+1}^{\delta}\xrightarrow[n\to\infty]{}-\infty.
\eea
Thus the finite-dimensional G\"artner--Ellis theorem and exponential tightness give the LDP of $\Phi_n^{-1}(\Pi_nm)$ with rate $J_{G,n}$, and the continuity of $\Phi_n$ and \corref{thm:abstract-lift} give the claimed LDP.

\emph{Step 3: small balls.} Fix $\phi$ with $\CI_w(\phi)<\infty$ and $\eta>0$. Suppose that \eqref{cond-Sch-low} fails for some $a>0$. Then there are $\e_k\downarrow0$ such that the event $\Omega_{\e_k}\cap\{\Pcond_{\e_k}(\zeta,B(\phi,a))<\exp\{-h^2(\e_k)(\CI_w(\phi)+\eta)\}\}$ has positive probability. Since $\zeta\in\mathscr{A}_{\e_k}$ on $\Omega_{\e_k}$, we find $\omega_k\in\mathscr{A}_{\e_k}$ with
$
 h^{-2}(\e_k)\log\Pcond_{\e_k}\big(\omega_k,B(\phi,a)\big)<-\CI_w(\phi)-\eta .
$
Since $B(\phi,a)$ is open and contains $\phi$, the lower bound of Step~2 gives
\[
\liminf_k h^{-2}(\e_k)\log\Pcond_{\e_k}(\omega_k,B(\phi,a))
\ge-\CI_w(\phi),
\]
a contradiction.

For the upper bound, the lower semicontinuity of $\CI_w$ gives $a_0>0$ such that $\CI_w(x)>\CI_w(\phi)-\eta/2$ whenever $\|x-\phi\|_{\infty}\leq a_0$. Let $a\in(0,a_0]$, and suppose that \eqref{cond-Sch-up} fails. As above, we find $\e_k\downarrow0$ and $\omega_k\in\mathscr{A}_{\e_k}$ with
$
 h^{-2}(\e_k)\log\Pcond_{\e_k}\big(\omega_k,\overline B(\phi,a)\big)>-\CI_w(\phi)+\eta .
$
The upper bound of Step 2 for the closed ball gives
\[
\limsup_k h^{-2}(\e_k)\log\Pcond_{\e_k}(\omega_k,\overline B(\phi,a))
\le-\inf_{\overline B(\phi,a)}\CI_w
\le-\CI_w(\phi)+\eta/2,
\]
which is
again a contradiction. \qed

\subsection{Proof of Proposition~\ref{prop:eta-equivalence}}
\label{app:exp-equivalence}
Denote by $\wdt X^\e(t)$ the solution of
\beq{til-X}\left\{ \barray
d\wdt X^\e(t) = [A \wdt X^\e(t)+F(\wdt X^\e(t),\zeta_{t/\e})]dt + \sqrt{\e} G(\bar X(t),\zeta_{t/\e}) dW^{\mathsf{Q}}(t), \\
\wdt X^\e(0)=x_0.
\earray\right.\eeq
Define the following auxiliary process
\beq{til-eta}
\wdt \eta^\e(t) =\frac{\wdt X^\e(t)-\bar X(t)}{\sqe h(\e)},\ t\in [0,1].
\eeq
We prove the exponential equivalence through the intermediate
process $\wdt\eta^\e$.
\smallskip

\noindent\textbf{Step 1}.  For every $r>0$, we first prove that
\beq{eta-wdt-eta}
\lim_{\e\to 0} h^{-2}(\e) \log \PP(\|\eta^\e-\wdt \eta^\e\|_\infty > r )=-\infty.
\eeq
Define
\beqa{M1-e}
M_1^\e(t):= h^{-1}(\e) \int_0^t e^{(t-s)A}\big[G(X^\e(s),\zeta_{s/\e})- G(\bar X(s),\zeta_{s/\e})\big] dW^{\mathsf{Q}}(s).
\eeqa
The uniform Lipschitz continuity of $F(\cdot,z)$ implies
\[
\begin{aligned}
|\eta^\e(t)-\wdt \eta^\e(t)|_H
& \leq\Big| \frac{1}{\sqe h(\e)}\int_0^t e^{(t-s)A}
[F(X^\e(s),\zeta_{s/\e})-F(\wdt X^\e(s),\zeta_{s/\e})] ds\Big|_H + |M_1^\e(t)|_H \\
&\leq \int_0^t K |\eta^\e(s)-\wdt \eta^\e(s)|_H ds + |M_1^\e(t)|_H.
\end{aligned}
\]
Gronwall's inequality gives
$
|\eta^\e(t)-\wdt \eta^\e(t)|_H \leq e^K \|M_1^\e\|_\infty$. Thus \eqref{eta-wdt-eta} follows from
\beq{neg-M1e}
\lim_{\e\to 0} h^{-2}(\e) \log\PP(\|M_1^\e\|_\infty > r)=-\infty, \quad \forall\;  r>0
\eeq
The term $M_1^\e$ is a stochastic convolution. Its estimate uses
only the Lipschitz property of $G$, without differentiating $G$, and $X^\e(s)-\bar X(s)=\sqrt{\e}h(\e)\eta^\e(s)$.

Fix $R_1,R_2>0$. Then
\bea
\PP(\|M_1^\e\|_\infty > r) \ad \leq \PP(\|M_1^\e\|_\infty > r; \|X^\e\|_\infty \leq R_1, \|\eta^\e\|_\infty \leq R_2) \\
\aad\quad + \PP(\|\eta^\e\|_\infty > R_2; \|X^\e\|_\infty \leq R_1) + 2 \PP(\|X^\e\|_\infty > R_1 )\\
\aad\mathrel{=:} \mathrm{I}_1^\e + \mathrm{I}_2^\e + 2 \mathrm{I}_3^\e.
\eea
For $\mathrm{I}_1^\e$ define the stopping time
$
 \tau_{R_1,R_2}^\e:=\inf\{t\in[0,1]:
 |X^\e(t)|_H>R_1\ \text{or}\ |\eta^\e(t)|_H>R_2\}\wedge 1$,
and the stopped convolution
\bea
 M_{1,R_1,R_2}^\e(t):=h^{-1}(\e)\int_0^te^{(t-s)A}
 \indi_{\{s\leq\tau_{R_1,R_2}^\e\}}
 [G(X^\e(s),\zeta_{s/\e})-G(\bar X(s),\zeta_{s/\e})]\mathsf{Q}\,dW_0(s).
\eea
On the event defining $\mathrm{I}_1^\e$,
$M_{1,R_1,R_2}^\e(t)=M_1^\e(t)$ for all $t\leq1$.  Moreover, by \eqref{G-diff},
\bea\ad
 \Big\|h^{-1}(\e)e^{(t-s)A}\indi_{\{s\leq\tau_{R_1,R_2}^\e\}}
 [G(X^\e(s),\zeta_{s/\e})-G(\bar X(s),\zeta_{s/\e})]\mathsf{Q}
 \Big\|_{\CL_2(H)}^2
\leq K\e R_2^2(t-s)^{-\vartheta}.
\eea
Lemma~\ref{lem:conv-exp} confirms
\begin{equation}\label{PM}
 h^{-2}(\e)\log \mathrm{I}_1^\e
 \leq h^{-2}(\e)\log K_1-\kappa_1 r^2/(\e h^2(\e)R_2^2) \to -\infty
\end{equation}
for fixed $R_1,R_2,r>0$, where $\kappa_1:=(K_2K)^{-1}$.

For $\mathrm{I}_2^\e$, the uniform Lipschitz bound on
$F(\cdot,z)$
implies that
on bounded sets of $X^\e$,
\bea
|\eta^\e(t)|_H \ad = \Big|\feh \int_0^t e^{(t-s)A}[ F(X^\e(s), \zeta_{s/\e})-\bar F(\bar X(s))]ds\\
\aad \quad + h^{-1}(\e) \int_0^t e^{(t-s)A} G(X^\e(s), \zeta_{s/\e})dW^{\mathsf{Q}}(s) \Big|_H \\
\aad = \Big|\feh \int_0^t e^{(t-s)A}[ F(X^\e(s), \zeta_{s/\e})-F(\bar X(s), \zeta_{s/\e})] ds \\
\aad \quad + \la^\e(t) + h^{-1}(\e) \int_0^t e^{(t-s)A} G(X^\e(s), \zeta_{s/\e})dW^{\mathsf{Q}}(s) \Big|_H \\
\aad \leq K \int_0^t |\eta^\e(s)|_H ds + \Big|\la^\e(t)+h^{-1}(\e)\int_0^t e^{(t-s)A} G(X^\e(s), \zeta_{s/\e})dW^{\mathsf{Q}}(s)\Big|_H.
\eea
Gronwall's inequality then gives
\bea\ad
\|\eta^\e\|_\infty \leq K \Big(\sup_{t\in [0,1]}|\la^\e(t)|_H + \sup_{t\in [0,1]} \Big|h^{-1}(\e) \int_0^t e^{(t-s)A} G(X^\e(s), \zeta_{s/\e})dW^{\mathsf{Q}}(s)\Big|_H \Big).
\eea
Consequently,
\bea
\mathrm{I}_2^\e \ad \leq \PP \Big(\sup_{t\in [0,1]}|\la^\e(t)|_H > \frac{R_2}{2K}; \sup_{t\in [0,1]} |X^\e(t)|_H \leq R_1 \Big)\\
\aad + \PP \Big(\sup_{t\in [0,1]}\Big|\frac1{h(\e)} \int_0^t e^{(t-s)A} G(X^\e(s),\zeta_{s/\e})dW^{\mathsf{Q}}(s) \Big|_H > \frac{R_2}{2K}; \sup_{t\in [0,1]} |X^\e(t)|_H \leq R_1\Big).
\eea
Thus the required estimate
\beq{P-eta}
\lim_{R_2 \to \infty} \limsup_{\e \to 0} h^{-2}(\e)  \log \mathrm{I}_2^\e = -\infty
\eeq
follows if we can establish
\beq{P-eta1}
\lim_{R_2 \to \infty} \limsup_{\e \to 0}  h^{-2}(\e)  \log\PP\Big(\sup_{t\in [0,1]}|\la^\e(t)|_H > \frac{R_2}{2 K}; \sup_{t\in [0,1]} |X^\e(t)|_H \leq R_1 \Big) = -\infty.
\eeq
\beq{P-eta2}\barray\disp
\lim_{R_2 \to \infty} \limsup_{\e \to 0}  h^{-2}(\e)  \log\PP\Big(\ad \sup_{t\in [0,1]} \Big|h^{-1}(\e) \int_0^t e^{(t-s)A} G(X^\e(s), \zeta_{s/\e})dW^{\mathsf{Q}}(s)\Big|_H > \frac{R_2}{2 K}; \\
\aad \sup_{t\in [0,1]} |X^\e(t)|_H \leq R_1\Big)=-\infty.
\earray\eeq
The two limits can be checked quantitatively.  Since $\CI_\zeta$ is
good, $ \kappa_\zeta(R):=\inf_{\|x\|_\infty\geq R}\CI_\zeta(x)
\to\infty$ as $R\to\infty$,
and the LDP upper bound gives
\[
 \limsup_{\e\to0}h^{-2}(\e)\log
 \PP(\|\lambda^\e\|_\infty>R_2/(2K))
 \leq-\kappa_\zeta(R_2/(2K)) \to-\infty.
\]
This proves \eqref{P-eta1}.  For the stochastic convolution in \eqref{P-eta2}, we stop at $\tau_{R_1,\infty}^\e=\inf\{t\in[0,1]:|X^\e(t)|_H>R_1\}\wedge1$, that is, $\tau_{R_1,R_2}^\e$ with $R_2=\infty$, which does not change the event in \eqref{P-eta2}. Then Assumption~\ref{ass:fg} (H2) gives
\[
\|h^{-1}(\e)e^{(t-s)A}\indi_{\{s\leq\tau_{R_1,\infty}^\e\}}G(X^\e(s),\zeta_{s/\e})\mathsf{Q}\|_{\CL_2(H)}^2
\le\frac{K(t-s)^{-\vartheta}}{h^2(\e)}.
\]
Lemma~\ref{lem:conv-exp} therefore yields
\[
 \PP\Big(\big\|h^{-1}(\e)\int_0^\cdot e^{(\cdot-s)A}
 \indi_{\{s\leq\tau_{R_1,\infty}^\e\}}G(X^\e(s),\zeta_{s/\e})\mathsf{Q}dW_0(s)\big\|_\infty
 >\frac{R_2}{2K}\Big)
 \leq K_1e^{-\kappa_2 R_2^2h^2(\e)}.
\]
Here $\kappa_2>0$ depends only on $K$ and $K_2$. Consequently the logarithmic upper limit in \eqref{P-eta2} is at most
$-\kappa_2 R_2^2$, which tends to $-\infty$.

As for $\mathrm{I}_3^\e$, Theorem~\ref{thm:XbarX} (3) supplies the
needed estimate directly. Choose $R_1$ large enough that
$\PP(\|X^\e\|_\infty>R_1)\leq e^{-\kappa_3/\e}$ for some $\kappa_3>0$.
Then
\begin{equation}\label{neg-PX}
 \limsup_{\e\to0}h^{-2}(\e)\log \mathrm{I}_3^\e
 \leq\limsup_{\e\to0}-\kappa_3/(\e h^2(\e))=-\infty.
\end{equation}
Letting $\e\to0$ and then $R_2\to\infty$ in \eqref{PM}, \eqref{P-eta}, \eqref{neg-PX}
gives \eqref{neg-M1e},
hence \eqref{eta-wdt-eta}.
\smallskip

\noindent\textbf{Step 2}.  For every $r>0$, we prove
\beq{neg-wdt-eta-hat}
\lim_{\e\to 0}  h^{-2}(\e)  \log\PP\Big(\sup_{t\in [0,1]}|\wdt \eta^\e(t)-\wdh \eta^\e(t)|_H > r \Big) = -\infty.
\eeq
Equation \eqref{til-eta} gives
\beq{eta-e}\barray
\wdt \eta^\e(t)
\ad =\la^\e(t)+h^{-1}(\e)\int_0^t e^{(t-s)A} G(\bar X(s), \zeta_{s/\e})dW^{\mathsf{Q}}(s)  \\
\aad + \frac{1}{\sqe h(\e)} \int_0^t  e^{(t-s)A} \big[F(\wdt X^\e(s),\zeta_{s/\e}) - F(\bar X(s), \zeta_{s/\e}) \big] ds \\
\ad = \la^\e(t)+h^{-1}(\e) \int_0^t e^{(t-s)A} G(\bar X(s), \zeta_{s/\e}) dW^{\mathsf{Q}}(s) \\
\aad +\frac{1}{\sqe h(\e)} \int_0^t e^{(t-s)A} \big[F(\wdt X^\e(s),\zeta_{s/\e})-F(\bar X(s),\zeta_{s/\e})\\
\aad\qquad\qquad\qquad\qquad\qquad -\sqe h(\e)\nabla F(\bar X(s),\zeta_{s/\e})\wdt \eta^\e(s) \big]ds \\
\aad+\int_0^t e^{(t-s)A} \big[\nabla F(\bar X(s),\zeta_{s/\e}) -\overline{\nabla F}(\bar X(s))\big] \wdt \eta^\e(s)ds \\
\aad +\int_0^t e^{(t-s)A} \overline{\nabla F}(\bar X(s))\wdt \eta^\e(s)ds \\
\ad\mathrel{=:} \la^\e(t)+M^\e(t) +\int_0^t e^{(t-s)A} \Psi_1(\wdt X^\e(s),\bar X(s),\zeta_{s/\e})ds \\
\aad+\int_0^t e^{(t-s)A}\Psi_2(\bar X(s),\zeta_{s/\e}) \wdt \eta^\e(s) ds + \int_0^t e^{(t-s)A} \overline{\nabla F}(\bar X(s))\wdt\eta^\e(s)ds,
\earray\eeq
where $\Psi_1$ and $\Psi_2$ are defined as
\beq{def-psi-1}
\barray \Psi_1(\wdt X^\e(s),\bar X(s),\zeta_{s/\e})\ad :=\frac{1}{\sqe h(\e)}\big[F(\wdt X^\e(s),\zeta_{s/\e})
-F(\bar X(s),\zeta_{s/\e})\\
\ad\qquad\qquad\qquad  -\sqe h(\e)\nabla F(\bar X(s),\zeta_{s/\e})\wdt \eta^\e(s)\big]
\earray\eeq
and
\beq{def-psi-2}
\Psi_2(\bar X(s),\zeta_{s/\e}):=\nabla F(\bar X(s),\zeta_{s/\e})-\overline{\nabla F}(\bar X(s)).
\eeq
The bound $\|\overline{\nabla F}(\bar X(s))\|_{\CL(H)}\leq K_F'$ yields
\bea\ad
|\wdt \eta^\e(t) -\wdh \eta^\e(t)|_H \\
\ad \leq \Big| \int_0^t e^{(t-s)A} \overline{\nabla F}(\bar X(s)) (\wdt \eta^\e(s)-\wdh \eta^\e(s)) ds \Big|_H \\
\aad\quad +\Big| \int_0^t e^{(t-s)A} \Psi_1(\wdt X^\e(s),\bar X(s),\zeta_{s/\e})ds \Big|_H + \Big|\int_0^t e^{(t-s)A} \Psi_2(\bar X(s),\zeta_{s/\e}) \wdt \eta^\e(s) ds\Big|_H\\
\ad \leq \int_0^t \|e^{(t-s)A}\|_{\CL(H)}\|\overline{\nabla F}(\bar X(s))\|_{\CL(H)} |\wdt \eta^\e(s)-\wdh \eta^\e(s)|_H ds \\
\aad \quad +\Big|\int_0^t e^{(t-s)A} \Psi_1(\wdt X^\e(s),\bar X(s),\zeta_{s/\e})ds \Big|_H +\Big|\int_0^t e^{(t-s)A} \Psi_2(\bar X(s),\zeta_{s/\e})\wdt \eta^\e(s)ds\Big|_H \\
\aad \leq K \int_0^t |\wdt \eta^\e(s)- \wdh\eta^\e(s)|_H ds + \Big| \int_0^t e^{(t-s)A}  \Psi_1(\wdt X^\e(s), \bar X(s), \zeta_{s/\e})ds \Big|_H  \\
\aad \quad + \Big| \int_0^t e^{(t-s)A} \Psi_2(\bar X(s), \zeta_{s/\e}) \wdt \eta^\e(s) ds \Big|_H.
\eea
By Gronwall's inequality, it suffices
to prove
\beq{neg-psi1}
\lim_{\e\to 0} h^{-2}(\e)  \log\PP\Big(\sup_{t\in [0,1]}\Big| \int_0^t e^{(t-s)A} \Psi_1(\wdt X^\e(s),\bar X(s), \zeta_{s/\e})ds \Big|_H > r \Big)=-\infty,\quad \forall r>0,
\eeq
and
\beq{neg-psi2}
\lim_{\e\to 0} h^{-2}(\e)  \log\PP\Big(\sup_{t\in [0,1]} \Big|\int_0^t e^{(t-s)A} \Psi_2(\bar X(s),\zeta_{s/\e})\wdt \eta^\e(s) ds\Big|_H  >r \Big) =-\infty, \quad \forall \ r>0.
\eeq
First,
\begin{equation}\label{psi-1}
 \|e^{(t-s)A}\Psi_2(\bar X(s),z)\|_{\mathcal L(H)}\leq2K_F'.
\end{equation}
For the Taylor remainder, the bound on $\nabla F$ gives
\begin{equation}\label{psi-3}
 |\Psi_1(\bar X+\sqrt{\e}h(\e) v,\bar X,z)|_H\leq2K_F'|v|_H.
\end{equation}
For every compact $\CO\subset C([0,1];H)$, the sets
$\{\bar X(t):t\leq1\}$ and $\{v(t):v\in \CO,t\leq1\}$ are compact in
$H$.  For $x,v\in H$ and $a>0$, the fundamental theorem of calculus gives
\[
\frac{F(x+av,z)-F(x,z)}{a}-\nabla F(x,z)v
=\int_0^1[\nabla F(x+uav,z)-\nabla F(x,z)]v\,du.
\]
By the compact-uniform continuity in (H1),
\begin{equation}\label{psi-4}
 \sup_{v\in \CO,t\leq1,z\in\SE}
 |\Psi_1(\bar X(t)+\sqrt{\e}h(\e) v(t),\bar X(t),z)|_H \to 0.
\end{equation}

We verify the compactness needed to use this estimate for the random
path $\widetilde\eta^\e$. Its equation and Lipschitz continuity give
\beqa{wdt-eta-laM}
 \|\widetilde\eta^\e\|_\infty
 \leq e^{K_F'}\|R^\e\|_\infty\leq e^{K_F'}(\|\lambda^\e\|_\infty+\|M^\e\|_\infty).
\eeqa
 Moreover, $\widetilde\eta^\e=R^\e+\Scv b^\e$, where $\Scv$ is the convolution operator in \appref{app:quad-lift}, $b^\e(s):=(\sqe h(\e))^{-1}[F(\wdt X^\e(s),\zeta_{s/\e})-F(\bar X(s),\zeta_{s/\e})]$, and
$|b^\e(s)|_H\leq K_F'|\widetilde\eta^\e(s)|_H$.
The map $\Scv$ from $L^2([0,1];H)$ to $C^0([0,1];H)$ is compact by \appref{app:quad-lift}.
If $R^\e$ lies in a fixed compact set $\CO$, the norm in
\eqref{wdt-eta-laM} is bounded and $\Scv b^\e$ lies in the
compact closure of the image of a fixed $L^2$ ball. Thus
$\widetilde\eta^\e$ lies in another compact set, namely the sum of
these two compact sets. Exponential tightness of $R^\e$ then gives that of $\widetilde\eta^\e$.

For each $\kappa$, choose such a compact set $\wdt\CO_\kappa$ with
logarithmic exceptional probability at most $-\kappa$.
On $\{\widetilde\eta^\e\in\wdt\CO_\kappa\}$, \eqref{psi-4} bounds
$\sup_t|\int_0^te^{(t-s)A}\Psi_1(s)ds|_H$ by a deterministic
quantity tending to zero. Its probability of exceeding $r$ therefore
has logarithmic upper limit at most $-\kappa$. Let $\kappa \to\infty$ to obtain
\eqref{neg-psi1}.

For \eqref{neg-psi2}, equation \eqref{eta-e} gives
\bea
\ad \int_0^t e^{(t-s)A}\Psi_2(\bar X(s), \zeta_{s/\e})\wdt \eta^\e(s) ds \\
\ad =\int_0^t e^{(t-s)A}\Psi_2(\bar X(s),\zeta_{s/\e})(\la^\e(s)+M^\e(s))ds\\
\aad +\int_0^t e^{(t-s)A} \Psi_2(\bar X(s),\zeta_{s/\e}) \int_0^s e^{(s-u)A}\overline{\nabla F}(\bar X(u))\wdt \eta^\e(u)duds \\
\aad +\int_0^t e^{(t-s)A}\Psi_2(\bar X(s),\zeta_{s/\e})\int_0^s e^{(s-u)A}\Psi_1(\wdt X^\e(u),\bar X(u),\zeta_{u/\e})du ds\\
\aad +\int_0^t e^{(t-s)A}\Psi_2(\bar X(s),\zeta_{s/\e})\int_0^s e^{(s-u)A} \Psi_2(\bar X(u),\zeta_{u/\e}) \wdt\eta^\e(u) du ds.
\eea
The estimate \eqref{psi-1} and Gronwall's inequality yield
\bea\ad
\Big|\int_0^t e^{(t-s)A} \Psi_2(\bar X(s),\zeta_{s/\e})\wdt \eta^\e(s)ds\Big|_H \\
\aad\leq  e^{2K_F'} \sup_{t\in [0,1]}\Big(\Big|\int_0^t e^{(t-s)A}\Psi_2(\bar X(s),\zeta_{s/\e})\la^\e(s)ds \Big|_H \\
\aad\quad +\Big|\int_0^t e^{(t-s)A} \Psi_2(\bar X(s),\zeta_{s/\e}) M^\e(s)ds\Big|_H \\
\aad\quad + \Big|\int_0^t e^{(t-s)A} \Psi_2(\bar X(s),\zeta_{s/\e}) \int_0^s e^{(s-u)A}\overline{\nabla F}(\bar X(u))\wdt\eta^\e(u)du ds\Big|_H  \\
\aad\quad + \Big|\int_0^t e^{(t-s)A}\Psi_2(\bar X(s),\zeta_{s/\e})\int_0^s e^{(s-u)A} \Psi_1(\wdt X^\e(u),\bar X(u), \zeta_{u/\e})du ds\Big|_H
 \Big) \\
\ad=: e^{2K_F'} \sup_{t\in [0,1]}\Big(|\text{II}_1^\e(t)|_H + |\text{II}_2^\e(t)|_H + |\text{II}_3^\e(t)|_H +|\text{II}_4^\e(t)|_H\Big).
\eea
It remains to prove that $\text{II}_i^\e$, $1\leq i\leq4$,
are negligible at speed $h^{-2}(\e)$. For $\text{II}_4^\e$, \eqref{psi-1} gives
\bea\ad
\sup_{t\in [0,1]}|\text{II}_4^\e(t)|_H \leq 2K_F' \sup_{t\in [0,1]}\Big|\int_0^t e^{(t-s)A}\Psi_1(\wdt X^\e(s), \bar X(s), \zeta_{s/\e}) ds\Big|_H.
\eea
Thus \eqref{neg-psi1} implies the negligibility of $\text{II}_4^\e$.
\smallskip

\noindent (i) \textbf{Negligibility of $\mathrm{II}_3^\e$.}
We apply integration by parts after projection onto finite-dimensional
spaces. Define $V^\e(s)=\int_0^se^{(s-u)A}
\overline{\nabla F}(\bar X(u))\widetilde\eta^\e(u)du$.
On $\{\|\widetilde\eta^\e\|_\infty\leq R\}$,
$\|V^\e\|_\infty\leq K_F'R$ and
$\|(I-\Pi_n)V^\e\|_\infty\leq K_F'R/\alpha_{n+1}$, while the part of
$\mathrm{II}_3^\e$ outside $H_n$ is bounded by
$2(K_F')^2R/\alpha_{n+1}$.  Consequently,
\beqa{II3-projection}\ad
 \|\mathrm{II}_3^\e-
 \int_0^\cdot e^{(\cdot-s)A_n}\Pi_n
 \Psi_2(\bar X(s),\zeta_{s/\e})\Pi_nV^\e(s)ds \|_\infty
 \leq KR\alpha_{n+1}^{-1}.
\eeqa
For fixed $n$, set $v_n=\Pi_nV^\e$ and
$\Psi_{2,n}^\e(t)=\int_0^t\Pi_n\Psi_2(\bar X(s),\zeta_{s/\e})\Pi_n\,ds$.
Every matrix entry of the integrand defining $\Psi_{2,n}^\e$ is bounded,
centered under $\pi$, and uniformly continuous in time by
Assumption~\ref{ass:fg}(H1).  Lemma~\ref{lem:averaging} applied to
the finitely many entries shows that
\begin{equation}\label{II3-primitive}
 \limsup_{\e\to0}h^{-2}(\e)\log
 \PP(\|\Psi_{2,n}^\e\|_\infty>r)=-\infty,\quad\forall r>0.
\end{equation}
The $v_n$ is absolutely continuous and
$v_n'=A_nv_n+\Pi_n\overline{\nabla F}(\bar X)\widetilde\eta^\e$.
Integration by parts gives
\begin{align}
 &\int_0^te^{(t-s)A_n}\Pi_n
 \Psi_2(\bar X(s),\zeta_{s/\e})\Pi_nV^\e(s)ds\notag\\
 &\;=\Psi_{2,n}^\e(t)v_n(t)
 +\int_0^tA_ne^{(t-s)A_n}\Psi_{2,n}^\e(s)v_n(s)ds
 -\int_0^te^{(t-s)A_n}\Psi_{2,n}^\e(s)v_n'(s)ds.       \label{II3-IBP}
\end{align}
On $\{\|\widetilde\eta^\e\|_\infty\leq R\}$, the norm of the right-hand side of \eqref{II3-IBP}, uniformly in $t$, is at most $K_{n,R}\|\Psi_{2,n}^\e\|_\infty$.  Fix $R,r$ and choose $n$ so that the right-hand side of \eqref{II3-projection} is smaller than $r/2$. Then \eqref{II3-primitive} gives $\limsup_{\e\to0}h^{-2}(\e)\log\PP(\|\mathrm{II}_3^\e\|_\infty>r,\,\|\wdt\eta^\e\|_\infty\leq R)=-\infty$. Finally, letting $R\to\infty$ and using \eqref{wdt-eta-laM} and the marginal LDPs with good rate functions, we have the negligibility of $\mathrm{II}_3^\e$.
\smallskip

\noindent (ii) \textbf{Negligibility of $\text{II}_1^\e$ and $\text{II}_2^\e$.} We treat $\text{II}_2^\e$ first.
 For every $r>0$,
\bea\ad
\limsup_{\e \to 0}  h^{-2}(\e)  \log \PP\Big(\sup_{t\in [0,1]} |\text{II}_2^\e(t)|_H > r\Big)=-\infty.
\eea
By definition, $
\text{II}_2^\e(t):= \int_0^t e^{(t-s)A} \Psi_2(\bar X(s), \zeta_{s/\e}) M^\e(s)ds$. By Theorem~\ref{thm:mdp-comb}, $M^\e$ has the good
rate $\CI_w$. For $\kappa >0$, let $\CO_\kappa=\{\CI_w\leq \kappa \}$.
For any $\dl>0$, choose finitely many $\phi_1,\ldots,\phi_{N_\kappa}\in \CO_\kappa$ such that
\begin{equation}\label{modul-zi}
 \CO_\kappa\subset\bigcup_{i=1}^{N_\kappa} B(\phi_i,\dl).
\end{equation}
Every $\phi_i$ is a fixed continuous $H$-valued path. For such a
path define $\ell_i(s,z)=\Psi_2(\bar X(s),z)\phi_i(s)$.
It is bounded and centered. Each scalar coordinate is uniformly continuous in $s$: the first
term follows from the compact-uniform continuity of
$(x,v)\mapsto \nabla F(x,z)v$ in Assumption~\ref{ass:fg}(H1), applied to the
compact ranges of $\bar X$ and $\phi_i$, and the increment of $\phi_i$ is
bounded by $2K_F'|\phi_i(t)-\phi_i(s)|_H$.
Lemma~\ref{lem:averaging} applies to each coordinate primitive.
Integration by parts with $e^{-\alpha_j(t-s)}$ bounds the corresponding
convolution coordinate by twice its primitive supremum. Moreover, $
 \sup_t |(I-\Pi_n)\int_0^te^{(t-s)A}\ell_i(s,\zeta_{s/\e})ds |_H
 \leq 2K_F'\|\phi_i\|_\infty/\alpha_{n+1}$.
Choose $n$ first and use a finite union bound for its coordinates. We obtain
\begin{equation}\label{mod-z}
 \limsup_{\e\to0}h^{-2}(\e)\log\PP\Big(
 \Big\|\int_0^\cdot e^{(\cdot-s)A}\ell_i(s,\zeta_{s/\e})ds\Big\|_\infty>r\Big)
 =-\infty.
\end{equation}
Then for all positive $r$, with $\CO_\kappa^\dl:=\bigcup_{i=1}^{N_\kappa}B(\phi_i,\dl)$,
\bea\ad
\PP(\|\text{II}_2^\e\|_\infty>r) \leq \PP(\|\text{II}_2^\e\|_\infty>r; M^\e\in \CO_\kappa^\dl)+ \PP(M^\e\notin \CO_\kappa^\dl) \\
\aad \leq \sum_{i=1}^{N_\kappa}  \PP\big(\|\text{II}_2^\e\|_\infty > r; M^\e\in B(\phi_i,\dl)\big)+\PP\big(M^\e\notin \CO_\kappa^\dl\big)\\
\aad \leq \sum_{i=1}^{N_\kappa}\PP\Big(\sup_{t
\in [0,1]} \Big|\text{II}_2^\e(t)-\int_0^t e^{(t-s)A} \Psi_2(\bar X(s),\zeta_{s/\e}) \phi_i(s)ds \Big|_H> \frac{r}{2}; M^\e\in B(\phi_i,\dl)\Big)\\
\aad\quad+\sum_{i=1}^{N_\kappa}\PP\Big(\sup_{t\in [0,1]}\Big|\int_0^t e^{(t-s)A}\Psi_2(\bar X(s),\zeta_{s/\e})\phi_i(s)ds\Big|_H> \frac{r}{2}\Big)+\PP\Big(M^\e\notin \CO_\kappa^\dl\Big)
\eea
On the event $M^\e\in B(\phi_i,\dl)$,
\bea\ad
\sup_{t\in [0,1]}\Big|\text{II}_2^\e(t)-\int_0^t e^{(t-s)A}\Psi_2(\bar X(s),\zeta_{s/\e})\phi_i(s)ds\Big|_H \\
\aad\quad \leq \dl\int_0^t\|e^{(t-s)A}\|\,
 \|\Psi_2(\bar X(s),\zeta_{s/\e})\|_{\mathcal L(H)}ds
 \leq 2K_F'\dl.
\eea
Thus for $0< \dl<r/(4K_F')$,
\beq{I-21}
\sum_{i=1}^{N_\kappa}\PP\Big(\sup_{t
\in [0,1]} \Big|\text{II}_2^\e(t)-\int_0^t e^{(t-s)A}\Psi_2(\bar X(s),\zeta_{s/\e})\phi_i(s)ds \Big|_H> \frac{r}{2}; M^\e\in B(\phi_i,\dl)\Big)=0.
\eeq
For each of the finitely many $\phi_i$, equation \eqref{mod-z} gives
\begin{equation}\label{I-22}
 \limsup_{\e\to0}h^{-2}(\e)\log\PP\Big(
 \Big\|\int_0^\cdot e^{(\cdot-s)A}\Psi_2(\bar X(s),\zeta_{s/\e})\phi_i(s)ds
 \Big\|_\infty>r/2\Big)=-\infty.
\end{equation}
The LDP upper bound for $M^\e$ on the closed set $(\CO_\kappa^\dl)^c$, which does not meet $\CO_\kappa$, gives
\beq{I-23}
\limsup_{\e\to 0} h^{-2}(\e) \log\PP\Big(M^\e\notin \CO_\kappa^\dl\Big) \leq-\inf_{\phi\notin \CO_\kappa^\dl}\CI_w(\phi)\leq-\kappa .
\eeq
Combining \eqref{I-21}, \eqref{I-22} and \eqref{I-23}
gives
\bea\ad
\limsup_{\e\to 0} h^{-2}(\e) \log\PP(\|
\text{II}_2^\e\|_\infty > r)\leq -\kappa .
\eea
Taking $\kappa$ to infinity concludes the negligibility of $\text{II}_2^\e$. For $\mathrm{II}_1^\e$, let
\[
\mathrm{II}_1^\e(t)=\int_0^t e^{(t-s)A}
\Psi_2(\bar X(s),\zeta_{s/\e})\lambda^\e(s)\,ds,
\qquad \CO_{\zeta,\kappa}:=\{\CI_\zeta\le\kappa\}.
\]
Choose $y_1,\ldots,y_{N_\kappa}\in \CO_{\zeta,\kappa}$ so that
$\CO_{\zeta,\kappa}\subset\bigcup_{j=1}^{N_\kappa}B(y_j,\dl)$, where
$0<\dl<r/(4K_F')$, and define
$\mathrm{II}_{1,j}^\e(t)=\int_0^te^{(t-s)A}
\Psi_2(\bar X(s),\zeta_{s/\e})y_j(s)ds$.
The coordinatewise averaging and spatial-tail calculation used in
\eqref{mod-z} applies to each fixed $y_j$ and gives
\begin{equation}\label{II1-fixed-net}
 \limsup_{\e\to0}h^{-2}(\e)\log
 \PP(\|\mathrm{II}_{1,j}^\e\|_\infty>r/2)=-\infty,
 \qquad 1\leq j\leq N_\kappa.
\end{equation}
On $\{\lambda^\e\in B(y_j,\dl)\}$, the operator bound
$\|\Psi_2\|_{\mathcal L(H)}\leq2K_F'$ yields
\begin{align*}
 \|\mathrm{II}_1^\e-\mathrm{II}_{1,j}^\e\|_\infty
 &\leq\sup_{t\leq1}\int_0^t
 \|e^{(t-s)A}\|\,\|\Psi_2(\bar X(s),\zeta_{s/\e})\|
 |\lambda^\e(s)-y_j(s)|_Hds \leq2K_F'\dl<r/2.
\end{align*}
Consequently,
\[
\PP(\|\mathrm{II}_1^\e\|_\infty>r)
\le\sum_{j=1}^{N_\kappa}\PP(\|\mathrm{II}_{1,j}^\e\|_\infty>r/2)
+\PP\big(\lambda^\e\notin \textstyle\bigcup_{j=1}^{N_\kappa}B(y_j,\dl)\big).
\]
Equation~\eqref{II1-fixed-net} and the LDP upper bound for $\la^\e$ on the complement of $\bigcup_{j}B(y_j,\dl)$ lead to
\[
\limsup_{\e\to0}h^{-2}(\e)\log
\PP(\|\mathrm{II}_1^\e\|_\infty>r)\le-\kappa.
\]  Letting $\kappa \to\infty$ proves
 \bea\ad
 \limsup_{\e\to0}h^{-2}(\e)\log\PP(\|\mathrm{II}_1^\e\|_\infty>r)=-\infty.
\eea
Together with the estimates for
$\mathrm{II}_2^\e$, $\mathrm{II}_3^\e$, and $\mathrm{II}_4^\e$, this
yields \eqref{neg-psi2}.  Equations \eqref{neg-psi1} and
\eqref{neg-psi2}, followed by Gronwall's inequality, prove
\eqref{neg-wdt-eta-hat}.  Combining it with \eqref{eta-wdt-eta}
proves the asserted exponential equivalence. \qed

\section{State-dependent regime switching}\label{app:feedback-auxiliary}
\subsection{Proof of Lemma~\ref{lem:sd-average}}\label{app:lem43}
First suppose that each $c(\cdot,i)$ is continuously differentiable
on $(0,1]$ and that $\max_i|\partial_t c(t,i)|\le K(1+t^{-1})$ there.
Let $\varphi(t)$ be the centered solution of
$-Q(t,\bar y_t)\varphi(t)=c(t)$, $\nu(t)\varphi(t)=0$.
By \eqref{bd-pi}, \lemref{lem:avg-al}, and the Poisson equation,
\beq{sd-varphi-regularity}
\max_i|\varphi(t,i)|\le K_c,
\qquad\max_i|\partial_t\varphi(t,i)|
\le K_c(1+t^{-1}),\quad 0<t\le1.
\eeq
Set $t_R=t\wedge\tau_R^\e$. With the same compensated Poisson random
measure as in \eqref{sd-NU}, define
\begin{equation}\label{sd-Nc}
N_c^\e(t)=\int_\e^{t\vee\e}\int_{\rr}
\bigl[\varphi(s,\al^\e(s-)+\mathfrak h^\e(s,X^\e(s-),\al^\e(s-),z))-\varphi(s,\al^\e(s-))\bigr]
\,\widetilde{\mathfrak p}(ds,dz).
\end{equation}
For $t_R\ge\e$, the same generalized It\^o formula used in
\eqref{Ito-U} gives
\begin{align}
\int_0^{t_R}c(s,\alpha^\e(s))\,ds
&=\int_0^\e c(s,\alpha^\e(s))\,ds
+\e N_c^\e(t_R)
+\e\int_\e^{t_R}\partial_s\varphi(s,\alpha^\e(s))\,ds\notag\\
&\quad-\e[\varphi(t_R,\alpha^\e(t_R))
-\varphi(\e,\alpha^\e(\e))]\notag\\
&\quad+\int_\e^{t_R}\sum_{j=1}^{m_0}
[q_{\alpha^\e(s)j}(s,\Upsilon X^\e(s))
-q_{\alpha^\e(s)j}(s,\bar y_s)]\varphi(s,j)\,ds.
\label{sd-dynkin-average}
\end{align}
For $t_R<\e$, the integral is bounded by $\e\|c\|_\infty$.
The jumps of $N_c^\e$ have norm at most $2K_c$ and
predictable total intensity bounded by $
\kappa_q/\e$.
Before $\tau_R^\e$,
$\|Q(s,\Upsilon X^\e(s))-Q(s,\bar y_s)\|
\le K\|\Upsilon\|\sqe h(\e)R$.
Therefore the nonmartingale terms satisfy
\beq{sd-average-rem}
\e\|c\|_\infty+K_c\e(3+|\log\e|)
+KK_c\sqe h(\e)R\longrightarrow0.
\eeq
 In \eqref{fb-jump-exp}, take
$
\fg(s,i,j)=\indi_{\{s>\e\}}[\varphi(s,j)-\varphi(s,i)]$.
Then $\kappa_{\fg}=2K_c$,
$\kappa_V=4K_c^2\kappa_q/\e$, and
\[
\PP\left(\e\sup_{t\le1}|N_c^\e(t)|>r/2\right)
\le2\exp\left(-\frac{r^2}
{\e(32eK_c^2\kappa_q+8K_c r)}\right),
\]
Since $\e h^2(\e)\to0$, this proves \eqref{sd-average-stopped}
for $C^1$ functions.

For continuous $c$, choose $c^{(m)}(\cdot,i)\in C^1([0,1])$
with $\max_i\|c^{(m)}(\cdot,i)-c(\cdot,i)\|_\infty\to0$. Let
\[
\widetilde c^{(m)}(t,i)=c^{(m)}(t,i)
-\sum_{j=1}^{m_0}\nu_j(t)c^{(m)}(t,j).
\]
Then
\beq{sd-c-approx-error}
\sum_i\nu_i(t)\widetilde c^{(m)}(t,i)=0,
\quad \max_i\|\widetilde c^{(m)}(\cdot,i)-c(\cdot,i)\|_\infty
\le2\max_i\|c^{(m)}(\cdot,i)-c(\cdot,i)\|_\infty\to0.
\eeq
Furthermore, by \lemref{lem:avg-al},
\[
\begin{aligned}
\partial_t\widetilde c^{(m)}(t,i)
&=\partial_t c^{(m)}(t,i)
-\sum_j\nu_j(t)\partial_t c^{(m)}(t,j)
-\sum_j\nu_j'(t)c^{(m)}(t,j),\\
\max_i|\partial_t\widetilde c^{(m)}(t,i)|
&\le K_{m}(1+t^{-1}),\qquad 0<t\le1.
\end{aligned}
\]
Applying the scalar estimate to $\widetilde c^{(m)}$ and then using \eqref{sd-c-approx-error} proves \eqref{sd-average-stopped} for continuous $c$. For any fixed $n$, apply this estimate to the $n$ scalar coordinates of a bounded $H_n$-valued function. Since the dimension is fixed, a union bound proves the projected assertion.

For the last assertion, let $K_\CO$ be a common Lipschitz constant of
$v\mapsto \ell_v$ on $\CO$ and choose a finite $\delta$-net
$v_1,\ldots,v_k$ of $\CO$.  On
$\{\eta^\e\in \CO\cap B(v_j,\delta)\}$,
\beqa{sd-net-error}
\sup_{t\leq1}\big|\int_0^{t\wedge\tau_R^\e}e^{(t-s)A}
 \{\ell_{\eta^\e}(s,\alpha^\e(s))-\ell_{v_j}(s,\alpha^\e(s))\}\,ds
 \big|_H \leq K_\CO\delta.
\eeqa
For fixed $v_j$, set
$\wdt V_n^\e(t)=\int_0^{t\wedge\tau_R^\e}
\Pi_n\ell_{v_j}(s,\alpha^\e(s))\,ds$. Then integration by parts gives
\bea\ad
\int_0^{t\wedge\tau_R^\e}e^{(t-s)A}\Pi_n\ell_{v_j}(s,\alpha^\e(s))\,ds
=(\Phi_n \wdt V_n^\e)(t).
\eea
By \eqref{def-Phi-n},
the uniform norm of this convolution is at most
$(1+\|A_n\|)\|\wdt V_n^\e\|_\infty$.
The projected assertion controls the last term. The remaining
spatial modes satisfy
\beqa{sd-average-tail}\ad
 \sup_{t\leq1}\big|(I-\Pi_n)\int_0^{t\wedge \tau_R^\e} e^{(t-s)A}
 \ell_v(s,\alpha^\e(s))ds\big|_H
 \leq\alpha_{n+1}^{-1}\sup_{s,i,v\in \CO}|\ell_v(s,i)|_H,
\eeqa
because
$\|(I-\Pi_n)e^{rA}\|_{\CL(H)}\leq e^{-\alpha_{n+1}r}$.  First choose
$n$ so that the right-hand side of \eqref{sd-average-tail} is below
$r/3$, then choose $\delta<r/(3K_\CO)$, and finally apply a finite union
bound over $j=1,\ldots,k$.  Equations \eqref{sd-net-error} and
\eqref{sd-average-tail} prove \eqref{sd-random-average}. \qed

\subsection{Proof of \propref{prop:equi-eta-al}}\label{app:prop-equi-eta-al}
For $i=\al^\e(s)\in \CM$, we define
\bea
E_1^\e(s) \ad:=(\sqe h(\e))^{-1}[F(\bar X(s)+\sqe h(\e)\eta^\e(s),i)-F(\bar X(s),i)]-\nabla F(\bar X(s),i) \eta^\e(s) \\
E_2^\e(s)\ad: = \big[\nabla F(\bar X(s),i)-\sum_{j=1}^{m_0} \nu_j(s)\nabla F(\bar X(s),j)\big]\eta^\e(s)\\
E_3^\e(s)\ad:= \indi_{\{s\ge\e\}}\sum_{j=1}^{m_0} \nabla_y q_{ij}(s,\bar y_s)[\Upsilon \eta^\e(s)] U(s,j)-B_\pi(s)\eta^\e(s) \\
E_4^\e(s)\ad := \indi_{\{s\ge\e\}}(\sqe h(\e))^{-1}\sum_{j=1}^{m_0} \big[q_{ij}(s,\Upsilon X^\e(s))-q_{ij}(s,\bar y_s) \\
\aad\qquad -\sqe h(\e) \nabla_y q_{ij}(s,\bar y_s)[\Upsilon \eta^\e(s)]\big]U(s,j).
\eea
Set $\CE_k^\e(t):=\int_0^t e^{(t-s)A} E_k^\e(s)ds$ for $k\leq4$ and $\Dl_\eta^\e:=\eta^\e-\wdh\eta_\al^\e$. By \eqref{decomp-eta}, \eqref{nab-f}, \eqref{wdt-la},
\bea\ad
\Dl_\eta^\e(t)
=\int_0^t e^{(t-s)A}\nabla\bar F_{\rm sd}(s,\bar X(s))
\Dl_\eta^\e(s)\,ds
+\sum_{k=1}^4\CE_k^\e(t)+\CR_2^\e(t).
\eea
 With $K_{\rm sd}:=\sup_{s\in[0,1]} \|\nabla \bar F_{\rm sd}(s,\bar X(s))\|_{\CL(H)}$, Gronwall's inequality yields
\beqa{est-Dl-eta}
\|\Dl_\eta^\e\|_\infty \leq e^{K_{\rm sd}}\big(\sum_{k=1}^4 \|\CE_k^\e\|_\infty + \|\CR_2^\e\|_\infty \big).
\eeqa
We now estimate the five terms on the right-hand side of the above inequality.

On the event $\{\tau_R^\e=1\}$, the bounded second derivative of $q_{ij}$ gives
\beq{est-E4}
\|E_4^\e\|_{L^1([0,1];H)}+\|\CE_4^\e\|_\infty \leq K \sqe h(\e)R^2.
\eeq
For every compact set $\CO\subset C^0([0,1];H)$, assumption (H1) implies
\bea\ad
\sup_{v\in \CO} \int_0^{1} \sup_{i\in\CM} \Big|\frac{F(\bar X(s)+\sqe h(\e) v(s),i)-F(\bar X(s),i)}{\sqe h(\e)} -\nabla F(\bar X(s),i)v(s)\Big|_H ds \to 0.
\eea
Hence for every $r>0$ and sufficiently small $\e>0$, $\eta^\e\in \CO$ implies $\|\CE_1^\e\|_\infty \leq r$.

For a fixed $v\in C^0([0,1];H)$, let
\bea
\ell_{2,v}(s,i)\ad:= \big[\nabla F(\bar X(s),i)-\sum_{j} \nu_j(s) \nabla F(\bar X(s),j)\big]v(s) \\
\ell_{3,v}(s,i)\ad :=\sum_{j=1}^{m_0} \nabla_{y}q_{ij}(s,\bar y_s)[\Upsilon v(s)] U(s,j)-B_\pi(s) v(s).
\eea
The function $\ell_{2,v}$ is centered by definition. To verify the same property for $\ell_{3,v}$, differentiating $\pi(t,y)Q(t,y)=0$ at $y=\bar y_t$ in the direction $w \in \rr^{d_\Upsilon}$ with respect to $y$ and multiplying by $U$ gives $\nu(t)\nabla_y Q(t,\bar y_t)[w]U(t)+ \nabla_y \pi(t,\bar y_t)[w] Q(t,\bar y_t) U(t) =0 $. Hence, by \eqref{Poi-al}, we have
\bea\ad
\nu(t)\nabla_y Q(t,\bar y_t)[w]U(t) = \nabla_y \pi(t,\bar y_t)[w] b_\al(t)= \sum_{i=1}^{m_0} \nabla_y \pi_i(t,\bar y_t)[w] F(\bar X(t),i)
\eea
where the last equality uses the fact that $\sum_{i=1}^{m_0} \nabla_y \pi_i=0$. Therefore, taking $w=\Upsilon v(s)$ and using the definition of $B_\pi$ in \eqref{def-Bpi} gives $\sum_{i\in \CM} \nu_i(s) \ell_{3,v}(s,i)=0$. On $\{\tau_R^\e=1\}$, the truncation introduces only
$
\|E_3^\e-\ell_{3,\eta^\e}(\cdot,\alpha^\e(\cdot))\|_{L^1([0,1];H)}
\le KR\e.
$
Consequently, \lemref{lem:sd-average} implies that for $k=2,3$, every compact set $\CO\subset C^0([0,1];H)$ and all $R,r>0$,
\beqa{err-1}\ad
\limsup_{\e\to 0} h^{-2}(\e)\log\PP\big(\eta^\e\in \CO, \tau_R^\e=1, \|\CE_k^\e\|_\infty > r\big)=-\infty.
\eeqa

Finally, (H2) and \lemref{lem:conv-exp} imply that
\beqa{err-fR}
\PP(\|\CR_2^\e\|_\infty >r, \tau_R^\e=1)\leq K_1 \exp(-r^2/(K_2K\e R^2)).
\eeqa
We next prove the compact restriction in \eqref{eta-com}. Set
\bea\disp
 b_R^\e(s)=\frac{\indi_{\{s\leq\tau_R^\e\}}}{\sqe h(\e)}
 \Big[&F(X^\e(s),\al^\e(s))-F(\bar X(s),\al^\e(s))\\
 &+\indi_{\{s\geq\e\}}\sum_j
 \big(q_{\al^\e(s)j}(s,\Upsilon X^\e(s))
 -q_{\al^\e(s)j}(s,\bar y_s)\big)U(s,j)\Big].
\eea
Then $\|b_R^\e\|_{L^2([0,1];H)}\leq K_R$ and, on
$\{\tau_R^\e=1\}$,
$
 \eta^\e=R_\al^\e+\CR_{2,R}^\e+\Scv b_R^\e.
$
By \appref{app:quad-lift}, $\Scv$ is compact from $L^2([0,1];H)$ to $C^0([0,1];H)$.

Fix $\beta,\sigma,a>0$ with
$\beta+\sigma<a<(1-\vartheta)/2$ and $p_1\geq2$ with $a-\beta-\sigma>1/p_1$.
Define
\bea\ad
 \Phi_R^\e(s)=\frac{\indi_{\{s\leq\tau_R^\e\}}}{h(\e)}
 [G(X^\e(s),\al^\e(s))-G(\bar X(s),\al^\e(s))]\mathsf{Q}.
\eea
By \eqref{G-diff},
$
 \|e^{uA}\Phi_R^\e(s)\|_{\CL_2(H)}^2
 \leq K\e R^2u^{-\vartheta},\, 0<u\leq1.
$
The factorization formula gives
\begin{align*}
 \CR_{2,R}^\e(t)
 &=\frac{\sin(\pi a)}{\pi}
 \int_0^t(t-s)^{a-1}e^{(t-s)A}Y_a^\e(s)\,ds,\\
 Y_a^\e(t)
 &=\int_0^t(t-s)^{-a}e^{(t-s)A}\Phi_R^\e(s)\,dW_0(s).
\end{align*}
Thus, for $p\geq p_1$, the Hilbert-valued BDG inequality yields
\bea
 \big(\EE|Y_a^\e(t)|_H^p\big)^{2/p}
 &\leq Kp\e R^2\int_0^t u^{-2a-\vartheta}du
 \leq Kp\e R^2,\\
 \big(\EE\|Y_a^\e\|_{L^p([0,1];H)}^p\big)^{1/p}
 &\leq K\sqrt{p\e}\,R.
\eea
By $\|(-A)^\beta e^{uA}\|\leq Ku^{-\beta}$ and
$\|(-A)^\beta(e^{vA}-I)e^{uA}\|\leq Kv^{\sigma} u^{-\beta-\sigma}$,
H\"older's inequality establishes that
uniformly for $p\geq p_1$,
$
 \|\CR_{2,R}^\e\|_{C^{\sigma}([0,1];\mathcal H_\beta)}
 \leq K\|Y_a^\e\|_{L^p([0,1];H)}$.
Optimizing Markov's inequality over $p\geq p_1$ gives that for fixed
$R,r'>0$ and all sufficiently small $\e$,
\begin{equation}\label{sd-R2-compact-tail}
 \PP\big(\|\CR_{2,R}^\e\|_{C^{\sigma}([0,1];\mathcal H_\beta)}>r'\big)
 \leq 2\exp\{-r'^2/(K\e R^2)\}.
\end{equation}
By the compact embedding $\mathcal H_\beta\hookrightarrow H$,
the unit ball of $C^{\sigma}([0,1];\mathcal H_\beta)$ is relatively
compact in $C^0([0,1];H)$. For $R,\kappa>0$, choose a compact
$\SK_\kappa\subset C^0([0,1];H)$ with
$\limsup_{\e\to0}h^{-2}(\e)\log\PP(R_\al^\e\notin \SK_\kappa)\leq-\kappa$
using \propref{prop:sd-driving}, and set
\[
 \CO_{R,\kappa}=\SK_\kappa
 +\overline{\{x:\|x\|_{C^{\sigma}([0,1];\mathcal H_\beta)}\leq1\}}
 +\overline{\{\Scv b:\|b\|_{L^2([0,1];H)}\leq K_R\}},
\]
with closures in $C^0([0,1];H)$. The three summands are compact, and
\[
 \PP(\eta^\e\notin\CO_{R,\kappa},\tau_R^\e=1)
 \leq\PP(R_\al^\e\notin \SK_\kappa)
 +\PP(\|\CR_{2,R}^\e\|_{C^{\sigma}([0,1];\mathcal H_\beta)}>1).
\]
Thus \eqref{sd-R2-compact-tail} gives
\beq{eta-com}
\limsup_{\e\to0}h^{-2}(\e)
\log\PP(\eta^\e\notin\CO_{R,\kappa},\tau_R^\e=1)\leq-\kappa.
\eeq

Fix $r>0$ and define $r_0:=r\, e^{-K_{\rm sd}}/5$. From \eqref{est-Dl-eta}, we have
\bea
\PP(\|\Dl_\eta^\e\|_\infty > r) \ad \leq \PP(\tau_R^\e<1)+\PP(\eta^\e \notin \CO_{R,\kappa},\tau_R^\e=1) +\PP(\|\CR_2^\e\|_\infty >r_0, \tau_R^\e=1) \\
\aad \quad + \sum_{k=1}^4 \PP(\|\CE_k^\e\|_\infty >r_0, \eta^\e \in \CO_{R,\kappa}, \tau_R^\e=1).
\eea
Applying \eqref{est-E4}--\eqref{err-fR} and \eqref{lim-tau} and letting $\kappa \to \infty$ in \eqref{eta-com} and $R\to \infty$, we obtain \eqref{equi-eta-al}. \qed

\section{Stochastic reaction--diffusion verification}\label{app:rd-verification}
\begin{proof}[Proof of Proposition~\ref{prop:rd-realization}]
The Nemytskii map $F$ is globally Lipschitz on $L^2$, with G\^ateaux
derivative $\nabla F(x,z)v=f_r(\cdot,x(\cdot),z)v$.
If $\SK\subset L^2(\SO)$ is compact, then
$\{|v|^2:v\in\SK\}$ is uniformly integrable.
Splitting $\SO$ into $\{|v|\le R\}$ and $\{|v|>R\}$ and using
the bounds on $f_r$ and $f_{rr}$ show that
$(x,v)\mapsto\nabla F(x,z)v$ is uniformly continuous on every
compact subset of $L^2(\SO)\times L^2(\SO)$, uniformly in
$z\in\SE$. This verifies (H1).

It remains to prove the smoothing bounds.  For
$b\in L^2(\SO)$, let $M_b$ be the multiplication by $b$. It is self-adjoint and $\|M_b\|_{\CL(L^\infty,L^2)}+\|M_b\|_{\CL(L^2,L^1)}\leq2|b|_{L^2}$. Conditions \eqref{kap-zeta} and \eqref{kap-zeta2} are (2.1)--(2.3) in Cerrai \cite[Hypothesis~1]{Cer09}, with $p_0$ in place of $\rho$ and $\rho=\infty$ if $d=1$. Hence (3.3)--(3.5) in the proof of \cite[Lemma~3.1]{Cer09}, with $J=M_b$, give
$
 \|e^{rA}M_b\mathsf{Q}\|_{\CL_2}^2\leq Kr^{-\vartheta}|b|_{L^2}^2,\, 0<r\leq1.$
Taking $b=g(\cdot,0,z)$ and
$b=g(\cdot,x,z)-g(\cdot,y,z)$ proves
\eqref{G0} and \eqref{G-diff}.  Boundedness of $g$ gives $\|G(x,z)\mathsf{Q}\|_{\CL(H)}\leq K_g\|\mathsf{Q}\|_{\CL(H)}$, which is \eqref{G-op}. The bound \eqref{abstract-Gtail} and the Lipschitz continuity of the entries in \remref{rem:G-bounds} follow from (H2).
\end{proof}

\section*{Acknowledgments}
The first author would like to thank Professor Arnaud Guillin
for sharing his early paper \cite{Gui00}. Part of the work was completed while H. Qian was a graduate student in the Department of Mathematics at the University of Connecticut. The first author is the corresponding author.

\section*{Funding}
The research of the first author was supported by a Wayne State University start-up fund.
The research of the second author was supported in part by the National Science Foundation under grant DMS-2204240.
\section*{Statement on AI use}
The authors learned from an interaction with ChatGPT (GPT-5.6 Sol, OpenAI) that a Galerkin approximation can be used in Section 3.3.2. ChatGPT also drew Figures 1 and 2 and compiled the list of notation. Apart from this, the problem, the main results and the method of proof are due to the authors, who also wrote the first complete version of the proofs. Claude Fable 5.2 (Anthropic) was used to check the proofs. It pointed out gaps, inaccurate references to the literature, clashes of notation and errors in constants, and it suggested corrections. Both were also used to reorganize the exposition, to polish and proofread the language, and to check the bibliographic data. The authors have verified all statements and proofs and take full responsibility for the content of the paper.

\end{document}